\documentclass[12pt]{amsart}

\usepackage{mathtools}
\usepackage{amsmath,amssymb,amsthm}
\usepackage{enumerate}
\usepackage[noadjust]{cite}

\usepackage[margin=1in]{geometry}

\usepackage{hyperref}
\usepackage[spanish, english]{babel}

\usepackage[normalem]{ulem}

\makeatletter % for internal macros with @
\renewcommand\p@enumii{}
\makeatother

\usepackage[ansinew]{inputenc}
\usepackage{here}

\newcommand{\A}{\mathcal{A}}
\newcommand{\BB}{\mathcal{B}}
\newcommand{\CC}{\mathcal{C}}
\newcommand{\FF}{\mathcal{F}}

\newcommand{\RR}{\mathcal{R}}
\newcommand{\SSS}{\mathcal{S}}
\newcommand{\UU}{\mathcal{U}}
\newcommand{\XX}{\mathcal{X}}
\newcommand{\YY}{\mathcal{Y}}

\newcommand{\CCone}{\CC^{T_3}_{\dim 1}}
\newcommand{\CCCsep}{\CC^{T_3,\, \text{sep}}_{\dim 1}}
\newcommand{\CCCfdi}{\CC^{T_2,\, \text{fdi}}_{\dim 1}}

\newcommand{\al}{\langle}
\newcommand{\ar}{\rangle}

\newcommand{\R}{\mathbb{R}}

\newcommand{\tauc}{\tau_s}
\newcommand{\mytau}{\hat{\tau}}
\newcommand{\mytauu}{\tilde{\tau}}
\newcommand{\taulex}{\tau_{lex}}

\newcommand{\taualex}{\tau_{Alex}}

\DeclareMathOperator{\taulim}{\tau{\operatorname{-}}lim}

\DeclareMathOperator{\elim}{{\it e}{\operatorname{-}}lim}
\DeclareMathOperator{\mytaulim}{\mytau{\operatorname{-}}lim}

\newcommand{\precc}{\preceq}

\newcommand{\econtinuous}{$e$-continuous }
\newcommand{\econtinuousp}{$e$-continuous}
\newcommand{\ehomeomorphism}{$e$-homeomorphism }
\newcommand{\ehomeomorphismp}{$e$-homeomorphism}

\newcommand{\Xopc}{X^c}
\newcommand{\tauopc}{\tau^c}

\newcommand{\Xcp}{X^c}
\newcommand{\taucp}{\tau^c}

\newcommand{\exR}{R_{\pm\infty}}
\newcommand{\num}{\textbf{n}}

\newcommand{\ssim}{\sim_\tau}

\newcommand{\lex}{<_{lex}}

\newcommand{\Xone}{X_n^{(1)}}
\newcommand{\Xtwo}{X_n^{(2)}}

\newcommand{\Xinf}{\XX_{\text{open}}}
\newcommand{\Xsing}{\XX_{\text{sgl}}}
\newcommand{\Asing}{F_{\text{sgl}}}
\newcommand{\Xisol}{A_{\text{isol}}}
\newcommand{\XXisol}{\XX_{\text{isol}}}
\newcommand{\Xgood}{X'}

\newcommand{\E}{E}

\newcommand{\cond}{(\dag)}

\theoremstyle{definition}
\newtheorem{definition}{Definition}[section]
\newtheorem{remark}[definition]{Remark}
\newtheorem{remark*}{Remark}
\newtheorem{example}[definition]{Example}
\newtheorem{exmpl}[definition]{Example}

\theoremstyle{theorem}
\newtheorem{theorem}[definition]{Theorem}
\newtheorem*{theorem*}{Theorem}
\newtheorem{introtheorem}{Theorem}
\newtheorem{lemma}[definition]{Lemma}
\newtheorem{proposition}[definition]{Proposition}
\newtheorem{corollary}[definition]{Corollary}
\newtheorem{claim}{Claim}[definition]

\newtheorem{question}[definition]{Question}

\numberwithin{equation}{section}

\begin{document}
\title[One-dimensional definable topologies in o-minimal structures]{One-dimensional definable topological spaces in o-minimal structures}
\author[P. And\'ujar Guerrero]{Pablo And\'ujar Guerrero}
%\thanks{}
\address{School of Mathematics, University of Leeds, Leeds, LS2 9JT, UK}
\email{pa377@cantab.net}

\author[M. E. M. Thomas]{Margaret E. M. Thomas}
%\thanks{The fourth author is supported by NSF grant DMS-2154328.}
\address{Department of Mathematics, Purdue University, 150 N. University Street, West Lafayette, IN 47907-2067, USA}
\email{memthomas@purdue.edu}

\begin{abstract}
We study the properties of topological spaces $(X,\tau)$, where $X$ is a definable set in an o-minimal structure and the topology $\tau$ on $X$ has a basis that is (uniformly) definable. Examples of such spaces include the canonical euclidean topology on definable sets, definable order topologies, definable quotient spaces and definable metric spaces. We use o-minimality to undertake their study in topological terms, focussing here in particular on spaces of dimension one. We present several results, given in terms of piecewise decompositions and existence of definable embeddings and homeomorphisms, for various classes of spaces that are described in terms of classical separation axioms and definable analogues of properties such as separability, compactness and metrizability. For example, we prove that all Hausdorff one-dimensional definable topologies are piecewise the euclidean, discrete, or upper or lower limit topology; we give a characterization of all one-dimensional, regular, Hausdorff definable topologies in terms of spaces that have a lexicographic ordering or a topology generalizing the Alexandrov double of the euclidean topology; and we show that, if the underlying structure expands an ordered field, then any one-dimensional Hausdorff definable topology that is piecewise euclidean is definably homeomorphic to a euclidean space. As applications of these results, we prove definable versions of several open conjectures from set-theoretic topology, due to Gruenhage and Fremlin, on the existence of a 3-element basis for regular, Hausdorff topologies and on the nature of perfectly normal, compact, Hausdorff spaces; we obtain universality results for some classes of Hausdorff and regular topologies; and we characterize when certain metrizable definable topologies admit a definable metric.

\end{abstract}



\maketitle

\section{Introduction}
\renewcommand{\theintrotheorem}{\Alph{introtheorem}}

A topological space $(X,\tau)$ is definable in a first order structure $\RR=(R,\ldots)$ if $X\subseteq R^n$, for some natural number $n$, and the topology $\tau$ has a basis that is (uniformly) definable in $\RR$. In other words, the topology $\tau$ is explicitly definable in the sense of Flum and Ziegler~\cite{flum_ziegler_80}. This definition generalizes the notion of first order topological structure in~\cite{pillay87}, which addresses the case where $X=R$.
A basic example of a topological space which is definable in any structure is the discrete topology on a definable set. Other examples of definable topological spaces arise from valued fields and the field of complex numbers with a predicate for the reals~\cite{pillay87}.

In this paper, we begin a detailed study of the general theory of topological spaces that are definable in o-minimal structures. A structure is \emph{o-minimal} if it is an expansion of a linear order $(R,<)$ satisfying that every unary definable set is a finite union of points and intervals with endpoints in $R\cup\{-\infty, +\infty\}$. There is a canonical topology in this setting given by the order topology on $R$ and, for every $n>1$, the induced product topology on $R^n$. In the present paper we refer to these definable topologies collectively as the \emph{euclidean topology}. Much of the research in o-minimality as a tame topological setting has centered on the analysis of the euclidean topology, as this has long been recognized as a suitable setting for the ``topologie mod\'er\'ee" of Grothendieck (i.e. a setting that avoids the pathologies of general set-theoretic topology) \cite{gro84}. Nevertheless, some examples of other o-minimal definable topological spaces have been explored in the literature, including definable manifold spaces~\cite[Chapter 10, Section 1]{dries98}, which encompass definable groups~\cite{pillay88}; definable metric spaces~\cite{walsberg15}; and definable orders (which generate definable order topological spaces)~\cite{ram13}. We direct the reader to \cite{dries98} for the requisite background on the theory of o-minimal structures, which will be used extensively throughout this paper.

Our perspective is to consider definable topological spaces $(X,\tau)$ in o-minimal structures from a general standpoint. Here we focus predominantly on the case where $\dim X =1$ but, even in this setting,  our perspective brings to the o-minimal setting topological spaces that exhibit a wide variety of topological properties (see Appendix~\ref{section:examples} for a catalogue of examples), including classical spaces 
such as the Sorgenfrey Line, the Split Interval and the Alexandrov Double Circle (Examples~\ref{example:taur_taul}, \ref{example:split_interval} and~\ref{example:alex_double_circle} respectively). These are common counterexamples in topology, displaying properties that the more ``well-behaved" spaces lack. Nevertheless, we argue through the present work that this setting retains some form of tameness. 
Specifically, the axiom of o-minimality implies that the structure of one-dimensional definable topological spaces is rather restrictive, in particular when compared to spaces of higher dimensions.

The main contributions of this paper are a series of decomposition and embedding theorems for one-dimensional o-minimal definable topological spaces satisfying certain classical separation axioms. These results can be understood as classifying such spaces, in the sense that they can be described in terms of a few classical examples. Much of our approach is motivated by partition problems in set-theoretic topology, which seek to understand topological spaces in similar terms under various axioms of set theory (see for example~\cite{todor89}). In order to prove our results and build an extensive theory of o-minimal definable topology, we also introduce to our setting suitable definable analogues of classical topological properties, including separability, metrizability and compactness, and investigate them in depth. 

We now describe the main results of this paper. We begin with the following decomposition results for $T_1$ and Hausdorff spaces. Note that, in the following statement, the right half-open interval topology (also known as the lower limit topology) is the definable analogue of the topology of the classical Sorgenfrey Line. 
\begin{introtheorem} [Corollary~\ref{cor:general_them_main}, Corollary~\ref{cor 1.5}] \label{introthm:T_1_T_2_spaces}
    Let $\RR$ be an o-minimal structure and let $(X,\tau)$ be a definable topological space in $\RR$.
\begin{enumerate}[I.]
    \item If $(X,\tau)$ is infinite and $T_1$, then it has a subspace that is definably homeomorphic to an interval with either the euclidean, right half-open interval, left half-open interval, or discrete topology. \label{itm:introthm:3EBC}
    \item If $(X, \tau)$ is Hausdorff and $\dim(X) \leq 1$, then there exists a finite definable partition $\XX$ of $X$ such that, for every $Y\in\XX$, $(Y, \tau)$ is definably homeomorphic to a point or an interval with either the euclidean, discrete, right half-open interval or left half-open interval topology. \label{itm:introthm:T_2decomp}
\end{enumerate}
\end{introtheorem}
Theorem~\ref{introthm:T_1_T_2_spaces}.\ref{itm:introthm:3EBC} can be seen as an o-minimal version of an open conjecture of set-theoretic topology, due to Gruenhage, known as the 3-element basis conjecture (see \cite{gru90}, \cite{gm07}). This conjecture states that it is consistent with ZFC that every uncountable, first-countable, regular, Hausdorff topological space contains a subspace homeomorphic to either a fixed subset of the reals with the euclidean topology, or a fixed subset of the reals with the Sorgenfrey topology, or contains an uncountable discrete subspace. Some cases of this conjecture are known that are conditional upon various axioms of set theory (see Subsection~\ref{subsection: 3-el_basis_conj} for further discussion). From Theorem~\ref{introthm:T_1_T_2_spaces}.\ref{itm:introthm:3EBC} we derive that the conclusion of the conjecture holds unconditionally for all infinite $T_1$ topological spaces definable in any o-minimal expansion of $(\R,<)$, and that a definable generalization holds in any o-minimal structure (see Subsection~\ref{subsection: 3-el_basis_conj}).

We also use Theorem~\ref{introthm:T_1_T_2_spaces}.\ref{itm:introthm:3EBC} to show that no space homeomorphic to the Cantor space $2^{\omega}$ is definable in any o-minimal structure (Corollary~\ref{cor:no_cantor_set}), supporting the thesis that o-minimality provides a tame topological setting in the general definable topological context.

Further to the above, we also investigate o-minimal one-dimensional definable topological spaces that are regular as well as Hausdorff, improving the conclusion of Theorem~\ref{introthm:T_1_T_2_spaces}.\ref{itm:introthm:T_2decomp} for these spaces by showing that they can be definably decomposed in terms of suitable generalizations of the Split Interval and the Alexandrov Double Circle, the latter with a topology which we call the `Alexandrov topology' (see Examples~\ref{example:n-split} and~\ref{example_n_line}). The precise statement in the case of spaces $(X,\tau)$, with $X$ a unary definable set, is as follows.
\begin{introtheorem}[Theorem \ref{them_ADC_or_DOTS}]\label{introthm:ADC_or_DOTS}
Let $\RR$ be an o-minimal structure and let $(X,\tau)$, $X\subseteq R$, be a regular and Hausdorff definable topological space in $\RR$. Then there exist disjoint definable open sets $Y, Z \subseteq X$ with \mbox{$X\setminus (Y\cup Z)$} finite, and $n_Y, n_Z>0$, such that $(Y,\tau)$ embeds definably into $R \times \{0,\ldots, n_Y\}$ endowed with the lexicographic order topology, and $(Z,\tau)$ embeds definably into $R \times \{0,\ldots, n_Z\}$ endowed with the Alexandrov topology. 
\end{introtheorem} 
These decomposition results, Theorems~\ref{introthm:T_1_T_2_spaces} and~\ref{introthm:ADC_or_DOTS}, allow us to address universality questions in our setting. Our motivation is classical universality literature in Banach space theory (e.g. the Banach--Mazur theorem). 
We prove that a number of classes of spaces admit an `almost definably universal' space, in a sense that we make precise (Definition~\ref{dfn:(almost)universal}). Such classes include the class of euclidean spaces of dimension at most $n$, the class of one-dimensional regular Hausdorff definable topological spaces, and the class of Hausdorff definable topological spaces $(X,\tau)$, with $X$ a unary definable set, that satisfy the `frontier dimension inequality' (Definition~\ref{dfn:fdi}). 
We moreover prove some negative results about the existence of definably universal spaces of certain kinds, in particular that there does not exist a one-dimensional, $T_1$ definable topological space that is almost definably universal for the class of one-dimensional regular Hausdorff definable topological spaces.

We also use Theorem~\ref{introthm:ADC_or_DOTS} to study definable Hausdorff compactifications of o-minimal one-dimensional definable topological spaces. Our main result in this respect (Theorem~\ref{them_compactification}) is that such a Hausdorff space $(X,\tau)$, with $X$ a unary definable set, admits a definable embedding into a one-dimensional Hausdorff definably compact space if and only if $(X,\tau)$ is regular.

A main line of research concerning various classes of definable topological spaces has been the study of affineness (the property of being definably homeomorphic to a space with the euclidean topology). In the setting of o-minimal expansions of ordered fields, van den Dries showed that any definable manifold space is affine if and only if it is regular~\cite{dries98}, while Walsberg showed that a definable metric space is affine if and only if it does not contain an infinite definable discrete subspace. Using Theorem~\ref{introthm:T_1_T_2_spaces}.\ref{itm:introthm:T_2decomp}, as well as our work on existence of definable compactifications, we prove, in the same setting, the following characterization of affineness for Hausdorff one-dimensional definable topological spaces.
\begin{introtheorem}[Theorem \ref{them_main_2}] \label{introthm:affine}
Let $\RR$ be an o-minimal expansion of an ordered field, and let $(X, \tau)$ be a Hausdorff topological space definable in $\RR$ with $\dim(X) \leq 1$. Then $(X,\tau)$ is affine if and only if it does not contain a subspace that is definably homeomorphic to an interval with either the discrete or the right half-open interval topology.
\end{introtheorem}
As an application of this theorem, in combination with Theorem~\ref{introthm:ADC_or_DOTS}, we also address in our setting a further set-theoretic question that is closely related to the aforementioned 3-element basis conjecture, namely the possible nature of non-metrizable, perfectly normal, compact, Hausdorff spaces in ZFC. Fremlin asked whether or not it is consistent with ZFC that every perfectly normal, compact, Hausdorff space admits a continuous, at most $2$-to-$1$ map onto a metric space (see~\cite{gm07}), which in turn led Gruenhage to ask if every such space is either metrizable or contains a copy of $A \times \{0,1\}$ equipped with the lexicographic order topology, for some uncountable $A \subseteq [0,1]$ (see~\cite{gru88}). 
It is indicated in~\cite{gru90} and~\cite{gm07} that positive answers to these questions follow from the 3-element basis conjecture (see also Subsection~\ref{subsection: Fremlin}). We consider both of these questions in our setting and show in particular that positive answers to both questions hold, in a definable sense, in any o-minimal expansion of $(\R, +, \cdot, <)$. More specifically, we show that, in such a structure, Fremlin's question has a positive answer for any regular, Hausdorff one-dimensional definable topological space that is either perfectly normal or separable (and, in particular, we do not require that such a space be compact) (Corollary~\ref{cor_Fremlin_2-1}). We also show that Gruenhage's question has a positive answer for any one-dimensional, perfectly normal, compact, Hausdorff topological space definable in such a structure, where we can strengthen the `metrizable' conclusion to `being affine', and in the other case take $A$ to be a subinterval of $[0,1]$ (Corollary~\ref{cor:Fremlin-Q2.2}).

We also use Theorem~\ref{introthm:affine}, this time together with \ref{introthm:T_1_T_2_spaces}.\ref{itm:introthm:T_2decomp}, to investigate a notion of `definable metrizability' extracted from the work of Walsberg~\cite{walsberg15}. Our main result (Theorem~\ref{thm:metrizability}) is that definable metrizability is equivalent to metrizability for one-dimensional topological spaces definable in certain o-minimal expansions of ordered fields (including o-minimal expansions of the field of reals).

Unsurprisingly, studying o-minimal definable topological spaces of higher dimensions is a lot less straightforward than studying those of dimension one. Key results that we present here fail to generalize to higher dimensions. In particular, Theorem~\ref{introthm:affine} does not, which we illustrate by providing an example of a two-dimensional Hausdorff definable topological space that is not affine, but which also does not contain a subspace that is definably homeomorphic to an interval with either the discrete or the right half-open interval topology (Example~\ref{example:fdi_Hausdorff_not_regular}). 
Moreover, we provide an example of a two-dimensional Hausdorff definable topological space that fails to be affine, yet all of its one-dimensional subspaces are affine (Example~\ref{example_line-wise_euclidean_not_euclidean}). Moreover, for regular Hausdorff definable topological spaces, it also does not follow in general that admitting a finite definable partition into euclidean subspaces implies being affine (see Example~\ref{example_space_cell-wise_euclidean_not_metrizable}), although these properties are equivalent for one-dimensional Hausdorff spaces. 

The outline of this paper is as follows. In Section~\ref{section:definitions}, we include many of the necessary definitions. In Section~\ref{section: metric spaces}, we introduce definable metric spaces, studied in~\cite{walsberg15}, to our setting, and discuss the properties of definable metrizability and definable separability in further detail. Section~\ref{section: prel. results} contains preliminary results, both for spaces of all dimensions and for spaces of dimension one in particular. In Section~\ref{section:T1_T2_spaces}, we focus on $T_1$ and Hausdorff spaces, in particular proving Theorem~\ref{introthm:T_1_T_2_spaces} and deducing, from Theorem~\ref{introthm:T_1_T_2_spaces}.\ref{itm:introthm:3EBC}, both the Gruenhage 3-element basis conjecture in our setting and the fact that the Cantor space is not definable. In Section~\ref{section: universal spaces}, we consider regular Hausdorff spaces, and in particular prove Theorem~\ref{introthm:ADC_or_DOTS}. In Section~\ref{section:universal_spaces_2}, we answer various universality questions, some as applications of Theorems~\ref{introthm:T_1_T_2_spaces} and~\ref{introthm:ADC_or_DOTS}. In Section~\ref{section:compactifications}, we prove, as a consequence of Theorem~\ref{introthm:ADC_or_DOTS}, that all regular Hausdorff definable topologies in the line can be Hausdorff compactified in a definable sense. In Section~\ref{section: affine}, we work in an o-minimal expansion of an ordered field, and, using results of previous sections, in particular Theorem~\ref{introthm:T_1_T_2_spaces}.\ref{itm:introthm:T_2decomp} and results on definable compactification, we prove Theorem~\ref{introthm:affine}. We also then prove various statements that address Fremlin's and Gruenhage's questions on perfectly normal, compact, Hausdorff spaces as applications of Theorems~\ref{introthm:ADC_or_DOTS} and~\ref{introthm:affine}.
In Section~\ref{section:metrizability} we prove, as an application of Theorems~\ref{introthm:T_1_T_2_spaces}.\ref{itm:introthm:T_2decomp} and~\ref{introthm:affine}, a theorem implying that, in an o-minimal expansion of the field of reals, any one-dimensional definable topology that is metrizable (with respect to the structure) also admits a definable metric. Finally, Appendix~\ref{section:examples} provides a catalogue of relevant examples, both key examples considered throughout the paper, as well as a number of further examples illustrating the necessity of the hypotheses of many of the results presented here. 

As we were completing this paper, we became aware that Peterzil and Rosel were working independently on similar questions. Section~\ref{sec:Peterzil_Rosel} is a note addressing their paper~\cite{pet_rosel_18}, in which we describe how their main result relates to some of ours, and answer some of their open questions. 

\section*{Acknowledgements}

It is a great pleasure to thank Erik Walsberg for many informative discussions and motivating insights that helped to shape the direction of this paper, as well as for feedback on early drafts. Thanks also to Matthias Aschenbrenner for his considered comments on an earlier version of this work. Thanks as well to Ilijas Farah for pointing us towards some additional background on the set-theoretic questions that we consider in our work (especially for alerting us to the reference~\cite{todorfarah95}) and to Niels J. Diepeveen for a useful discussion via StackExchange (see \cite{diep17}). 

The support of the following while the work on this paper was carried out is gratefully acknowledged. The first author is supported under the UK Engineering and Physical Sciences Research Council (EPSRC) fellowship EP/V003291/1, and was previously supported under the Canada Natural Sciences and Engineering Research Council (NSERC) Discovery Grant RGPIN-06555-2018 and by the Fields Institute for Research in Mathematical Sciences, Toronto, Canada (during the Thematic Programs on \textquotedblleft Trends in Pure and Applied Model Theory'' and \textquotedblleft Tame Geometry, Transseries and Applications to Analysis and Geometry''). The second author is supported by NSF grant DMS-2154328, and was previously supported by the Ontario Baden--W\"urttemberg Foundation, and under the Canada Natural Sciences and Engineering Research Council (NSERC) Discovery Grant RGPIN 261961. Both authors were supported by German Research Council (DFG) Grant TH 1781/2-1; the Zukunftskolleg, Universit\"at Konstanz; and the Fields Institute for Research in Mathematical Sciences, Toronto, Canada (during the Thematic Program on \textquotedblleft Unlikely Intersections, Heights and Efficient Congruencing").

\section{Definitions}\label{section:definitions}

We begin by laying out a number of conventions that we will use throughout the paper. Since our goal is to study various topological spaces and their properties from the perspective of definability in o-minimal structures, and we wish to avoid any ambiguity between conventional concepts and those which we will introduce here, we are careful to make our terminology explicit, even in the case of certain notions that can often be taken as read.

Throughout this paper, $\RR=(R,<,\ldots)$ denotes an o-minimal expansion of a dense linear order without endpoints, possibly with extra assumptions that we make explicit in context. Unless stated otherwise, by \textquotedblleft definable" we will mean \textquotedblleft definable in $\RR$, possibly with parameters from $R$". 
Throughout, $n$ and $m$ denote natural numbers. 
By the euclidean topology on $R^n$, we refer to the canonical topology in an o-minimal structure, which is given by the order topology when $n=1$ and by its induced product topology when $n>1$.
We let $\exR = R \cup \{+\infty, -\infty\}$, and extend the euclidean topology to $\exR^n$ in the natural way. Without reference to a particular topology, any topological notion is to be understood with respect to the euclidean topology.
Unless stated otherwise, by interval we mean an open interval with endpoints in $\exR$. We fix infinitely many parameters $0<1<2<\ldots$ in $R$ in such a way that it will be clear from context when these numerals denote elements of $R$ and when they are just natural numbers. At times, we assume that $\RR$ expands an ordered group $(R,0,+,<)$ or field $(R,0,1,+,\cdot,<)$, in which case these parameters have their natural interpretations. Throughout, let $\pi: R^{n}\rightarrow R$ denote the projection to the first coordinate. We abuse terminology as follows: we say that a relation $\Phi(x,y)\subseteq R^n\times \exR^m$ is definable if its restriction to $R^n\times R^m$ is definable and the family of fibers $\{x\in R^n : \Phi(x,y)\}$ for $y\in \exR^m\setminus R^m$ is definable. Note that in this sense any definable partial function $R\rightarrow \exR$ satisfies o-minimal monotonicity. 

We now define the central object of study in this paper.

\begin{definition}\label{dfn:dtopology}
A \emph{definable topological space} is a tuple $(X,\tau)$, where $X\subseteq R^n$ is a definable set and $\tau$ is a topology on $X$ such that there exists a definable family of sets $\BB_\tau$ that is a basis for $\tau$. We call $\BB_\tau$ a \emph{definable basis for $\tau$} and say that the topology $\tau$ is definable. 
\end{definition}

Clearly there is some redundancy in this definition, as the definability of the basis $\BB_\tau$ implies the definability of the set $X$.

The following are some basic facts about definable topological spaces which are true regardless of the axiom of o-minimality or even the fact that $\RR$ expands a linear order. Familiarity with them will be assumed throughout the paper. 

\begin{proposition}\label{prop_basic_facts_dfn_top_spaces}
Let $(X,\tau)$ and $(Y,\mu)$ be definable topological spaces. 
\begin{enumerate}[(a)]
\item If $\BB_\tau$ is a definable basis for $\tau$, then the family $\BB_\tau(x)=\{A\in\BB_\tau : x\in A\}$ is a basis of open neighbourhoods of $x$ that is definable uniformly on $x\in X$. 
\item Let $Z\subseteq X$ be a definable set. The closure $cl_\tau Z$, interior $int_\tau Z$ and frontier $\partial_\tau Z$ $(:=cl_\tau Z \setminus Z)$ of $Z$ in $(X,\tau)$ are also definable. 
\item Let $f:(X,\tau) \rightarrow (Y,\mu)$ be a definable function. The set of points where $f$ is continuous is definable.
\item If $Z\subseteq X$ is a definable set, then the subspace $(Z,\tau|_Z)$ is a definable topological space. 
\item The product space $(X\times Y, \tau \times \mu)$ is a definable topological space. 
\end{enumerate}
\end{proposition}

When $(X,\tau)$ is a definable topological space and $Y\subseteq X$, we abuse notation by writing $(Y,\tau)$ to mean $(Y,\tau|_Y)$.

Given a definable set $X$, we denote the euclidean and discrete topologies on $X$ by $\tau_e$ and $\tauc$, respectively, in such a way that the notation remains unambiguous. We generally write the letter $e$ in place of $\tau_e$ when used as a prefix or subscript, for example as in $cl_e$ or $e$-neighbourhood, and adopt analogous conventions in writing the letter $s$ in place of $\tauc$. Moreover, given a definable function $f:X\rightarrow Y$, we say that $f$ is \econtinuous (respectively, an \ehomeomorphismp) if, as a map $(X,\tau_e)\rightarrow (Y,\tau_e)$,  $f$ is continuous (respectively, a homeomorphism). 

Recall that a topological space $X$ is $T_1$ if every singleton is closed, $T_2$ if it is Hausdorff and $T_3$ if it is Hausdorff and regular, where regular means that any point $x$ $\in X$ and closed set $C \subseteq X$ with $C \not\ni x$ are separated by neighbourhoods, i.e. there exist disjoint open sets $U, V$ $\subseteq X$ with $x\in U$ and $C\subseteq V$. We approach the study of definable topologies in terms of these three separation axioms. 

We now introduce two definable topologies that are highly relevant to this paper and which are immediate generalizations of classical topologies definable in \mbox{$(\R,<)$}. Since there will be no ambiguity, we use some standard terminology to refer to them as understood in our setting.    

The right half-open interval topology (Appendix \ref{section:examples}, Example~\ref{example:taur_taul}), also called the lower limit topology, on $R$, denoted $\tau_r$, is the topology with definable basis
\[
[x,y)\,\text{ for } x, y\in R,\, x<y.
\]
We reserve the name \textquotedblleft Sorgenfrey Line" to refer to the classical space with this topology, namely $(R,\tau_r)$ when $\RR$ expands $(\R,<)$.
Similarly, the left half-open interval topology (or upper limit topology) on $R$, denoted $\tau_l$, is the topology with definable basis 
\[
(x,y] \, \text{ for } x, y\in R, \, x<y.
\]
These spaces are clearly $T_3$. Much as in general topology, they work as counterexamples to a number of otherwise plausible conjectures in our setting. We adopt all notational conventions with respect to $\tau_r$ and $\tau_l$ that were previously set for $\tau_e$ and $\tauc$ (i.e. we write $r$ and $l$ in place of $\tau_r$ and $\tau_l$, respectively, when used as subscripts or prefixes).

\begin{definition}\label{dfn:connected}
A definable topological space $(X,\tau)$ is \emph{definably connected} if it is not the union of two disjoint non-empty definable $\tau$-open sets. 

For a fixed definable topological space $(X,\tau)$, we say that a definable subset $Y\subseteq X$ is definably connected if $(Y,\tau)$ is. A \emph{definably connected component} of $(X,\tau)$ is a maximal definably connected definable subset of $X$.  
\end{definition}

Note that another way of stating that the structure $\RR$ is o-minimal is by saying that $(R,\tau_e)$ is definably connected (equivalently, $\RR$ is definably complete) and every definable subset of $R$ can be partitioned into finitely many definably connected euclidean spaces.  

\begin{remark}\label{remark_no_homeomorphism}
Clearly any order reversing bijection $R\rightarrow R$ is a homeomorphism $(R,\tau_r)\rightarrow (R,\tau_l)$. 
Let $\tau_*$ denote either $\tau_r$ or $\tau_l$. Then, for any distinct pair $\tau,\mu\in \{\tau_e, \tau_*, \tauc\}$ and intervals $I,J\subseteq X$, there is no definable homeomorphism $(I,\tau)\rightarrow (J,\mu)$. 
This is obvious if one of the topologies is discrete. 
 If $\{\tau,\mu\}=\{\tau_e,\tau_*\}$, then this follows from the fact that any interval with the euclidean topology is definably connected, while any interval with either the right or left half-open interval topology is totally definably disconnected (i.e. every definably connected subspace is trivial).

\end{remark}

We now introduce a notion of definable separability that generalizes the one given by Walsberg for definable metric spaces in~\cite{walsberg15}. We provide further justification for our definition and its relationship to that of~\cite{walsberg15} in Subsection~\ref{subsection: definable separability}.

\begin{definition}\label{dfn:separable}
We say that a definable topological space $(X,\tau)$ is \emph{definably separable} if there exists no infinite definable family of $\tau$-open pairwise disjoint subsets of $X$.  
\end{definition}

The reader will note the similarity between our definition of definable separability and the countable chain condition (\textbf{ccc}, or Suslin's condition) for topological spaces: a topological space has the \textbf{ccc} if it does not contain an uncountable family of pairwise disjoint open sets. In general, every separable topological spaces has the \textbf{ccc}, but the converse is not true. 

The main justification for our terminology (we discuss further justifications in Subsection~\ref{subsection: definable separability}) is that, in o-minimal expansions of $(\R,<)$, separability, definable separability, and having the \textbf{ccc} are
equivalent notions. To make the present paper self-contained we give a proof of this equivalence for one-dimensional spaces. (The proof in the general case is given in forthcoming work by the first author \cite{ag_separability}.) 

\begin{proposition}\label{prop:sep-equiv}
Suppose that $\RR$ expands $(\R,<)$. Let $(X,\tau)$ be a definable topological space with $\dim X \leq 1$. The following are equivalent. 
\begin{enumerate}
    \item $(X,\tau)$ is definably separable. 
    \item $(X,\tau)$ is separable.
    \item $(X,\tau)$ has the \textbf{ccc}.
\end{enumerate}
\end{proposition}
\begin{proof}
\sloppy Since $\dim X \leq 1$, applying o-minimal cell decomposition we fix throughout a finite family of triples $\{ (I_i, C_i, f_i) : 1 \leq i \leq m\}$ where, for each $1 \leq i \leq m$, $I_i\subseteq \R$ is a singleton or an interval, $C_i$ is a definable subset of $X$, $f_i \colon I_i \rightarrow C_i$ is a definable bijection, and moreover the family $\{ C_i : 1 \leq i\leq m\}$ is a partition of $X$. 

Implication $(2)\Rightarrow(3)$ is a routine exercise. We prove $(1)\Rightarrow(2)$ and $(3)\Rightarrow(1)$.

\textbf{Proof of $(1)\Rightarrow (2)$.}

Suppose that $(X,\tau)$ is definably separable. Let $\BB$ be a definable basis for $\tau$. Observe that any given $x\in X$ has a finite $\tau$-neighbourhood if and only if there exists a finite set $A\in \BB$ with $x\in A$. By o-minimality (uniform finiteness), there exists some $k<\omega$ such that said neighbourhood $A$ is always of size at most $k$. It follows that the set $X_0$ of all $x\in X$ with a finite $\tau$-neighbourhood is definable. Furthermore, observe that, for every $x\in X_0$, there exists a (necessarily unique) $\tau$-neighbourhood of $x$ in $\BB$ of minimum size among all $\tau$-neighbourhoods of $x$. We denote this $\tau$-neighbourhood by $A(x)$. Observe that the family $\{A(x) : x\in X_0\}$ is definable. 

Let $\preccurlyeq_\tau$ denote the classical specialization preorder (reflexive transitive relation) on $X$, which is given by $x \preccurlyeq_\tau y$ whenever $x \in cl_\tau \{y\}$.
Let $\sim_\tau^0$ be the equivalence relation on $X$ where $x\sim_\tau^0 y$ whenever $x$ and $y$ are topologically indistinguishable, i.e. they have the same $\tau$-neighbourhoods. Clearly $\sim_\tau^0$ and $\preccurlyeq_\tau$ are both definable. Observe that $x \sim_\tau^0 y$ holds exactly when $x \preccurlyeq_\tau y$ and $y \preccurlyeq_\tau x$ both do. Let us use the notation $[\cdot]_\tau^0$ to denote the equivalence classes by the relation $\sim_\tau^0$.

Let $X_{\text{max}} \subseteq X_0$ denote the set of points in $X_0$ that are $\preccurlyeq_\tau$-maximal, in the sense that there does not exist some $y\in X$ with $x \preccurlyeq_\tau y$ and $y \not\preccurlyeq_\tau x$. Since $\preccurlyeq_\tau$ is definable then $X_{\text{max}}$ is definable too. 

\begin{claim}\label{claim1:sep-prop}
For every $x\in X_0$, it holds that $X_{\text{max}} \cap A(x) \neq \emptyset$.     
\end{claim}
\begin{proof}[Proof of claim]
The claim follows easily from the observation that, for any $x\in X_0$, it holds that $A(x)=\{ y \in X : x \preccurlyeq_\tau y \}$. In particular the inclusion $A(x)\supseteq \{ y \in X : x \preccurlyeq_\tau y \}$ follows from the definition of $\preccurlyeq_\tau$ and $A(x)\subseteq \{ y \in X : x \preccurlyeq_\tau y \}$ holds by minimality of $A(x)$.
\renewcommand{\qedsymbol}{$\square$ (claim)}
\end{proof}

\begin{claim}\label{claim2:sep-prop}
$X_{\text{max}}$ is finite. 
\end{claim}
\begin{proof}[Proof of claim]
We first show that, for any $x \in X_{\text{max}}$, it holds that $A(x)=[x]_\tau^0$, and hence, in particular, that $[x]_\tau^0$ is $\tau$-open and finite.

Let $x\in X_{\text{max}}$. Clearly $[x]_\tau^0 \subseteq A(x)$. On the other hand, suppose, towards a contradiction, that there exists some $y\in A(x)$ with $y \notin [x]_\tau^0$. However, then $x \preccurlyeq_\tau y$ but $y \not\preccurlyeq_\tau x$, contradicting that $x$ is $\preccurlyeq_\tau$-maximal.

We derive that the family $\{A(x) : x \in X_{\text{max}}\}=\{ [x]_\tau^0 : x \in X_{\text{max}}\}$ is a definable family of identical or pairwise disjoint $\tau$-open sets. Since $(X,\tau)$ is definably separable, this definable family contains finitely many pairwise disjoint sets, each of which is finite, the union of which contains $X_{\text{max}}$. Hence $X_{\text{max}}$ is finite.
\renewcommand{\qedsymbol}{$\square$ (claim)}
\end{proof}

Now recall the notation $\{ (I_i, C_i, f_i) : 1 \leq i \leq m\}$. Let $J\subseteq \{1,\ldots, m\}$ be the set of those $i$ such that $I_i$ is an interval. We define 
\[
D=\bigcup_{i\in J} f_i(\mathbb{Q} \cap I_i) \cup  X_{\text{max}}.
\]

Since $(X,\tau)$ is definably separable, by Claim~\ref{claim2:sep-prop} the set $D$ is countable. We show that it is dense in $(X,\tau)$. Let $A$ be a $\tau$-open set. By making $A$ smaller if necessary we may assume that it is definable. If $A$ is infinite, then there is some $i \in \{1,\ldots, m\}$ such that $A\cap C_i$ is infinite, and so $f_i^{-1}(A\cap C_i)$ contains an interval in $I_i$, and thus $\emptyset \neq f_i(\mathbb{Q} \cap I_i) \cap A \subseteq D\cap A$. Now suppose that $A$ is finite and fix $x\in A$. Then $x \in X_0$ and $A(x)\subseteq A$ and so, by Claim~\ref{claim1:sep-prop}, $\emptyset \neq X_{\text{max}} \cap A(x) \subseteq D \cap A$. 

\textbf{Proof of $(3)\Rightarrow (1)$.}

Suppose that $(X,\tau)$ is not definably separable, witnessed by an infinite definable family of pairwise disjoint $\tau$-open sets $\A$. We show that $\A$ is uncountable. Let $i\leq m$ be such that the family $\{ C_i \cap A : A \in \A\}$ is infinite. We show that the infinite definable family of pairwise disjoint sets $\YY=\{f^{-1}_i(C_i \cap A) : A\in \A\}$ is uncountable. 

For each non-empty set $Y\in \YY$ consider its infimum $\inf Y \in \R \cup\{-\infty\}$. By o-minimality and the fact that the sets in $\YY$ are pairwise disjoint, for any given $t \in \R \cup\{-\infty\}$ at most two sets in $\YY$ can have infinimum equal to $t$. So the set $\R\cap \{ \inf Y : Y\in \YY\setminus \{\emptyset\}\}$ is infinite. Observe that this set is also definable and so, by o-minimality, it contains an interval, which implies that it is uncountable, and it follows that $\YY$ is uncountable too.     
\end{proof}

We now give a simple characterization of definable separability for $T_1$ one-dimensional definable topological spaces.

\begin{lemma}\label{lem:T1-sep}
Let $(X,\tau)$ be a $T_1$ definable topological space with $\dim X\leq 1$. Then $(X,\tau)$ is definably separable if and only if it has finitely many $\tau$-isolated points. 
\end{lemma}
\begin{proof}
The set of all $\tau$-isolated points in $(X,\tau)$ is clearly definable. If it is infinite, then $(X,\tau)$ fails to be definably separable. 

Conversely, suppose that $(X,\tau)$ is not definably separable, witnessed by an infinite definable family of pairwise disjoint $\tau$-open sets $\A$. Since $\dim X\leq 1$, by the Fiber Lemma for o-minimal dimension \cite[Chapter 4, Proposition 1.5 and Corollary 1.6]{dries98}, the family $\A$ contains only finitely many infinite sets, hence infinitely many finite sets. Since $\tau$ is $T_1$, any finite $\tau$-open subset of $X$ must contain only $\tau$-isolated points. It follows that $(X,\tau)$ has infinitely many $\tau$-isolated points.  
\end{proof}

Following the above lemma, we now consider definable separability in the context of the three fundamental topologies that have been introduced so far.
\begin{proposition} \label{prop:def-sep_eucl_disc_llt}
Let $X\subseteq R^n$ be a definable set.
\begin{enumerate}[(a)]
    \item  \label{itm1:sep-spaces-examples} The space $(X,\tau_e)$ is definably separable. 
    \item \label{itm2:sep-spaces-examples} The space $(X,\tauc)$ is definably separable if and only if $X$ is finite. 
    \item \label{itm3:sep-spaces-examples} Suppose that $n=1$. The spaces $(X,\tau_r)$ and $(X,\tau_r)$ are definably separable. 
\end{enumerate}
\end{proposition}
\begin{proof}
Statement~\eqref{itm2:sep-spaces-examples} is obvious and
statement~\eqref{itm3:sep-spaces-examples} follows immediately from Lemma~\ref{lem:T1-sep}, o-minimality and the fact that the topologies $\tau_r$ and $\tau_l$ are $T_1$.

Statement~\eqref{itm1:sep-spaces-examples} also follows immediately in the case $n=1$, by a similar argument to case~\eqref{itm3:sep-spaces-examples}. In order to prove~\eqref{itm1:sep-spaces-examples} in general, suppose that $(X,\tau_e)$, $X \subseteq R^n$, is not definably separable, witnessed by an infinite definable family of open pairwise disjoint sets $\mathcal{A}$. By o-minimal cell decomposition there exists at least one cell $C\subseteq X$ such that the family $\{ A \cap C : A\in \A\}$ is infinite, and in particular the subspace $(C,\tau_e)$ is not definably separable. Recall that every cell is definably $e$-homeomorphic to an open cell. Hence to prove~\eqref{itm1:sep-spaces-examples} it suffices to show that open cells with the euclidean topology are definably separable. 

Towards a contradiction let $C$ be an open cell such that $(C,\tau_e)$ is not definably separable, witnessed by an infinite definable family of open pairwise disjoint sets $\mathcal{A}$. By o-minimality every definable non-empty open subset of $C$ has dimension equal to $\dim C$, and in particular this holds for every non-empty set in $\A$. Applying the Fiber Lemma for o-minimal dimension \cite[Chapter 4, Proposition 1.5 and Corollary 1.6]{dries98} we conclude that $\dim \bigcup\mathcal{A} > \dim C$, which contradicts that $\bigcup \A \subseteq C$. 
\end{proof}

We now introduce definable curves, which play a crucial role in the study of o-minimal definable topologies, often taking the role that sequences have in general topology.

\begin{definition} \label{dfn:curve}
Let $(X,\tau)$ be a definable topological space. A \emph{curve in $X$} is a map $\gamma: (a,b)\rightarrow X$, where $a,b\in \exR$, $a<b$. 

We say that \emph{$\gamma$ converges in the $\tau$-topology} (or \emph{converges in $(X,\tau)$}, or \emph{$\tau$-converges})  \emph{as $t$ tends to $a$} to a point $x\in X$ if, for every $\tau$-neighbourhood $U$ of $x$, there exists some $a<t_U<b$ such that $\gamma(t)\in U$ for all $a<t<t_U$. 
In this case, we say that $x$ is a \emph{$\tau$-limit of $\gamma$ as $t$ tends to $a$} and, if this limit is unique (which will certainly be the case if $\tau$ is Hausdorff), then we write $x=\taulim_{t\rightarrow a}\gamma(t)$. Convergence \emph{as $t$ tends to $b$} is defined analogously. When we say that $\gamma$ \emph{$\tau$-converges} to $x\in X$, without reference to $a$ or $b$, it should be understood that we have already implicitly fixed an endpoint $c\in \{a,b\}$, which we will call the \emph{convergence endpoint} of $\gamma$, and are saying that $\gamma$ $\tau$-converges to $x$ as $t$ tends to $c$. We say that $\gamma$ is \emph{$\tau$-convergent} if it $\tau$-converges to some $x\in X$ (as $t$ tends to its convergence endpoint). 
\end{definition}

\begin{remark}\label{remark_assumptions_curves}
We adopt some further conventions regarding definable curves.
Let $\gamma:(a,b)\rightarrow R^n$ be a definable curve. Frequently, we are only concerned with the behaviour of $\gamma$ near its convergence endpoint $c \in \{a,b\}$. 
By o-minimality, for any definable set $X\subseteq R^n$ there exists $a'>a$ such that either $\gamma[(a,a')]\subseteq X$ or $\gamma[(a,a')]\subseteq R^n\setminus X$, and the analogous statement holds for $b$. If (say) $c=a$ and there is some $a'>a$ such that $\gamma[(a,a')]\subseteq X$, we may treat $\gamma$ as a curve in $X$ by implicitly identifying it with its restriction $\gamma|_{(a,a')}$.

Similarly, we will adopt the convention of saying that $\gamma$ is constant or injective (or some other property) if it has this property when restricting its domain to an appropriate interval as above. By o-minimality, every definable curve $\gamma:(a,b)\rightarrow R^n$ can be assumed to be either constant or injective (strictly monotonic if $n=1$) and \econtinuousp.

Whenever we say that $\gamma$ $\tau$-converges and $\mu$-converges, for two definable topologies $\tau$ and $\mu$, and without explicit reference to a convergence endpoint $c\in \{a,b\}$, it should be understood that the same endpoint is being considered in both cases.
\end{remark}

\begin{definition}\label{dfn:compact}
A definable topological space $(X,\tau)$ is \emph{definably compact} if any definable curve $\gamma:(a,b)\rightarrow X$ $\tau$-converges as $t$ tends to $a$ and as $t$ tends to $b$. 
\end{definition}

Definition~\ref{dfn:compact} is motivated by the definition of definable compactness introduced in~\cite{pet_stein_99}. Adapting the terminology in~\cite{pet_stein_99}, we would say that $(X,\tau)$ is definably compact if and only if every definable curve $\gamma$ in $X$ is $\tau$-completable. In our terminology, this corresponds to the property that every definable curve (with any convergence endpoint) $\tau$-converges.

It is easy to see that, for a given infinite definable set $X\subseteq R$, the spaces $(X,\tau_r)$, $(X,\tau_l)$ and $(X,\tauc)$ are not definably compact, and the space $(X,\tau_e)$ is definably compact if and only if $X$ is $e$-closed and bounded. The fact that euclidean spaces are definably compact if and only if they are closed and bounded was proved for sets of all dimensions in~\cite[Theorem 2.1]{pet_stein_99}. 

\sloppy Other notions of definable compactness besides the one in Definition~\ref{dfn:compact} have been studied in the o-minimal and other model-theoretic contexts. One of them is the following (recall that a family of sets $\BB$ is downward directed if for every pair $B, B'\in\BB$ there exists some $B''\in\BB$ such that $B''\subseteq B \cap B'$).
\begin{equation}\label{dfn:directed-compact}
\parbox{0.9\textwidth}{
    Every downward directed definable family of non-empty closed sets has non-empty intersection.
}
\end{equation}
This notion has been studied in the o-minimal context by Johnson~\cite{johnson14}, the first author\cite{ag_FTT}, and the authors and Walsberg~\cite{atw1}. It has also been studied in the setting of p-adically closed fields by the first author and Johnson~\cite{ag-j-22}. In the more general model-theoretic context it has been approached by Fornasiero~\cite{fornasiero}. Johnson proved that for o-minimal euclidean spaces it is equivalent to being closed and bounded~\cite[Proposition 3.10]{johnson14}. 

\begin{remark}\label{rem:compactness}
\sloppy Definition~\ref{dfn:compact} and condition~\eqref{dfn:directed-compact} are equivalent for one-dimensional definable topological spaces. The fact that condition~\eqref{dfn:directed-compact} implies Definition~\ref{dfn:compact} can be seen from \cite[Corollary 44]{atw1}, (3) $\Rightarrow$ (1) (that corollary is presented in the setting of o-minimal expansions of ordered groups, but this part of the statement does not require that assumption). The fact that Definition~\ref{dfn:compact} implies condition~\eqref{dfn:directed-compact} can be obtained from the results of the present paper (see Remark~\ref{remark_curve_selection}). The details are all presented in~\cite[Proposition 6.2.4]{andujar_thesis}.

For spaces of any dimension, we show together with Walsberg in~\cite{atw1} that the equivalence holds whenever $\RR$ expands an ordered group, and more generally in~\cite{ag_FTT} the first author shows that they are equivalent whenever $\RR$ has definable choice, and also without any assumption on $\RR$ beyond o-minimality whenever the topology of the space is Hausdorff. 

Since the focus on the present paper is that of definable topological spaces of dimension at most one, we will usually be in a position to assume equivalence of the two definitions. We will largely work with Definition~\ref{dfn:compact} (but make clear when we are instead working for convenience with condition~(\ref{dfn:directed-compact})). 
\end{remark}

While it is easy to see that compactness implies definable compactness, the converse is not true in general. Nevertheless, both notions are equivalent whenever $\RR$ expands the field of reals~\cite[Corollary 48]{atw1}. In~\cite{ag_FTT} the first author shows that, if $\RR$ expands $(\R,<)$, then a definable topological space is compact if and only if it is definably compact in the sense of condition~\eqref{dfn:directed-compact}. By Remark~\ref{rem:compactness} it follows that, if $\RR$ expands $(\R,<)$, then definable compactness and compactness are equivalent notions for definable topological spaces of dimension at most one. 

We will use the following notions frequently throughout this paper.
\begin{definition}\label{dfn:push-forward}
Let $(X,\tau)$ be a definable topological space with definable basis $\BB$ and let $f:X\rightarrow R^m$ be an injective definable map. We define the \emph{push-forward of $(X,\tau)$ by $f$} to be the definable topological space $(f(X),f(\tau))$, where $f(\tau)$ is the topology on $f(X)$ with definable basis $\{f(A) : A \in \BB\}$. Thus $f(\tau)$ is the topology satisfying that \mbox{$f:(X,\tau)\rightarrow (f(X),f(\tau))$} is a homeomorphism.
Given a bijective definable partial map $g:R^n\rightarrow X$, we define the \emph{pull-back of $(X,\tau)$ by $g$} to be its push-forward by $g^{-1}$.  

\end{definition}

Suppose that we have finitely many topological spaces $(X_i,\tau_i)$, where, for each $0\leq i \leq k$, we have $X_i\subseteq R^n$ and $(X_i,\tau_i)$ has basis $\BB_i$. Then their disjoint union is the set $\bigcup_{0\leq i \leq k} (\{i\}\times X_i )$ with topology given by the basis $\bigcup_{0\leq i\leq k} (\{i\}\times \BB_i )$. Note that, for each $i$, the map $X_i\rightarrow \bigcup_{0\leq i\leq k} (\{i\}\times X_i)$, given by $x\mapsto \al i,x\ar$, is an open embedding. The union of finitely many definable topological spaces is clearly definable. If the sets $X_i$ are not all part of the same ambient space $R^n$ and we wish to consider their disjoint union, we first identify them with their product with singletons through the natural push-forward in order to assume that they are.

\section{Definable metric spaces}\label{section: metric spaces}

In this section, we recall the definition of definable metric spaces, introduce the notion of definable metrizability, and discuss further our definition of definable separability (Definition~\ref{dfn:separable}) in light of the notion with the same name introduced for definable metric spaces by Walsberg~\cite{walsberg15}.

Throughout this section, we suppose that $\RR$ expands an ordered group \mbox{$(R,0,+,<)$}. In the spirit of the definition of $\mathcal{M}$-norm introduced by the second author in~\cite{thomas12} we include the following definition. 

\begin{definition}\label{dfn:R-metric}
Let $X$ be a set. An \emph{$\RR$-metric on $X$} is a map $d:X\times X\rightarrow R^{\geq 0}$ that satisfies the metric axioms, i.e. identity of indiscernibles, symmetry and subadditivity.
\end{definition}

We now recall the following definition of a  definable metric space from \cite{walsberg15}. Although Walsberg works under the assumption that $\RR$ is an o-minimal expansion of an ordered field, the following definition, as well as any other notion that we borrow from~\cite{walsberg15}, still makes sense in the ordered group setting.

\begin{definition}[\cite{walsberg15}]
A \emph{definable metric space} is a tuple $(X,d)$, where $X$ is a definable set and $d$ is a definable $\RR$-metric on $X$. 
\end{definition}

Any $\RR$-metric $d$ generates a topology in the usual way, which we denote by $\tau_d$. Following the conventions set in Section~\ref{section:definitions} for the euclidean,  discrete, right and left half-open interval topologies, we sometimes abuse notation and write $d$ in place of $\tau_d$. 

\subsection{Definable metrizability} \label{subsection: definable metrizability}

It is easy to prove that any topology generated by an $\RR$-metric is Hausdorff and regular. Any  definable $\RR$-metric induces a definable topology, and so by this identification every definable metric space is a definable topological space. Much as with the notion of metrizability in general topology, the converse is not true, i.e. there are definable topologies (including Hausdorff regular topologies) that do not arise from definable $\RR$-metrics. Hence it is reasonable to investigate which topologies have this property. This motivates the following definitions.

\begin{definition}\label{dfn:metrizable}
A topological space $(X,\tau)$ is \emph{$\RR$-metrizable} if there exists an $\RR$-metric $d$ on $X$ such that $\tau_d=\tau$ and \emph{definably metrizable} if there exists some definable $\RR$-metric $d$ on $X$ such that $\tau_d=\tau$. 
\end{definition}

We shall simplify our terminology throughout to refer to metrics, rather than $\RR$-metrics, and similarly to metrizability, rather than $\RR$-metrizability, without any loss of clarity.

Both the euclidean and discrete topologies are definably metrizable (and hence metrizable) on any definable set. 

A basic example of a definable topological space that is not metrizable would be any non-Hausdorff definable topological space, e.g. the Sierpinski space $X=\{0,1\}$, $\tau=\{\emptyset, \{1\},\{0,1\}\}$.  

An example of a space that is not definably metrizable but displays all the separation axioms of definable metric spaces would be the space $(R,\tau_r)$, as we now show. 
\begin{proposition}\label{prop:tau_r not metrizable}
The space $(R,\tau_r)$ is not definably metrizable. 
\end{proposition}
\begin{proof}
If $(R,\tau_r)$ were definably metrizable with definable metric $d$, then there would exist, for every $x\in R$, some $\varepsilon_x>0$ such that $B_d(x,\varepsilon_x) \subseteq [x,+\infty)$, where $B_d(x,\varepsilon_x)$ is the $d$-ball of radius $\varepsilon_x$ and center $x$. Let $1$ denote some fixed positive element of $R$ and let $f:X\rightarrow (0,\infty)$ be the definable map given by $f(x)=\sup\{ t\leq 1 : B_d(x,t)\subseteq [x,+\infty)\}$. By o-minimality, there exists an interval $I\subseteq R$ such that, for some $\varepsilon>0$, we have $f(x)>\varepsilon$, for all $x\in I$. Hence, for any distinct $x,y\in I$, it holds that $d(x,y) \geq \varepsilon$, i.e. $(I,\tau_r)$ is a discrete space, which contradicts the definition of $\tau_r$. 
\end{proof}
The above result still holds if we consider any infinite definable set $X\subseteq R$ in place of $R$ and also if we put $\tau_l$ in place of $\tau_r$. It is worth noting that, in the particular case where $R=\mathbb{R}$, the Sorgenfrey Line is separable but not second countable, and thus it is not even metrizable. On the other hand, if $\RR=(\mathbb{Q},+,<)$, then the space $(\mathbb{Q},\tau_r)$ is metrizable (meaning $(\mathbb{Q},+,<)$-metrizable) (see \cite{diep17}). 

In Section~\ref{section:metrizability}, we address the question of which metrizable definable topological spaces are definably metrizable, and give a characterization for o-minimal expansions of certain ordered fields including $(\R,+,\cdot,<)$ (Theorem~\ref{thm:metrizability}).

\subsection{Definable separability} \label{subsection: definable separability}

We now turn to justifying our notion of definable separability (Definition~\ref{dfn:separable}), in light of a similar definition given by Walsberg for definable metric spaces in~\cite[Section 7.1]{walsberg15}. There, it is stated that a definable metric space $(X,d)$ is definably separable if there exists no infinite definable subset $Y\subseteq X$ such that the subspace topology $\tau_d|_Y$ is discrete.

First, we note that this definition is similar to ours, in that it simply asks that the property we define be hereditary (i.e. that every definable subspace of a definable metric space is definably separable in our sense). In fact, when restricted to the context of definable metric spaces, Walsberg's definition and the one in this paper are equivalent, as shown by the following. 

\begin{lemma}\label{remark_def_sep}
Let $(X,d)$ be a definable metric space. Then $(X,d)$ is definably separable (in the sense of Definition~\ref{dfn:separable}) if and only if there exists no infinite definable discrete subspace. 
\end{lemma}
\begin{proof}

Let $(X,d)$ be a definable metric space. Let $Y \subseteq X$ be an infinite definable discrete subspace of $X$. Since $\RR$ expands an ordered group, it has definable choice, so one may definably select, for each $x\in Y$, some $\varepsilon_x>0$ such that, for every $y\in Y\setminus \{x\}$, $2\varepsilon_x \leq d(x,y)$. From the triangle inequality, it follows that the infinite definable family of open $d$-balls $\{B_d(x,\varepsilon_x) : x\in Y\}$ is pairwise disjoint, hence $(X,d)$ is not definably separable. The other direction follows immediately from definable choice.
\end{proof}

From the above lemma it follows that, for the class of definable metric spaces, definable separability is a hereditary property. Moreover note that the proof of the lemma relies solely on definable choice and not on the fact that the structure $\RR$ is o-minimal.

Further to our discussion of definable separability in Section~\ref{section:definitions} (see Proposition~\ref{prop:sep-equiv}), our notion of definable separability is further justified by the following observation. If we drop the assumption of definable metrizability, then there are definable topological spaces that are definably separable (according to our definition) but still contain an infinite definable discrete subspace. Some examples are the generalizations to $\RR$ of the Moore Plane (defined assuming ordered field structure in $\RR$, see Appendix \ref{section:examples}, Example~\ref{example:Moore_plane}) or the Sorgenfrey Plane (the product of the Sorgenfrey Line with itself); 
see also Example~\ref{example: dfbly sep not hereditary} in Appendix~\ref{section:examples} for an example of dimension one. Hence definable separability is not in general a hereditary property, and Walsberg's definition turns out to be strictly stronger than ours in the general context. This is in accordance with general topology, where every subspace of a separable metric space is separable, but where the Moore Plane and the Sorgenfrey Plane are examples of separable topological spaces with uncountable discrete subspaces.

Walsberg showed in~\cite[Theorem 7.1]{walsberg15}, that, whenever $\RR$ expands an ordered field, any definably separable definable metric space is definably homeomorphic to a euclidean space. This result does not generalize to all $T_3$ definable topologies, as witnessed by the right and left half-open interval topologies. In Section~\ref{section: affine}, we address the question of which one-dimensional definable topological spaces are, up to definable homeomorphism, euclidean.

\section{Preliminary results} \label{section: prel. results}

From now until the end of Section~\ref{section:compactifications}, we return to the general setting in which $\RR$ is an o-minimal expansion of a dense linear order without endpoints.

\subsection{Spaces of all dimensions}\label{subsection: prelim_results_alldim}

In this section we include some preliminary results concerning definable topologies of all dimensions. 

The purpose of the following definition is to have a tool that allows us to study definable topologies in the o-minimal setting. This will in particular allow us to provide a description of a basis of neighbourhoods of any given point in a $T_1$ definable topological space (see Lemma~\ref{lemma_explaining_neighbourhoods}).

\begin{definition}
Let $(X,\tau)$, with $X\subseteq R^n$, be a definable topological space. Let $x\in X$ and let $\BB(x)$ be a basis of $\tau$-neighbourhoods of $x$. We define the \emph{$e$-accumulation set of $x$ in $(X,\tau)$}, namely $\E^{(X,\tau)}_x$, to be:
\[
\E^{(X,\tau)}_x:= \bigcap_{A\in\BB(x)} \{ y\in \exR^n : A\cap B\setminus \{y\} \neq \emptyset, \text{ for all } B\in\tau_e \text{ with } y\in B\},
\]
where $\tau_e$ refers to the euclidean topology in $\exR^n$.
So $\E^{(X,\tau)}_x$ is the intersection of the set of $e$-accumulation points of every $\tau$-neighbourhood of $x$. 

If $(X,\tau)$ is $T_1$ and $x, y \in X$ with $x\neq y$, then $y\in \E^{(X,\tau)}_x$ is equivalent to stating that, for every $\tau$-neighbourhood $A$ of $x$ and every $e$-neighbourhood $B$ of $y$, $A\cap B\neq \emptyset$, i.e. $y \in cl_e A$, for every $\tau$-neighbourhood $A$ of $x$. 

\end{definition} 

The definition of $\E^{(X,\tau)}_x$ is clearly independent of the choice of basis of neighbourhoods. Note that, if an element $x\in X$ is $\tau$-isolated, then it satisfies $\E^{(X,\tau)}_x=\emptyset$. The converse is not true in general; it is however true for $T_1$ spaces. In the case where $X\subseteq R$, this will follow immediately from Lemma~\ref{lemma_explaining_neighbourhoods} and Proposition~\ref{prop_basic_facts_P_x_2}(\ref{itm2:basic_facts_P_x_2}).
We leave it to the reader to check that the implication holds in higher dimensions, using Remark~\ref{rem:compactness} and the fact that the space $(\exR^n, \tau_e)$ for any $n$ is definably compact.

For any point $x$ in a euclidean space $(X,\tau_e)$, with $X \subseteq R^n$, it holds that $\E^{(X,\tau_e)}_x=\{x\}$.

Generally, since there will be no room for confusion, once a definable topological space $(X,\tau)$ is fixed, then, for any $x\in X$, we will write $\E_x$ in place of $\E^{(X,\tau)}_x$, and will only resort to the latter when we also intend to address the $e$-accumulation set $\E^{(Y,\tau)}_x$ for some definable subspace $Y$ containing $x$.

The following are facts regarding $e$-accumulation sets that follow immediately from the definition. Recall that the euclidean topology is understood in $\exR$. 

\begin{proposition}\label{basic_facts_P_x}
Let $(X,\tau)$ be a definable topological space. 
\begin{enumerate}[(a)]
\item \label{itm: basic_facts_1} $\E_x$ is $e$-closed and $\E_x \subseteq cl_e X$ for every $x\in X$.
\item \label{itm: basic_facts_2} The relation $\{\al x, y\ar: y\in \E_x \}\subseteq X\times \exR^m$ is definable.  
\end{enumerate}
\end{proposition}

\sloppy We now prove a bound (Lemma~\ref{lemma_2}) on the dimension of $e$-accumulation sets (in particular, that they are finite when $\dim(X) \leq 1$). In order to do so, we need two technical lemmas. For the first of these, recall that the dimension of a definable set $X$ is equal to the dimension of the interpretation $X^*$ of $X$ in any elementary extension $\RR^*$ of $\RR$. 
  
\begin{lemma}\label{lemma_1}
Let $\RR^*$ be an elementary extension of $\RR$. Let $X\subseteq R^n$ be an $\RR$-definable set and let $X^*$ denote the interpretation of $X$ in $\RR^*$. If $Y\subseteq X^*$ is $\RR^*$-definable and $X\subseteq Y$ (in $(R^*)^n$) then $\dim Y = \dim X$ (where $\dim X$ denotes the dimension of $X$ as a definable set in $\RR$). 
\end{lemma}
\begin{proof}
Since $Y\subseteq X^*$, we have that $\dim Y \leq \dim X^*=\dim X$. We show, by induction on $n$, that from $X\subseteq Y$ it follows that $\dim X \leq \dim Y$. Since otherwise $X=Y=X^*$, we may assume that $X$ is infinite.

Suppose that $n=1$. Since $X$ is infinite, we have that $Y$ is infinite and so, by o-minimality, it must contain an interval. In particular, $\dim Y \geq 1$, and so the result follows.

Now suppose that $n>1$. By passing to a cell inside $X$ of maximal dimension if necessary, we may assume that $X$ is a cell. For any given set $S\subseteq R^n$, let $S'$ denote the projection to the first $n-1$ coordinates. For any $x\in Y'$, let $Y_x=\{t : \langle x,t \rangle\in Y\}$.  

If $\dim X'=\dim X$, then we are done, since $\dim Y' \leq \dim Y$ and by induction hypothesis $\dim X'\leq \dim Y'$. Otherwise, $\dim X= \dim X'+1$ and $X$ is of the form $(f,g)_{X'}$, for continuous functions $f,g:X'\rightarrow \exR$ with $f<g$. In this case, consider the definable set $Y_{\text{inf}}=\{x\in Y': Y_x \text{ is infinite} \}$. Then $X'\subseteq Y_{\text{inf}}$ and, by induction hypothesis, $\dim X' \leq \dim Y_{\text{inf}}$. By the Fiber Lemma for o-minimal dimension \cite[Chapter 4, Proposition 1.5 and Corollary 1.6]{dries98}, $\dim Y_{\text{inf}} + 1 \leq \dim Y$. We conclude that $\dim X= \dim X' +1 \leq \dim Y_{\text{inf}} +1 \leq \dim Y$. 
\end{proof}

For the next lemma, recall that a family of sets $\SSS$ has the finite intersection property if $\bigcap \FF\neq \emptyset$ for every finite subfamily $\FF\subseteq \SSS$. 

\begin{lemma}\label{lemma_FIP}
Let $\{ S_u\subseteq R^n : u\in \Omega\}$ be a definable family with the finite intersection property. Then there exists $\Sigma \subseteq \Omega$ with $\dim \Sigma =\dim \Omega$ such that $\bigcap \{S_u : u\in \Sigma\}\neq \emptyset$.
\end{lemma}
\begin{proof}
Let $\RR^*=(R^*,<,\ldots)$ denote an $|R|^+$-saturated elementary extension of $\RR$. Then, by saturation, there exists $x_0\in (R^*)^n$ such that $x_0\in S^*_u$, for every $u\in \Omega$. By Lemma~\ref{lemma_1}, the $\RR^*$-definable set $\Sigma^*_{x_0}=\{u\in \Omega^* : x_0\in S^*_u\}$ has dimension equal to $\dim \Omega$. For each $x\in R^n$, let $\Sigma_x=\{ u\in \Omega : x\in S_u\}$. Since $\dim \Sigma^*_{x_0} = \dim \Omega$ and $\RR \preceq \RR^*$, there must exist some $x\in R^n$ such that $\dim \Sigma_x=\dim \Omega$. 
\end{proof}

Let $(X,\tau)$ be a definable topological space. Let $x\in X$ and let $\UU$ be a definable basis of neighbourhoods of $x$ in $(X,\tau)$. For the next lemma, we define the \emph{local dimension of $(X,\tau)$ at $x$} to be 
\[
\dim_x(X,\tau)=\min\{\dim A : A\in\UU\}.
\]
Clearly the definition of local dimension does not depend on the choice of basis of neighbourhoods. This definition generalizes the definition of local dimension of a definable metric space at a point that was introduced by Walsberg in~\cite{walsberg15}. 

\begin{lemma}\label{lemma_2}
Let $(X,\tau)$ be a $T_1$ definable topological space. For any $x\in X$, $\dim(E_x)<\dim_x(X,\tau)$. In particular, when $\dim X\leq 1$, the set $\E_x$ is finite for every $x\in X$.
\end{lemma} 
\begin{proof}
Towards a contradiction, suppose that there exists $x\in X$ such that $\dim \E_x \geq\dim_x(X,\tau)$. Let $\{A_u : u\in \Omega\}$ be a definable basis of $\tau$-neighbourhoods of $x$. If $\dim_x(X,\tau)=0$, then, by definition of $\E_x$, we have that $\E_x=\emptyset$, which contradicts the fact that $\dim \E_x \geq\dim_x(X,\tau)$, so we may assume that $\dim_x(X,\tau)>0$. 

Let $n=\dim_x(X,\tau)$. We have that $\dim \E_x \geq n >0$, and in particular that $\dim E_x = \dim (E_x \setminus \{x\})$. For any $y\in E_x \setminus \{x\}$, let $\Omega_y=\{u\in \Omega : y\notin A_u\}$. Since $(X,\tau)$ is $T_1$, the sets $\Omega_y$ are non-empty and in fact the definable family $\{\Omega_y : y\in E_x \setminus \{x\}\}$ has the finite intersection property. 

By Lemma~\ref{lemma_FIP}, there exists a definable set $B\subseteq \E_x\setminus\{x\}$ with $\dim B=\dim E_x$ and there exists $u\in \Omega$ such that $A_u \cap B=\emptyset$. By shrinking $A:=A_u$ if necessary, we may assume that $\dim A = n$. Note however that, by definition of $\E_x$, $B \subseteq cl_e A $, and so $B\subseteq \partial_e A $. In particular $\dim \E_x = \dim B \leq \dim \partial_e A$. However, by o-minimality, $\dim \partial_e A < \dim A = n$, a contradiction.  
\end{proof}

We end this subsection with a remark which we will use extensively throughout the paper. In particular, it allows us to make the assumption, whenever $\RR$ expands an ordered field, that any definable topological space of dimension at most one is, up to definable homeomorphism, a bounded subset of $R$.

\begin{remark}\label{remark_assumption_X}

Let $X$ be a definable set and $n>0$ be such that $\dim(X) \leq n$. If $\RR$ expands an ordered group and $X$ is bounded, then there exists a definable injection $f: X \rightarrow R^n$. In particular, if $\tau$ is a definable topology on $X$ then, by passing to the push-forward of $(X,\tau)$ by $f$ if necessary, one may always assume, up to definable homeomorphism, that $X \subseteq R^n$.

If, moreover, $\RR$ expands an ordered field, then such an injection $f$ exists without the assumption that $X$ is bounded, and with the added condition that $f(X) \subseteq (0,1)^n$. In particular one may always assume, up to definable homeomorphism, that $X \subseteq (0,1)^n$. 

The existence of such injections in each case can been seen as follows.  
Suppose that $\RR$ expands an ordered group and that $X$ is bounded. 
Let $\mathcal{X}$ be a finite partition of $X$ into cells. By o-minimality, each cell in $C \in \mathcal{X}$ is in bijection, under an appropriate projection $\pi_C$, with a subset of $R^{\dim C}$. Let $f_C:C\rightarrow R^n$ be the definable map given by $x \mapsto \pi_C(x)$ if $\dim C=n$, and otherwise $x \mapsto \pi_C(x) \times \{ a_C \}$, where $a_C \in R^{n-\dim C}$ is a given fixed parameter. Note that $f_C$ is an injection into $R^n$. Since $X$ is bounded, then so are the sets $f_C(C)$ for $C \in\mathcal{X}$. As $\RR$ is an expansion of a group, we can find appropriate translations of these sets such that they do not intersect each other. The union of the image of these translations is then in definable bijection with $X$.

Furthermore note that, whenever $\RR$ expands an ordered field, the map given coordinate-wise, for every $n$, by $x_i\mapsto \frac{2x_i-1}{(2x_i-1)^2-1}$, for $i=1,\ldots,n$, gives a definable \ehomeomorphism from $(0,1)^n$ to $R^n$.
\end{remark}

If $\RR$ does not expand an ordered field, then it is not true in general that a definable set $X$ with $\dim X = 1$ is in definable bijection with a subset of $R$. In this case, however, o-minimal cell decomposition and the observations made in Remark~\ref{remark_generalizing_Rx_Lx} below will suffice to generalize many of the results in this paper concerning spaces $(X,\tau)$ with $X\subseteq R$ to one-dimensional definable topological spaces. 

\subsection{One-dimensional spaces} \label{subsection: prelim_results_dim1}

We now focus on preliminary results about definable topological spaces $(X,\tau)$ where $\dim X\leq 1$, or at times more specifically when $X\subseteq R$, which we informally refer to as \textquotedblleft spaces in the line".

The following lemma shows that definable curve selection is a property of all one-dimensional definable topological spaces (if we drop the requirement that the curves be continuous). In~\cite{atw1}, we prove together with Walsberg that it holds for definable topological spaces of all dimensions when $\RR$ expands an ordered field, but not in general e.g. when $\RR$ expands an ordered group (although in this case one may use definable choice to show that it still holds for definable metric spaces, arguing similarly to the proof of \cite[Proposition 41]{atw1}).

\begin{lemma}[Definable curve selection]\label{lemma_curve_selection_R}
Let $(X,\tau)$, $\dim X\leq 1$, be a definable topological space. Then $(X,\tau)$ has definable curve selection, that is, for any $x\in X$ and definable set $Y\subseteq X$, $x\in cl_\tau Y$ if and only if there exists a definable curve $\gamma:I \rightarrow Y$ in $Y$ that $\tau$-converges to $x$.
\end{lemma}
\begin{proof}
It follows readily from the definition of curve convergence that, if there exists a curve in $Y$ $\tau$-converging to $x \in X$, then $x \in cl_{\tau}Y$. We prove the converse.

Fix $x\in X$ and a definable set $Y\subseteq X\subseteq R^n$ with $x\in cl_\tau Y$.
Let $\UU$ denote a definable basis of $\tau$-neighbourhoods of $x$ and set $\BB:=\{U\cap Y : U\in\UU\}$. 
It suffices to prove the existence of an interval $I$ with an endpoint $c\in\exR$ and a definable curve $\gamma: I \rightarrow R^n$ such that, for every $B\in \BB$, $\gamma(t)\in B$ for $t\in I$ close enough to $c$. (We may then restrict $\gamma$ to a suitable subinterval close to $c$ to ensure that it maps only into $Y$.) Since otherwise the proof is immediate, we may assume that $x \notin Y$ and hence, for any $y\in Y$, there is $B\in\BB$ such that $y\notin B$.

First we consider the case where $X\subseteq R$, finding an interval $I$ on which we may take $\gamma$ to be the identity. Consider the definable set \mbox{$H=\{ t\in R : \exists B\in\BB, B\cap (-\infty,t]= \emptyset \}$}. If $H$ is empty then, by o-minimality, for every $B\in\BB$, there is $t_B\in R$ such that $(-\infty, t_B)\subseteq B$, in which case we may take $I=R$ and $c=-\infty$. 

Now suppose that $H$ is non-empty. Note that $H$ is an interval in $R$ (possibly right closed) which is unbounded from below. Let $c=\sup H\in R\cup\{+\infty\}$.

If $c=\max H$, then there exists $B_c\in \BB$ such that $B_c\cap (-\infty,c]= \emptyset$. Let $B\in \BB$. If there exists $s_B>c$ such that $(c,s_B)\cap B=\emptyset$, then any set $B'\in \BB$ with $B'\subseteq B_c \cap B$ satisfies that $(-\infty, s_B)\cap B' =\emptyset$, contradicting that $c=\sup H < s_B$. Hence, by o-minimality, for any $B\in\BB$, there exists $t_B>c$ such that $(c,t_B)\subseteq B$. So let $I=(c,+\infty)$.  

If $c\notin H$, then, for every $B\in \BB$, $B \cap (-\infty, c]\neq \emptyset$. Let $B\in \BB$. Suppose that there exists $s_B<c$ such that $(s_B,c)\cap B=\emptyset$. By assumptions on $\BB$, we may assume that $c\notin B$. Since $s_B\in H$ there is $B'\in \BB$ such that $B' \cap (-\infty,s_B]=\emptyset$. But then any $B''\in\BB$ with $B''\subseteq B \cap B'$ satisfies that $(-\infty, c]\cap B'' =\emptyset$, contradicting that $c\notin H$. Hence, by o-minimality, for any $B\in\BB$ there exists $t_B<c$ such that $(t_B,c)\subseteq B$. So let $I=(-\infty,c)$. 

Now, in the case where $X$ is not a subset of $R$, let us pass, by o-minimal cell decomposition, to a cell $Y'\subseteq Y$ such that $x\in cl_\tau Y'$, where $Y'$ is definably homeomorphic to an interval $J$ by a function $f:J\rightarrow Y'$. Then we may apply the above argument taking the family $\{f^{-1}(U\cap Y') : U\in\UU\}$ in place of $\BB$ to reach an interval $I$ and an endpoint $c$ of $I$ such that, for every $U\in \UU$, we have $t\in f^{-1}(U\cap Y')$, for $t\in I$ close enough to $c$. Finally, we take $\gamma=f|_{I\cap J}: I\cap J \rightarrow Y'$. 
\end{proof}

Note that, in Lemma~\ref{lemma_curve_selection_R}, we may relax the condition $\dim X \leq 1$ and instead draw the same conclusion, for any $x \in X$, as long as $\dim_x(X, \tau) \leq  1$ (i.e. the local dimension of $(X, \tau)$ at $x$ is
at most one).

\begin{remark}\label{remark_curve_selection}
Let $\BB$ be a downward directed family of subsets of $R^n$ (i.e. for every pair $B, B'\in\BB$ there is $B''\in\BB$ such that $B''\subseteq B \cap B'$). We say that a curve $\gamma:(a,b)\rightarrow R^n$ with convergence endpoint $c\in\{a,b\}$ is cofinal in $\BB$ if, for every $B\in \BB$, $\gamma(t)\in B$ for all $t\in (a,b)$ close enough to $c$.  

Given this terminology, what the proof of Lemma~\ref{lemma_curve_selection_R} shows is that any definable downward directed family of non-empty sets of dimension at most one admits a definable cofinal curve.   
In~\cite{atw1}, the authors and Walsberg study definable topologies through an analysis of definable directed sets and existence of cofinal maps. We show that definable directed sets of all dimensions admit definable cofinal curves whenever $\RR$ expands an ordered field, but that this property does not hold for o-minimal structures in general.

These results can be used in characterizing definable compactness as described in Remark~\ref{rem:compactness}. Specifically, the proof of Lemma~\ref{lemma_curve_selection_R} can be expanded to show that, for one-dimensional definable topological spaces, definable compactness in the sense of Definition~\ref{dfn:compact} implies the definition given in condition (\ref{dfn:directed-compact}) (see~\cite[Lemma 6.2.1 and Proposition 6.2.4]{andujar_thesis}).
\end{remark}

Definable curve selection allows us to understand continuity in terms of convergence of definable curves. This was already shown in~\cite{atw1}, Proposition 42. We state here explicitly the analogous statement for one-dimensional spaces that we will use throughout this paper. Note, however, that the conclusion of the following statement (as well as that of~\cite{atw1}, Proposition 42) holds more generally, in that the proof does not specifically require that $(X,\tau)$ have dimension at most one (or be definable in an expansion of a field), only that it have definable curve selection.

\begin{proposition}\label{prop_cont_lim}
Let $(X,\tau)$ and $(Y,\mu)$ be definable topological spaces, where $\dim X \leq 1$. Let $f:(X,\tau) \rightarrow (Y,\mu)$ be a definable map. Then $f$ is continuous at $x\in X$ if and only if, for every definable curve $\gamma:(a,b) \rightarrow X$ and $c\in\{a,b\}$, if $\gamma$ $\tau$-converges to $x$ as $t$ tends to  $c$, then $f \circ \gamma$ $\mu$-converges to $f(x)$ as $t$ tends to $c$. 
\end{proposition}
\begin{proof}
The proof is identical to that of~\cite{atw1}, Proposition 42, except that, in order to invoke the property of definable curve selection here, the appeal to Proposition 41 therein should be replaced by one to Lemma~\ref{lemma_curve_selection_R} in this paper.
\end{proof} 

Definable curve selection also allows us to prove the following lemma, which we will make use of in proving our characterization (Theorem~\ref{them_main_2}) of which one-dimensional definable topological spaces are homeomorphic to euclidean space (see in particular Lemma~\ref{lemma_walsberg_homeomorphism}). 
By Remark~\ref{rem:compactness}, this is equivalently \cite[Lemma 3.11]{johnson14}, but we include this direct proof for the sake of completeness.

\begin{lemma}\label{lemma:compact_homeomorphism}
Let $f:(X,\tau)\rightarrow (Y,\mu)$ be definable continuous bijection between one-dimensional definable topological spaces. If $(X,\tau)$ is definably compact and $(Y,\mu)$ is Hausdorff, then $f$ is a homeomorphism. 
\end{lemma}
\begin{proof}
By Proposition~\ref{prop_cont_lim} it suffices to prove that, for any definable curve $\gamma$ in $Y$, if $\gamma$ $\mu$-converges to some $f(x)$, then $h^{-1}\circ \gamma$ $\tau$-converges to $x$. 

Let $\gamma$ be a definable curve in $Y$ $\mu$-converging to some $y\in Y$. By definable compactness of $(X,\tau)$, the curve $f^{-1}\circ \gamma$ $\tau$-converges to some $x\in X$. Then, by Proposition~\ref{prop_cont_lim} and continuity of $f$, the curve $\gamma=f\circ f^{-1}\circ \gamma$ $\mu$-converges to $f(x)$. Since $\mu$ is Hausdorff, any definable curve can converge to at most one point, so $f(x)=y$. This completes the proof.  
\end{proof}

Definable curve selection also allows us to prove the following facts regarding $e$-accumulation sets. For completeness and in accordance with the focus of this paper we only prove them for one-dimensional definable topological spaces. Nevertheless, one may show, using~\cite[Corollary 25]{atw1}, that Proposition~\ref{basic_facts_P_x_1}(\ref{itm: basic_facts_3}) holds for definable topological spaces of all dimensions whenever $\RR$ expands an ordered field and, using the fact that the space $(\exR^m,\tau_e)$ is definably compact for every $m$ and Remark~\ref{rem:compactness}, that Proposition~\ref{basic_facts_P_x_1}(\ref{itm: basic_facts_4}) holds in general for definable topological spaces of all dimensions. In the next proposition the euclidean closure means closure with respect to the euclidean topology on $\exR^m$.

\begin{proposition}\label{basic_facts_P_x_1}
Let $(X,\tau)$, $X \subseteq R^m$, $\dim X \leq 1$, be a definable topological space. 
\begin{enumerate}[(a)]
\item \label{itm: basic_facts_3} For any $x\in X$, $y\in \exR^m$, it holds that $y\in \E_x$ if and only if there exists an injective definable curve in $X$ $\tau$-converging to $x$ and $e$-converging to $y$.  
\item \label{itm: basic_facts_4} Let $Y\subseteq X$ be a definable set and $x\in \partial_\tau Y$. If $\tau$ is $T_1$, then $\E_x \cap cl_e Y \neq \emptyset$.     
\end{enumerate}
\end{proposition}
\begin{proof}
The right to left implication in~(\ref{itm: basic_facts_3}) is immediate. For the left to right implication, fix $x\in X$ and $y \in E_x$. Consider the definable topology $\mu$ on $X$ where every $z\neq x$ is isolated and where a basis of neighbourhoods of $x$ is given by the family $\{ \{x\} \cup (A \cap B \setminus \{y\}): x\in A\in \tau, \, y\in B\in \tau_e\}$. Clearly, $\mu$ is Hausdorff and finer than $\tau$. Since $y\in \E_x$, the sets $(A \cap B \setminus \{y\})$, where $x\in A\in \tau$ and $y\in B\in \tau_e$, are non-empty. 
In particular, they intersect $X \setminus \{x\}$: if not, then some such set $A \cap B \setminus \{y\}$ must be equal to $\{x\}$, so $y \neq x$, and there is an element $B' \in \tau_e$ with $B' \subseteq B$ that contains $y$ but does not contain $x$; in this case $x \notin A \cap B' \setminus \{y\} \subseteq A \cap B \setminus \{y\} = \{x\}$, so $A \cap B' \setminus \{y\} = \emptyset$, which is a contradication. Hence $x$ is in the $\mu$-closure of $X\setminus \{x\}$. 
Applying Lemma~\ref{lemma_curve_selection_R}, there necessarily exists an injective definable curve $\gamma$ in $X \setminus \{x\}$ $\mu$-converging (and thus $\tau$-converging) to $x$. By construction, $\gamma$ must $e$-converge to $y$.

To prove~(\ref{itm: basic_facts_4}), note that, if $x\in\partial_\tau Y$, then, by Lemma~\ref{lemma_curve_selection_R}, there is a definable curve $\gamma$ in $Y$ $\tau$-converging to $x$. If the topology is $T_1$, no such curve can be constant, and so, by o-minimality, $\gamma$ can be assumed to be injective. By o-minimality, $\gamma$ $e$-converges in $\exR^m$ and the result then follows from the right to left implication in~(\ref{itm: basic_facts_3}).    
\end{proof}

We now turn to the notion of $e$-accumulation set for definable topological spaces in the line. 
While for the rest of the section we deal almost exclusively with spaces in the line, recall (Remark~\ref{remark_assumption_X}) that, whenever $\RR$ expands an ordered field, any definable topological space of dimension at most one is definably homeomorphic to a space in the line. Furthermore, in Remark~\ref{remark_generalizing_Rx_Lx} below we describe how the definitions and results for spaces in the line that follow can be generalized to all spaces of dimension at most one, regardless of any assumption on $\RR$ beyond o-minimality.

\begin{lemma}\label{lemma:basic_facts_P_x_2} 
Let $(X,\tau)$, $X\subseteq R$, be a definable topological space. 
\begin{enumerate}[(a)]
\item\label{itma:lemma_basic_facts_P_x_2} Given $x, y\in X$, $y\in \E_x$ if and only if at least one of the following holds.
\begin{enumerate}[(i)]
\item \label{itm1:lemma_basic_facts_P_x_2} For any $\tau$-neighbourhood $A$ of $x$, there exists $z>y$ such that \mbox{$(y,z)\subseteq A$}. 
\item \label{itm2:lemma_basic_facts_P_x_2} For any $\tau$-neighbourhood $A$ of $x$, there exists $z<y$ such that \mbox{$(z,y)\subseteq A$}. 
\end{enumerate}
\item \label{itmb:lemma_basic_facts_P_x_2} It follows from (\ref{itma:lemma_basic_facts_P_x_2}) that, if $(X,\tau)$ is Hausdorff then, for any $y\in R$, there exist as most two points $x_0, x_1 \in X$ such that $y$ belongs in both $\E_{x_0}$ and $\E_{x_1}$ (i.e. for any distinct $x_0,x_1,x_2\in X$, $\E_{x_0} \cap \E_{x_1} \cap \E_{x_2}=\emptyset$).
\end{enumerate}
\end{lemma}
\begin{proof}

If~(\ref{itma:lemma_basic_facts_P_x_2})(\ref{itm1:lemma_basic_facts_P_x_2})   
fails then, by o-minimality, there exists a $\tau$-neighbourhood $A'$ of $x$ and $z'>y$ such that $(y,z')\cap A'=\emptyset$. Similarly if~(\ref{itma:lemma_basic_facts_P_x_2})(\ref{itm2:lemma_basic_facts_P_x_2})  
fails there is a $\tau$-neighbourhood $A''$ of $x$ and $z''<y$ with $(z'',y)\cap A''=\emptyset$. So $A'\cap A''$ is a $\tau$-neighbourhood of $x$ such that $(z'',z')\cap A'\cap A''\subseteq \{y\}$. This contradicts that $y\in \E_x$. The rest of the lemma is immediate.
\end{proof}

Lemma~\ref{lemma:basic_facts_P_x_2} motivates the following definition. 

\begin{definition}
Let $(X,\tau)$, $X \subseteq R$, be a definable topological space. For $x\in X$, we define the \emph{right $e$-accumulation set of $x$}, denoted $R_x \subseteq \E_x$, to be the set of points $y\in \exR$ satisfying that, for any $\tau$-neighbourhood $A$ of $x$, there exists $z>y$ such that $(y,z) \subseteq A$. In other words, if $\{A_u: u\in \Omega_x\}$ is a definable basis of $\tau$-neighbourhoods of $x$ in $(X,\tau)$, then
\[
R_x=\{y\in \exR : \forall u\in\Omega_x,\, \exists z>y, (y,z)\subseteq A_u\}.
\]
So the set $R_x\setminus \{-\infty\}$ is definable, for every $x \in X$. Similarly, the \emph{left $e$-accumulation set of $x$}, denoted $L_x \subseteq \E_x$, is defined to be the set of points $y\in \exR$ satisfying that, for any $\tau$-neighbourhood $A$ of $x$, there exists $z<y$ such that $(z,y) \subseteq A$. In other words, with $\{A_u: u\in \Omega_x\}$ as above,
\[
L_x=\{y\in \exR : \forall u\in\Omega_x,\, \exists z<y, (z,y)\subseteq A_u\},
\]
and $L_x \setminus \{ + \infty\}$ is likewise a definable set, for every $x \in X$.
\end{definition}

The following proposition follows from the definition of right and left $e$-accumulation set and Lemma~\ref{lemma:basic_facts_P_x_2}.

\begin{proposition}\label{prop_basic_facts_P_x_2}
Let $(X,\tau)$ be a definable topological space with $X\subseteq R$ and $x\in X$. Then
\begin{enumerate}[(a)]
\item\label{itm1:basic_facts_P_x_2} the relations $\{\al x,y\ar \in X\times \exR: y \in R_x\}$ and $\{\al x,y\ar \in X\times \exR: y\in L_x\}$ are definable;
\item\label{itm2:basic_facts_P_x_2} $\E_x= R_x \cup L_x$;
\item\label{itm3:basic_facts_P_x_2} if $(X,\tau)$ is Hausdorff then, for any $y\in X\setminus \{x\}$, $R_x \cap R_y = \emptyset$ and $L_x \cap L_y = \emptyset$. 
\end{enumerate} 
\end{proposition}
\begin{remark}\label{remark:R_L_finite}
By Lemma~\ref{lemma_2} and Proposition~\ref{prop_basic_facts_P_x_2}(\ref{itm2:basic_facts_P_x_2}), if $(X,\tau)$ is $T_1$, then $R_x$ and $L_x$ are finite for every $x\in X$. 
\end{remark}

\begin{remark}\label{remark_side_convergence}
Let $\gamma$ be a definable curve in $X \subseteq R$. By o-minimality, $\gamma$ $e$-converges to some $y\in \exR$. If $\gamma$ is injective, then we may assume that it lies in either $(y,+\infty)$ or $(-\infty,y)$ (recall Remark~\ref{remark_assumptions_curves}). In the former case, we say that $\gamma$ $e$-converges to $y$ from the right and in the latter that it does so from the left. Let $\tau$ be a definable topology on $X$. Note that, if $\gamma$ $e$-converges to $y$ from the right (respectively left) and $x\in X$, then $\gamma$ $\tau$-converges to $x$ if and only if $y\in R_x$ (respectively $y\in L_x$).  
\end{remark}

The following lemma will be useful in Section~\ref{section: affine}. It follows easily from Remark~\ref{remark_side_convergence}.

\begin{lemma} \label{lem:RL-compact}
A definable topological space $(X,\tau)$ with $X\subseteq R$ is definably compact if and only if, for every interval $(a,b)\subseteq X$, it holds that $[a,b) \subseteq \bigcup_{x\in X} R_x$ and $(a,b] \subseteq \bigcup_{x\in X} L_x$.
\end{lemma}

It turns out that, if $(X,\tau)$ is $T_1$, then, for any $x\in X$, the sets $R_x$ and $L_x$ characterize a definable basis of neighbourhoods for $x$. We show this in the next lemma.     

\begin{lemma}\label{lemma_explaining_neighbourhoods}
Let $(X,\tau)$, $X\subseteq R$, be a definable $T_1$ topological space. Let $x\in X$. By Remark~\ref{remark:R_L_finite}, the sets $R_x$ and $L_x$ are finite. Set $R_x:=\{y_1,\ldots, y_n\}$ and $L_x:=\{z_1,\ldots, z_m\}$. Define $\UU(x)$ to be the family of sets of the form
\[
\{x\} \cup \bigcup_{1\leq i \leq n} (y_i,y'_i) \cup \bigcup_{1 \leq j \leq m} (z'_j, z_j),
\] 
which is a family uniformly definable over $(y'_1, \ldots, y'_n, z'_1,\ldots, z'_m)\in R^{n+m}$, where $y_i<y'_i$ and $z'_j<z_j$. The definable family $\{ U\cap X : U\in \UU(x)\}$ is a basis of neighbourhoods of $x$ in $(X,\tau)$. 
  
In particular, in the case where $\RR$ expands an ordered group and $x$ has a bounded $\tau$-neighbourhood (implying $\E_x\cap \{-\infty, +\infty\}=\emptyset$), we may take $\UU(x)$ to be of the form 
\[
U(x,\varepsilon):=\{x\} \cup \bigcup_{y\in R_x} (y,y+\varepsilon) \cup \bigcup_{y\in L_x} (y-\varepsilon, y),
\]
for $\varepsilon>0$.

By passing to a subfamily if necessary, we may always assume that $\UU(x)$ is a family of subsets of $X$.  
\end{lemma}
\begin{proof}
Let $\UU(x)$ be as in the lemma. By definition of $R_x$ and $L_x$ it clearly holds that, for every $\tau$-neighbourhood $A$ of $x$, there exists $U\in\UU(x)$ such that $U \subseteq A \subseteq X$. It therefore remains to prove that all sets in $\UU(x)$ are $\tau$-neighbourhoods of $x$. 

Towards a contradiction, suppose that there exists $U\in\UU(x)$ that is not a $\tau$-neighbourhood of $x$. So $x\in \partial_\tau (X\setminus U)$. By Lemma~\ref{lemma_curve_selection_R}, there exists a definable curve $\gamma: I\rightarrow X\setminus U$ (which is necessarily injective, as $(X,\tau)$ is $T_1$) that $\tau$-converges to $x$ and that, by o-minimality, must $e$-converge to some $a\in \exR$. By Remark~\ref{remark_side_convergence}, if $\gamma$ $e$-converges from the right, then $a\in R_x$, and otherwise $a\in L_x$. Either way, by construction of $\UU(x)$, it follows that $\gamma(I) \cap U \neq \emptyset$, a contradiction. 
\end{proof}

From the above lemma, it follows that, if $(X,\tau)$ is a $T_1$ definable topological space with $X\subseteq R$, then a point $x\in X$ is $\tau$-isolated if and only if $\E_x=\emptyset$, and the identity map \mbox{$(X,\tau)\rightarrow (X,\tau_e)$} is continuous at $x\in X$ if and only if $\E_x\subseteq \{x\}$.  

The following lemma will be fundamental in proofs in later sections. 

\begin{lemma}\label{lemma_f_Rx_Lx}
Let $(X,\tau)$, $X\subseteq R$, be an infinite definable topological space. 
Let $f:I\subseteq X \rightarrow R$ be a function on an interval $I=(a,b)$, $a,b\in \exR$, such that, for every $x\in I$, $f(x)\in \E_x$. Suppose that $f$ is \econtinuous and strictly increasing (respectively decreasing). We extend $f$ to a function $[a,b] \rightarrow \exR$ by letting $f(a)=\elim_{x\rightarrow a} f(x)$ and $f(b)=\elim_{x\rightarrow b} f(x)$. 
For all $y\in X$, we have that
\begin{enumerate}[(a)]
\item \label{itma:lemma_f_Rx_Lx} for any $x\in [a,b)$, if $x\in R_y$, then $f(x)\in R_y$ (respectively $f(x)\in L_y$); 
\item \label{itmb:lemma_f_Rx_Lx}for any $x\in (a,b]$, if $x\in L_y$, then $f(x)\in L_y$ (respectively $f(x)\in R_y$). 
\end{enumerate} 

Under the additional assumption that $\tau$ is regular, the converse also holds. In other words,
\begin{enumerate}[(a)]
\setcounter{enumi}{2}
\item \label{itmc:lemma_f_Rx_Lx} for any $x\in [a,b)$, if $f(x)\in R_y$ (respectively $f(x)\in L_y$), then $x\in R_y$; 
\item \label{itmd:lemma_f_Rx_Lx} for any $x\in (a,b]$, if $f(x)\in L_y$ (respectively $f(x)\in R_y$), then $x\in L_y$. 
\end{enumerate} 
\end{lemma}
\begin{proof}
Let $y\in X$ and suppose that $f$ is strictly increasing. We prove (\ref{itma:lemma_f_Rx_Lx}) and (\ref{itmc:lemma_f_Rx_Lx}) in this case, with all other parts of the lemma proved analogously to one of these cases.

 We begin with case (\ref{itma:lemma_f_Rx_Lx}). Suppose that $x\in [a,b) \cap R_y$. If $f(x)\notin R_y$, then, by o-minimality, there is $z>f(x)$ and an open $\tau$-neighbourhood $A$ of $y$ such that $(f(x),z)\cap A =\emptyset$. Since $x\in R_y$, there is $x'>x$ in $I$ such that $(x,x')\subseteq A$. Since $f$ is \econtinuous and strictly increasing, there is $x'' \in (x,x')$ such that $f(x'')\in (f(x),z)$. So $A$ is a $\tau$-neighbourhood of $x''$ and $f(x'')\notin cl_e A$, which contradicts that $f(x'')\in \E_{x''}$. This completes the proof of (\ref{itma:lemma_f_Rx_Lx}) in the increasing case.

To prove~(\ref{itmc:lemma_f_Rx_Lx}), suppose that $f(x)\in R_y$, for $x\in [a,b)$, and let $A$ be a $\tau$-neighbourhood of $y$. Then there is some $z>f(x)$ such that $(f(x),z)\subseteq A$. Since $f$ is \econtinuous and strictly increasing, there is $x'>x$ such that $f[(x,x')]\subseteq (f(x),z)$. For any $x''\in (x,x')$, since $f(x'')\in \E_{x''}$ and $f(x'')\in(f(x),z)\subseteq A$, it follows that $x''\in cl_\tau(f(x),z) \subseteq  cl_\tau A$. Hence $(x,x')\subseteq cl_\tau A$. So we have shown that, for every $\tau$-neighbourhood $A$ of $y$, there exists $x'>x$ such that $(x,x')\subseteq cl_\tau A$. If $x\notin R_y$, then there must exist some $x''>x$ and some $\tau$-neighbourhood $A'$ of $y$ such that $(x,x'')\cap A'=\emptyset$. But then, by regularity, there is a $\tau$-neighbourhood $A''\subseteq A'$ of $y$ such that $cl_\tau A'' \subseteq A'$, and in particular such that $(x,x'') \cap cl_\tau A'' =\emptyset$, a contradiction by the above.
\end{proof} 

\begin{definition}\label{definition:dfbly_normal}
Let $(X,\tau)$ be a definable topological space. We say that $(X,\tau)$ is \emph{definably normal} if, given any pair of disjoint definable $\tau$-closed sets $B,C \subseteq X$, there exist definable disjoint open sets $U,V \subseteq X$ such that $B\subseteq U$ and $C\subseteq V$. 

We say that $(X,\tau)$ is \emph{definably completely normal} if any definable subspace of $(X,\tau)$ is definably normal. 
\end{definition}

\begin{proposition}\label{prop_reg_implies_normal}
Let $(X,\tau)$, $X\subseteq R$, be a definable topological space. If $(X,\tau)$ is $T_1$ and regular, then it is definably completely normal. 
\end{proposition}

\begin{proof}
We suppose that $(X,\tau)$, $X\subseteq R$, is $T_1$ and regular and prove that it is definably normal. Since being $T_1$ and regular are hereditary properties we conclude that $(X,\tau)$ is definably completely normal.

Let $B,C \subseteq X$ be disjoint $\tau$-closed definable sets in $(X,\tau)$. To prove the proposition it suffices to show the existence of a definable $\tau$-neighbourhood $U$ of $B$ such that the $\tau$-closure of $U$ is disjoint from $C$. We proceed by constructing a suitable partition of $B$ into two sets, $B=B'\cup B''$, where $B''$ is finite. By regularity of $(X,\tau)$, there clearly exists a definable $\tau$-neighbourhood $U''$ of $B''$ such that $cl_\tau U'' \cap C=\emptyset$. It is therefore enough to show, with $B'$ and $B''$ defined in this way, the existence of a definable $\tau$-neighbourhood $U'$ of $B'$ such that $cl_\tau U' \cap C=\emptyset$. The proof is then completed by taking $U=U'\cup U''$. 

First note that, since $(X,\tau)$ is $T_1$ and regular, the space is also Hausdorff. Set $E_B:=\bigcup_{x\in B} \E_x$. Let $int_e E_B$ be the euclidean interior of $E_B$ and set $B':=\{x\in B : \E_x \subseteq int_e E_B\}$. By o-minimality, $E_B\setminus int_e E_B$ is finite and so, by Hausdorffness and Lemma~\ref{lemma:basic_facts_P_x_2}(\ref{itmb:lemma_basic_facts_P_x_2}), $B''=B\setminus B'$ is also finite. Applying Lemma~\ref{lemma_explaining_neighbourhoods}, and using the fact that, by o-minimality, $int_e E_B$ is a finite union of intervals, observe that $B'\cup int_e E_B$ is a $\tau$-neighbourhood of $B'$.

Now note that, for any $x\in X$, if $\E_x \cap int_e C \neq \emptyset$ then, by definition of the set $\E_x$, it holds that $x\in cl_\tau C$, and so $x\in C$  (since $C$ is $\tau$-closed). Since $B\cap C=\emptyset$ it follows that $\E_B \cap int_e C =\emptyset$, i.e. $\E_B \cap C \subseteq C \setminus int_e C$, and in particular, by o-minimality, $\E_B \cap C$ is finite. It follows that $int_e \E_B \setminus C$ is cofinite in $int_e \E_B$. Now recall from the previous paragraph that $B'\cup int_e E_B$ is a $\tau$-neighbourhood of $B'$. 
Since $(X,\tau)$ is $T_1$ and $int_e E_B \setminus C$ is cofinite in $int_e B$, we conclude that $U'=B'\cup int_e E_B \setminus C$ is also a $\tau$-neighbourhood of $B'$. We complete the proof by showing that $cl_\tau U' \cap C=\emptyset$.

Towards a contradiction, suppose that some $x\in C$ is in the $\tau$-closure of $U'$. Then, since $B$ is $\tau$-closed and disjoint from $C$, $x$ must be in the $\tau$-closure of $int_e E_B \setminus C$. Set $E'_B:=int_e E_B \setminus C$. If there were some $\tau$-neighbourhood $A$ of $x$ such that $A\cap E'_B$ were finite, then, since the space is $T_1$, $A\setminus E'_B$ would also be a $\tau$-neighbourhood of $x$, which contradicts that $x$ is in the $\tau$-closure of $E'_B$. On the other hand, suppose that, for every $\tau$-neighbourhood $A$ of $x$, the intersection $A\cap E'_B$ is infinite. Then, for every such $A$, there exists an interval $I\subseteq A\cap E'_B$. 
Since, by definition of $\E'_B$, every point in $I$ lies in $E_y$, for some $y\in B$, we clearly have that there exists $y \in B$ with $E_y \cap I \neq \emptyset$, hence there exists $y \in B \cap cl_{\tau}I$, by definition of $E_y$, and in particular $y \in  B \cap cl_\tau A$. So, in this case, for every $\tau$-neighbourhood $A$ of $x$, it holds that $B \cap cl_\tau A \neq \emptyset$, which contradicts that $(X,\tau)$ is regular. 
\end{proof}
As we already indicated at the end of Subsection~\ref{subsection: prelim_results_alldim}, the following remark will allow us to generalize many results in this paper about spaces in the line to definable topological spaces of dimension at most one, even when $\R$ does not expand an ordered field.

\begin{remark}\label{remark_generalizing_Rx_Lx}
Let $(X,\tau)$, $X\subseteq R^n$, $n>1$, be a definable topological space, let $I=(a,b)\subseteq R$, $a,b\in\exR$, be an interval and let $f:I\rightarrow f(I)\subseteq X$ be an \ehomeomorphism (which we extend continuously to a function $f\colon [a,b]\rightarrow \exR^n$). Suppose moreover that $f(a)\neq f(b)$. Consider the definable total order $\prec$ on $cl_e f(I)=f([a,b])$ given by identifying $cl_e f(I)$ with $[a,b]$ through $f$, i.e. for every pair $x,y\in cl_e f(I)$, set $x\prec y$ if and only if $f^{-1}(x) < f^{-1}(y)$. Accordingly, for any $x\prec y$ in $cl_e f(I)$, let $(x,y)_\prec$ denote the corresponding interval with respect to $\prec$. 

By means of this identification, we may generalise the notion of right and left $e$-accumulation point to points in $cl_e f(I)$. That is, if $x\in X$ and $y\in cl_e f(I)$, then $y\in R_x$ if and only if, for every $\tau$-neighbourhood $U$ of $x$, there is $z\in f(I)$ such that $y \prec z$ and $(y,z)_{\prec}\subseteq U$. Similarly, $y\in L_x$ if and only if, for every $\tau$-neighbourhood $U$ of $x$, there is $z\in f(I)$ such that $z \prec y$ and $(x,y)_{\prec}\subseteq U$. Proposition~\ref{prop_basic_facts_P_x_2}~(\ref{itm1:basic_facts_P_x_2}) and~(\ref{itm3:basic_facts_P_x_2}) generalize to this setting.

Now suppose that $\dim X \leq 1$. By o-minimal cell decomposition, there is a finite definable partition $\XX$ of $X$ into cells such that, for every $C\in \XX$, there is a projection $\pi_C:C\rightarrow I_C\subseteq R$ that is an \ehomeomorphism onto a cell (i.e. a point or an open interval). Note that, for every one-dimensional cell $C \in \XX$, if we set $I_C = (a,b)$ and extend $(\pi_{C})^{-1}$ continuously (as in the first paragraph of this remark) to a function $(\pi_{C})^{-1}:[a,b] \to \exR^n$, then $(\pi_C)^{-1}(a) \neq (\pi_C)^{-1}(b)$. 
By passing to a pushforward of $(X,\tau)$ if necessary, we may assume that, for every distinct pair $C,C'\in\XX$, we have $cl_e C \cap cl_e C'= \emptyset$. Then, for any $C\in \XX$ such that $I_C$ is an interval, let $\prec_C$ be the order on $cl_e C$ given by identifying $cl_e C$ with $cl_e I_C$ through $\pi_C$ (as indicated above). Now let $\{n(C)<\omega : C\in\XX\}$ be an enumeration of the cells in $\XX$ and let $\prec$ be definable linear order on $cl_e X$ such that, for any $x\in cl_e C$ and $y\in cl_e C'$, where $C, C'\in \XX$, we have that $x\prec y$ if and only if $n(C)<n(C')$ or $n(C)=n(C')$ and $x\prec_C y$, that is, $\prec$ is the linear order induced by the lexicographic order given the push-forward $x\mapsto \{n(C)\}\times \pi_C(x)$ for $x\in cl_e C$.

Given this convention, the space $(X,\tau)$ behaves very much like a space in the line. The definitions of right and left $e$-accumulation set immediately generalise to points $x\in X$, by saying that $y\in cl_e C$ belongs in $R_x$ or $L_x$ if it does with respect to $\prec_C$. Note that, under this construction, the definitions of sets $R_x$ and $L_x$, for any $x\in X$, are dependent on the choice of cell decomposition $\XX$ of $X$. 

Under this correspondence, the statements and proofs of Lemma~\ref{lemma:basic_facts_P_x_2}, Proposition~\ref{prop_basic_facts_P_x_2}, Lemma~\ref{lem:RL-compact}, Lemma~\ref{lemma_explaining_neighbourhoods} and Proposition~\ref{prop_reg_implies_normal}  generalise to this setting. 
Moreover, suppose that, for any $C,C'\in \XX$ and partial function $f:C\rightarrow C'$ defined on an interval $(a,b)_{\prec_C}$, we consider that $f$ is increasing or decreasing to mean with respect to $\prec_C$ and $\prec_{C'}$. Then Lemma~\ref{lemma_f_Rx_Lx} and its proof also generalise to $(X,\tau)$.
\end{remark}

The following definition is not topological in flavour. Nevertheless, we introduce it as a natural property to consider when seeking to prove facts about definable topological spaces via inductive arguments.

\begin{definition} \label{dfn:fdi}
We say a definable topological space $(X,\tau)$ satisfies the \emph{frontier dimension inequality} (\textbf{fdi}) if, for every non-empty definable set $Y\subseteq X$, $\dim \partial_\tau Y < \dim Y$. 
\end{definition}
The \textbf{fdi} is clearly a hereditary property. The topologies $\tau_e$, $\tau_r$, $\tau_l$ and $\tauc$ all satisfy the \textbf{fdi}; however, we will show that any $T_1$ definable compactification of $(R,\tau)$, where $\tau\in\{\tau_r, \tau_l, \tauc\}$, does not (see the proof of Corollary~\ref{cor_them_2}). Observe that the \textbf{fdi} implies in particular that the frontier of any finite set is empty, and so any space with the \textbf{fdi} is $T_1$. Walsberg proved~\cite[Lemma 7.15]{walsberg15} that every definable metric space satisfies the \textbf{fdi}. By an inductive argument on dimension it is easy to show that in any space with this property every definable set is a boolean combination of definable open sets (i.e. property (A) in~\cite[Section 2]{pillay87} holds). 

The next proposition highlights a connection between the \textbf{fdi} and regularity.  

\begin{proposition}\label{prop_T2_frontier_ineq_regular}
Let $(X,\tau)$, $\dim X \leq 1$, be a Hausdorff definable topological space that satisfies the frontier dimension inequality. Then $(X,\tau)$ is regular.
\end{proposition}
\begin{proof}
We prove that, for any $x\in X$ and any $\tau$-neighbourhood $A$ of $x$, there exists a $\tau$-neighbourhood $U$ of $x$ such that $cl_\tau U\subseteq A$.

Let $x\in X$ and let $A$ be a $\tau$-neighbourhood of $x$. By passing to a subset of $A$ if necessary, we may assume that $A$ is definable. By the frontier dimension inequality, $\partial_\tau A$ is finite. Since $(X,\tau)$ is Hausdorff, there exists, for every $y\in \partial_\tau A$, a $\tau$-neighbourhood $A(y)$ of $x$ such that $y \notin cl_\tau A(y)$. Let 
\[
U= \bigcap_{y\in \partial_\tau A} A(y) \cap A.
\]
Then $U$ is a $\tau$-neighbourhood of $x$ and $cl_\tau U\subseteq A$.
\end{proof}

Recall that any $T_1$ regular topological space is Hausdorff. Since any definable topological space with the \textbf{fdi} is $T_1$, we derive from Proposition~\ref{prop_T2_frontier_ineq_regular} that a one-dimensional definable topological space with the \textbf{fdi} is regular if and only if it is Hausdorff.

The assumptions of Hausdorffness and the \textbf{fdi} in Proposition~\ref{prop_T2_frontier_ineq_regular} are justified respectively by Examples~\ref{example: T_1 fdi not regular} and~\ref{example: dfbly sep not hereditary} in Appendix~\ref{section:examples}, which describe non-regular definable topological spaces in the line, the first of which is non-Hausdorff and satisfies the \textbf{fdi}, and the second of which is Hausdorff but does not satisfy the \textbf{fdi}. 
In Appendix \ref{section:examples}, Example~\ref{example:fdi_Hausdorff_not_regular}, we construct a two-dimensional Hausdorff space which satisfies the \textbf{fdi} but again is not regular, showing that, in addition, Proposition~\ref{prop_T2_frontier_ineq_regular} does not generalize to spaces of dimension greater than one.

\section{$T_1$ and Hausdorff ($T_2$) spaces. Decomposition in terms of the $\tau_e$, $\tau_c$, $\tau_r$ and $\tau_l$ topologies}
\label{section:T1_T2_spaces}

This section focuses on the properties of $T_1$ and Hausdorff definable topological spaces of dimension one. 
The main results are Corollary~\ref{cor:general_them_main} and Theorem~\ref{them 1.5}. The first shows that every infinite $T_1$ definable topological space contains a definable copy of an interval with one of the $\tau_e$, $\tauc$, $\tau_r$ or $\tau_l$ topologies. We use this to give a positive answer to the Gruenhage 3-element basis conjecture of set-theoretic topology in our setting (see Subsection~\ref{subsection: 3-el_basis_conj}). Theorem~\ref{them 1.5} then improves Corollary~\ref{cor:general_them_main} in the setting of Hausdorff definable topological spaces, by showing that such spaces in the line can be definably partitioned into finitely many subspaces, each of which has one of the $\tau_e$, $\tauc$, $\tau_r$ or $\tau_l$ topologies, a result which immediately generalizes to all one-dimensional Hausdorff definable topological spaces (Corollary~\ref{cor 1.5}). 

We begin by approaching the first main result of the section only for spaces in the line.

\begin{proposition}\label{them_main}
Let $(X,\tau)$, $X\subseteq R$, be an infinite $T_1$ definable topological space. Then there exists an interval $J\subseteq X$ such that \mbox{$(J,\tau)=(J,\tau_\square)$}, where $\square$ is one of $e$, $r$, $l$ or $s$.
\end{proposition}

We prove this proposition below. First, however, we present a generalization which follows directly from Proposition~\ref{them_main} and the fact that, by o-minimal cell decomposition, every infinite definable set contains a one-dimensional cell.

\begin{corollary}\label{cor:general_them_main}
Every infinite $T_1$ definable topological space has a subspace that is definably homeomorphic to an interval with either the euclidean, right half-open interval, left half-open interval, or discrete topology. 
\end{corollary}

\begin{proof}[Proof of Proposition~\ref{them_main}]
By Lemma~\ref{lemma_2}, for each $x\in X$ the set $\E_x$ is finite. Suppose that there exist infinitely many points $x\in X$ satisfying $(x,\infty)\cap \E_x=\emptyset$. In that case, let $I\subseteq X$ be a bounded interval containing only such points and fix $C>I$. Otherwise, let $I'\subseteq X$ be an interval such that $(x,\infty)\cap E_x\neq \emptyset$, for every $x\in I'$, and consider the definable map $f$ on $I'$ taking each $x$ to the smallest $y>x$ such that $y\in \E_x$. The map $f$ satisfies $x<f(x)$ for all $x\in X$ and so, by o-minimality, after passing if necessary to a subinterval where $f$ is continuous and then applying continuity, there exists an interval $I\subseteq I'$ and $C>I$ such that, for all $x\in I$, $f(x)>C$. In either case, we have that, for all $x\in I$, $(x,C]\cap \E_x=\emptyset$. Similarly, we can isolate a bounded subinterval $J\subseteq I$ and some $c<J$ such that, for every $x\in J$, $[c, x)\cap \E_x=\emptyset$. Thus we have reached an interval $J$, and $c, C \in R$ with $c<J<C$, such that, for all $x\in J$, we have $[c,C]\cap \E_x \subseteq \{x\}$, and so in particular $cl_e J \cap \E_x \subseteq \{x\}$. 

For any $x\in J$, let $\UU(x)$ denote a family of $\tau$-neighbourhoods of $x$ as described in Lemma~\ref{lemma_explaining_neighbourhoods}. 
Note that, by construction of $J$, for any given $x\in J$ and $y<x<z$, there is $U\in\UU(x)$ such that, 
\begin{equation}\label{eqn_U}
U \cap J=\begin{cases}
(y, z) & \text{if } x\in R_x \cap L_x, \\
[x,z) & \text{if } x\in R_x \setminus L_x,\\
(y,x] & \text{if } x\in L_x \setminus R_x,\\
\{x\} & \text{if } x\notin R_x \cup L_x.
\end{cases}
\end{equation}

Recall that the families $\{R_x : x\in J\}$ and $\{L_x : x\in J\}$ are definable. Thus we may partition $J$ into four definable sets as follows: 
$$\begin{aligned}
J_1 &=\{ x\in J: x\in R_x \cap L_x\},              \\
J_2 &=\{ x\in J : x \in R_x, x\notin L_x\},        \\
J_3 &=\{ x\in J : x \notin R_x,  x\in L_x \},      \\
J_4 &= \{ x\in J : x\notin R_x \cup L_x \}.        \\
\end{aligned}
$$
By~\eqref{eqn_U} and the definitions of $R_x$ and $L_x$, the subspace topology on $J_1$ is $\tau_e$. Similarly, the subspace topologies on $J_2$, $J_3$ and $J_4$ are $\tau_r$, $\tau_l$ and $\tauc$, respectively. At least one of the four definable sets $J_1$, $J_2$, $J_2$ and $J_4$ must be infinite and thus contain an interval, and so the proposition follows.
\end{proof}

For a justification for the condition of $T_1$-ness in Proposition~\ref{them_main}, see Appendix \ref{section:examples}, Example~\ref{example: t_0 not t_1}, which describes a $T_0$ definable topological space in the line that fails to be $T_1$ and does not contain an interval with one of the $\tau_e$, $\tau_r$, $\tau_l$ or $\tauc$ topologies.

\subsection{3-element basis conjecture} \label{subsection: 3-el_basis_conj}
We now discuss how Corollary~\ref{cor:general_them_main} relates to an open conjecture of set-theoretic topology due to Gruenhage known as the 3-element basis conjecture.

Given a class of topological spaces $\mathcal{S}$, a basis for $\mathcal{S}$ is a subclass $\mathcal{F} \subseteq \mathcal{S}$ such that every member of $\mathcal{S}$ contains a subspace that is homeomorphic to one of the members of $\mathcal{F}$. When the class $\mathcal{S}$ consists of topological spaces definable in some structure, we will say that $\mathcal{F} \subseteq \mathcal{S}$ is a \emph{definable basis} for $\mathcal{S}$ if every member of $\mathcal{S}$ contains a subspace that is definably homeomorphic to one of the members of $\mathcal{F}$.

Consider the following statement.
\begin{equation}\label{star}
\parbox{0.9\textwidth}{The class $\mathcal{S}$ of uncountable, first-countable, regular, Hausdorff topological spaces has a  basis $\mathcal{F}$ consisting of an uncountable discrete subspace, a fixed uncountable subset of the reals with the euclidean topology, and a fixed uncountable subset of the reals with the Sorgenfrey topology.}
\end{equation} 

The \emph{$3$-element basis conjecture for uncountable, first-countable, regular, Hausdorff spaces} is an open conjecture stating that $\eqref{star}$ is consistent with ZFC. It arose in connection with various questions concerning perfectly normal, compact, Hausdorff spaces, in particular questions due to Fremlin and Gruenhage, which we discuss further in Subsection~\ref{subsection: Fremlin}. The $3$-element basis conjecture appears as Question 2 in~\cite{gm07}, following a series of works by Gruenhage in which related statements were put forward. In \cite{gru88} and \cite{gru89}, Gruenhage considered the statement that all uncountable, first-countable, regular, Hausdorff spaces contain either an uncountable metrizable subspace or a copy of an uncountable subspace of the Sorgenfrey Line. He proved under the Proper Forcing Axiom (PFA) that this statement holds for regular cometrizable spaces (spaces that admit a coarser metric topology such that each point has a neighbourhood base consisting of sets closed in the metric topology), by showing that all such spaces either contain an uncountable discrete subspace, contain a copy of an uncountable subset of the reals with the Sorgenfrey topology, or are cosmic (i.e. are the continuous image of a separable metric space). Todor\v{c}evi\'{c} then reproved this result under the Open Colouring Axiom (OCA), which follows from PFA~\cite{todor89}. In both cases, the assumption of first-countability was not in fact required. 

Following these results, Gruenhage put forward the question of the consistency (with ZFC) of the version of \eqref{star} in which the first-countability assumption is not imposed \cite{gru90}. This stronger statement was shown to fail by Moore in his solution to the $L$-space problem~\cite{moore06}, and so the first-countability assumption in \eqref{star} is necessary if the $3$-element basis conjecture is to have a positive answer. (Note that it is known that $L$-spaces, which are hereditarily Lindel\"of non-separable spaces, as well as $S$-spaces, which are hereditarily separable non-Lindel\"of spaces, provide counterexamples to \eqref{star}; these can be constructed, for example, with the continuum hypothesis (CH), but the existence of $S$-spaces is independent of ZFC \cite{todor83} and, to the best of our knowledge, to date no \emph{first-countable} $L$-space has been constructed in ZFC. See \cite{gru90} for further discussion.) 
Beyond the statements above, we note that Todor\v{c}evi\'{c} studied related questions in the context of the class of spaces that can be represented as relatively compact subsets of the class of all Baire class 1 functions on a Polish space, endowed with the topology of pointwise convergence~\cite{todor99}. In addition, the conclusion of the conjecture 
was shown by Farhat to be consistent under PFA for the class of uncountable subspaces of monotonically normal compacta, and under Souslin's Hypothesis (SH) for any uncountable space having a zero-dimensional monotonically normal compactification \cite{farhat15}. Recently, Peng and Todor\v{c}evi\'{c} gave an analysis of different possible approaches to proving or disproving the conjecture~\cite{peng_todor22}.

Here, we indicate how our work provides a number of results related to the 3-element basis conjecture in the context of topological spaces definable in o-minimal structures. In particular, Corollary~\ref{cor:general_them_main}, shown above, establishes the existence of a $3$-element basis, as is posited by (\ref{star}), for the class of infinite $T_1$ topological spaces definable in any o-minimal expansion of $( \mathbb{R}, < )$. Moreover, in such a basis, the fixed subsets of the reals that form the underlying sets of the basis elements (i.e. those which have either the discrete, euclidean or Sorgenfrey topology) can each be taken to be $\mathbb{R}$ itself.

\sloppy More specifically, Corollary~\ref{cor:general_them_main} states that, for any o-minimal structure $\RR=(R,<,\ldots)$, the collection of definable topological spaces \mbox{$\{ (I,\tau) : I \subseteq R \text{ is an interval}, \, \tau\in\{\tau_e, \tau_r, \tau_l, \tauc\} \}$} is a definable basis for the class of infinite $T_1$ topological spaces definable in $\RR$. Clearly, if $\RR$ defines an order-reversing bijection (e.g. when $\RR$ expands an ordered group), then this definable basis can be reduced to the family \mbox{$\{ (I,\tau) : I \subseteq R \text{ is an interval}, \, \tau\in\{\tau_e, \tau_r,  \tauc\} \}$}. In addition, whenever $\RR$ expands an ordered field, this can in fact be reduced further to a definable basis consisting only of the three fixed spaces $(R,\tau_e)$, $(R,\tau_r)$ and $(R, \tauc)$, since, in this case, any interval with one of the $\tau_e$, $\tau_r$, or $\tauc$ topologies is definably homeomorphic to one of these three spaces. If there is no requirement that the basis be definable (in the sense that the homeomorphisms involved are definable), then, as long as the ordered set $(R,<)$ can be expanded to an ordered field $(R, +, \cdot, <)$, we have that $\{ (R,\tau_e), (R,\tau_r), (R, \tauc)\}$ serves as a 3-element basis, in the sense of \eqref{star}, for the class of infinite $T_1$ topological spaces definable in any o-minimal expansion of $(R, < )$ (and in particular this holds, as indicated above, in the special case that $(R,<)=(\mathbb{R},<)$).

Note that, for none of these results do we require the assumptions of uncountability, first-countability, Hausdorffness or regularity (although it should be noted that, in the special case of o-minimal expansions of $( \mathbb{R}, < )$, every definable topological space is known to be first-countable; see \cite[Proposition 38]{atw1}).

\subsection{Non-definability of the Cantor space} \label{subsection: Cantor_space}

We now consider the Cantor Space $2^\omega$ (i.e. the product of countably many copies of the discrete space $\{0,1\}$). Whether or not such a space exists (up to homeomorphism) in the form of a definable topological space in a given structure could be considered a means of assessing the tameness of the structure. We make use of Proposition~\ref{them_main} to prove that, in the setting of o-minimal structures, no space homeomorphic to the Cantor Space is definable.

Note that there is a notion of a `Cantor set', as defined in \cite{fkms-10}, which has been studied by various authors in considering notions of tameness expanding o-minimality. A Cantor set in this sense is (in the terminology of the present paper) any non-empty compact subset of $(\mathbb{R}, \tau_e)$ with empty interior and no isolated points. The Cantor space is known to be homeomorphic to such a set (for example, to the classical middle thirds Cantor set with the euclidean topology). On the one hand, in any o-minimal expansion of $(\mathbb{R},<)$, such a Cantor set is clearly not definable. However, on the other hand, there are examples of definable topological spaces in this structure which do have non-empty definable subsets that are compact, have empty interior and do not have isolated points (for example, the subset $[0,1] \times \{0\}$ of the definable Alexandrov $2$-line, which will be defined in Section~\ref{section: universal spaces}, Definition~\ref{example_n_line_0}). We show that, nevertheless, the Cantor space is not homeomorphic to any definable topological space in an o-minimal structure.

The proof of this fact will follow from a series of auxiliary results, beginning with two concerning the weights of our key definable topological spaces. 
Recall the classical definition of the \emph{weight} of a topological space $(X,\tau)$, $w_\tau(X)$, namely the minimum cardinality of a basis for $\tau$.

\begin{lemma} \label{lem:weight-taur}
For any set $X\subseteq R$ it holds that $w_r(X)=w_l(X)=w_s(X)=|X|$. 
Furthermore, in the case where $X$ is infinite and definable and $(R,<)$ can be expanded to an ordered field, it holds that $w_r(X)=w_l(X)=w_s(X)=|R|$.
\end{lemma}
\begin{proof}
The second sentence of the lemma follows from the first one by applying o-minimality and the fact that any interval in an ordered field is in bijection with the whole field. We prove the first sentence.

In the case of the discrete topology $\tauc$ the statement $w_s(X)=|X|$ is immediate from the definition. We show that $w_{r}(X)=|X|$. (An analogous argument shows that $w_{l}(X)=|X|$.) We may assume that $X$ is infinite, since otherwise the topology is discrete. By the definition of $\tau_r$ (Appendix \ref{section:examples}, Example~\ref{example:taur_taul}) we have that $w_{r}(X)\leq |X|^2 = |X|$. We show the reverse inequality.
Let $\BB$ be a basis of $(X,\tau_r)$ of minimum cardinality. For every $x\in X$ there exists some $A\in\BB$ such that $x\in A \subseteq [x,+\infty)$. A map $X\rightarrow \BB$ that takes each $x\in X$ to one such neighbourhood $A$ must be injective, so $|X|\leq |\BB|= w_r(X)$. 
\end{proof}

\begin{proposition} \label{prop:weight}
Let $(X,\tau)$ be an infinite $T_1$ definable topological space. 
Let $\alpha_R:=\min\{ w_e(I) : I\subseteq R \text{ is an interval}\}$. The following hold. 
\begin{enumerate}[(a)]
    \item $\alpha_R\leq w_\tau(X) \leq |X|\leq 2^{w_e(X)}$.\label{itm_weight_1}
    \item If $\RR$ expands an ordered field, then $w_e(R)\leq w_\tau(X)$. \label{itm_weight_2}
\end{enumerate}
\end{proposition}

\begin{proof}
The inequality $w_\tau(X) \leq |X|$ is given by~\cite[Corollary 2]{ag_cardinality}. The inequality $|X|\leq 2^{w_e(X)}$ follows from noticing that, for any basis for the euclidean topology on $X$, due to this topology being $T_1$, the map taking each $x\in X$ to the collection of its basic neighbourhoods is injective.

To prove that $\alpha_R\leq w_\tau(X)$, note that, for any topological space $(Z,\mu)$ and subspace $Z'\subseteq Z$, it holds that $w_\mu(Z') \leq w_\mu(Z)$. Now, by Corollary~\ref{cor:general_them_main}, there exists an interval $I\subseteq R$ and a definable embedding $(I,\mu)\hookrightarrow (X,\tau)$, where $\mu$ is one of $\tau_e$, $\tau_r$, $\tau_l$ or $\tauc$. In particular we have that $w_\mu(I)\leq w_\tau(X)$. Furthermore by Lemma~\ref{lem:weight-taur} observe that $\alpha_R \leq w_\mu(I)$, and so we derive that $\alpha_R \leq w_\tau(X)$.

Statement (\ref{itm_weight_2}) follows from the inequality $\alpha_R\leq w_\tau(X)$ in (\ref{itm_weight_1}) and from noticing that, if $\RR$ expands an ordered field, then any two intervals are definably $e$-homeomorphic, and so $\alpha_R = w_e(R)$.
\end{proof}

\begin{lemma}\label{lemma_compact_implies_reals}
Let $(X,\tau)$ be a $T_1$ definable topological space, let $I\subseteq R$ be an interval, and suppose that there exists an injective definable curve $\gamma:I \to X$. If $(X,\tau)$ is compact, then the linear order $(I,<)$ is Dedekind complete (i.e. every non-empty subset of $I$ that is bounded above in $I$ admits a supremum).
\end{lemma} 
\begin{proof}
Towards a contradiction suppose that there exists a non-empty set $S\subseteq I$ bounded above in $I$ but with no supremum. Let $S'$ be the set of upper bounds of $S$ in $I$. Consider the following family of closed non-empty subsets of $(X,\tau)$:
\[
\SSS= \{ cl_\tau(\gamma[(t,s)]) : t\in S, s\in S'\}.
\]  
This family clearly has the finite intersection property. By compactness of $(X,\tau)$, there exists $x\in X$ belonging in $\bigcap \SSS$. We define a curve $\gamma_x \colon I \to X$ as follows. If $x\notin \gamma(I)$, then we fix some $t_x\in I$ and set 
\[
\gamma_{x}(t) =  \begin{cases} \gamma(t) &\text{ if } t \neq t_x, \\ x & \text{ if } t=t_x. \end{cases} 
\]
Otherwise, let $\gamma_x=\gamma$. Let $(I,\tau_x)$ be the pull-back of $(\gamma_x(I),\tau)$ by $\gamma_x$. 

Now note that, for every $t \in S$, every $s \in S'$ and every $\tau_x$-neighbourhood $U$ of $t_x$, we have  $U\cap (t,s) \neq \emptyset$. This is clear if $t_x$ itself also lies in $(t,s)$. 
If $t_{x} \notin (t,s)$, then $\gamma = \gamma_{x}$ on $(t,s)$, and so $x \in cl_{\tau}(\gamma[(t,s)]) = cl_{\tau}(\gamma_{x}[(t,s)])$, whence $\gamma_{x}(U) \cap \gamma_{x}[(t,s)] \neq \emptyset$, as $\gamma_{x}(U)$ is a $\tau$-neighbourhood of $x$ in $(\gamma_x(I),\tau)$; it follows that $U \cap (t,s) \neq \emptyset$, by injectivity of $\gamma_{x}$. 

Now, clearly $(I,\tau_x)$ is $T_1$, and so, by Lemma~\ref{lemma_2}, the set $E^{(I,\tau_x)}_{t_{x}}$ is finite. Since $S$ has no supremum in $I$, there exists $t\in S$ with \mbox{$t > \{ y \in E^{(I,\tau_x)}_{t_{x}} : \exists z \in S \text{ with } z > y\}$} and $s\in S'$ with $s<S'\cap E^{(I,\tau_x)}_{t_{x}}$. We may also choose them so that $t_x\notin (t,s)$. The fact that every $\tau_x$-neighbourhood of $t_x$ intersects the interval $(t,s)$ clearly contradicts Lemma~\ref{lemma_explaining_neighbourhoods}.  
\end{proof}

\begin{remark}\label{remark_compact_implies_reals}
Since any two intervals are order-isomorphic in an ordered field, it follows from Lemma~\ref{lemma_compact_implies_reals}  that, under the assumption that $\RR$ expands an ordered field, if there exists an infinite compact $T_1$ topological space definable in $\RR$, then $(R,<)$ is Dedekind complete. In particular, since $(\R,+,\cdot,<)$ is, up to (unique) field isomorphism, the only Dedekind complete ordered field, it must be that $\RR$ is an expansion of the field of reals. 
\end{remark}

\begin{proposition}\label{prop_no_cantor_set}
There exists no infinite definable topological space $(X,\tau)$ that is compact, totally disconnected, and that satisfies $w_\tau(X)<|Y|$ for every $Y\subseteq X$ that is infinite and definable. 
\end{proposition}

\begin{proof}
Let $(X,\tau)$ be an infinite compact totally disconnected definable topological space satisfying $w_\tau(X)<|Y|$ for every $Y\subseteq X$ that is infinite and definable. We reach a contradiction by showing that $(X,\tau)$ contains an infinite (and in fact definable) connected subspace.

First note that, since $(X,\tau)$ is totally disconnected, in particular it is $T_1$. By o-minimality, there exists an interval $I\subseteq R$ and  an injective definable curve $\gamma:I \rightarrow X$. 
Let $(I,\mu)$ be the pull-back of $(\gamma(I),\tau)$ by $\gamma$. Since $(I,\mu)$ is $T_1$, by Proposition~\ref{them_main} we may assume, after passing to a subinterval if necessary, that $\mu \in \{\tau_e, \tau_r, \tau_l, \tauc\}$. Since, by hypothesis, $w_{\mu}(I)=w_\tau(\gamma(I))\leq w_\tau(X) < |f(I)|=|I|$, by Lemma~\ref{lem:weight-taur} it must be that $\mu=\tau_e$. Now recall that, by compactness and Lemma~\ref{lemma_compact_implies_reals}, $(I,<)$ is Dedekind complete, and so $(I,\tau_e)$ is connected. Hence $(\gamma(I),\tau)$ is connected.
\end{proof}

\begin{corollary}\label{cor:no_cantor_set}
The Cantor space $2^\omega$ is not a definable topological space.  
\end{corollary}
\begin{proof}
Let $(X,\tau)$ be (homeomorphic to) the Cantor space and towards a contradiction suppose that it is a definable topological space. Let $\mathfrak{c}=|X|$ denote the cardinality of the continuum. 
Recall that the Cantor space is compact, totally disconnected (in particular $T_1$) and second countable (i.e. $w_\tau(X)\leq\omega$). Let $Y\subseteq X$ be an infinite definable set. We show that $|Y|=\mathfrak{c}$. The result then follows from Proposition~\ref{prop_no_cantor_set}.

Clearly $|Y|\leq \mathfrak{c}$. By o-minimality, there exists an interval $I\subseteq R$ and an injective definable curve $\gamma:I\rightarrow Y \subseteq X$. 
Since $(X,\tau)$ is compact and $T_1$, by Lemma~\ref{lemma_compact_implies_reals} we have that $(I,<)$ is Dedekind complete, so $\mathfrak{c}\leq |I|$, and hence $\mathfrak{c}\leq |\gamma(I)|\leq |Y|$.
\end{proof}

It follows from the above corollary that the class of definable topological spaces up to homeomorphism is not closed under countable products. 

\subsection{Hausdorff ($T_2$) spaces} \label{subsection: T2_spaces}
We now prove a strengthening of Proposition~\ref{them_main} for Hausdorff spaces. The assumption of Hausdorffness is necessary here; see Appendix \ref{section:examples}, Example~\ref{example: t_1 not t_2} for a $T_1$, non-Hausdorff space that cannot be decomposed as described in the statement.

\begin{theorem}\label{them 1.5}
Let $(X,\tau)$, $X\subseteq R$, be a Hausdorff definable topological space. Then there exists a finite partition $\XX$ of $X$ into points and intervals such that, for every $Y\in\XX$, $\tau|_Y \in \{\tau_e, \tau_r, \tau_l, \tauc\}$. 
\end{theorem}
\begin{proof}
We start by proving a simple case. Suppose that, for every $x\in X$, $\E^{(X,\tau)}_x\subseteq \{x\}$. We call this condition \cond. 
Then let us partition $X$ into four definable sets as follows. 

 \begin{align*}
&\{ x\in X: x\in R_x \cap L_x\},              \\
&\{ x\in X: x \in R_x, x\notin L_x\},        \\
&\{ x\in X: x \notin R_x,  x\in L_x \},      \\
&\{ x\in X: x\notin R_x \cup L_x \}.         \\ 
\end{align*}

By Lemma~\ref{lemma_explaining_neighbourhoods}, these correspond respectively to spaces with the $\tau_e$, $\tau_r$, $\tau_l$ and $\tauc$ topologies. By o-minimality, we can partition each of these into a finite number of points and intervals, and the result follows. 

In order to prove the theorem, it is enough to show that we may partition $(X,\tau)$ into finitely many definable subspaces where \cond\ holds. We do so as follows.  

Note that, for any definable subspace $S\subseteq X$ and any $x\in S$, $E^{(S,\tau)}_x \subseteq \E^{(X,\tau)}_x$. We prove the existence of a finite partition of $X$ formed by points and intervals such that, for any interval $I$ in the partition and any $x\in I$, $\E^{(X,\tau)}_x \cap cl_e I \subseteq \{x\}$. Since any element in $\E^{(I,\tau)}_x$ must belong in $cl_e I$ (Proposition~\ref{basic_facts_P_x}(\ref{itm: basic_facts_1})), it follows that, for any $x\in I$, $\E^{(I,\tau)}_x = \E^{(I,\tau)}_x \cap cl_e I \subseteq \E^{(X,\tau)}_x \cap cl_e I \subseteq \{x\}$, i.e. \cond\ holds in $(I,\tau)$, which completes the proof. 
 
From now on, for any $x\in X$, let $\E_x=\E^{(X,\tau)}_x$.
By Lemma~\ref{lemma_2}, for any $x\in X$, the set $\E_x$ is finite. By o-minimality (uniform finiteness), there exists some $n$ such that $|\E_x|\leq n$, for every $x\in X$. We may partition $X$ into finitely many definable subspaces $X_0,\ldots, X_n$, where $X_i=\{x\in X : |\E_x|=i\}$, for $0\leq i \leq n$. We fix $Y=X_m$, for some $0\leq m \leq n$ and prove the existence of a partition of $Y$ with the desired properties. Since otherwise the result is trivial we assume that $m>0$ and that $Y$ is infinite.

For $1\leq i \leq m$, let $f_i:Y\rightarrow \exR$ be the definable function taking each element in $x\in Y$ to the $i$-th smallest element in $\E_x$. Since the family $\{\E_x : x\in Y\}$ is definable, these maps are definable. Moreover, by Hausdorffness (see Lemma~\ref{lemma:basic_facts_P_x_2}(\ref{itma:lemma_basic_facts_P_x_2})(\ref{itm2:lemma_basic_facts_P_x_2}), these functions cannot be constant on any interval. By o-minimality, let $\YY$ be a partition of $Y$ into finitely many intervals and points such that, for every interval $I\in \YY$, the functions $f_i$, $1\leq i\leq m$, are \econtinuous and strictly monotone. 

Without loss of generality, we fix an interval $I\in \YY$ and show that, for any $x\in I$, $\E_x \cap cl_e I\subseteq \{x\}$, completing the proof. Let $x\in I$ and $y\neq x$ be such that $y\in \E_x \cap cl_e I$. If $y\in I$, then, by Lemma~\ref{lemma_f_Rx_Lx}, $\E_y\subseteq \E_x$. Since $|\E_y|=|\E_x|$, it follows that $\E_y =\E_x$, contradicting that the functions $f_i$ are injective. Suppose now that $y\in \partial_e I$. As $y \in E_x$, we have that $y=f_i(x)$ for some $1\leq i \leq m$. By $e$-continuity and strict monotonicity of $f_i$ on $I$ there exists a point $x'\in I$ such that $f_i(x')\in I$ and $f_i(x')\neq x'$. A contradiction then follows as before.
\end{proof}

By o-minimal cell decomposition, the above theorem can be immediately generalized to all one-dimensional spaces. 

\begin{corollary} \label{cor 1.5}
Let $(X,\tau)$, $\dim{X} \leq 1$, be a Hausdorff definable topological space. Then there exists a finite definable partition $\XX$ of $X$ such that, for every $Y\in\XX$, $(Y, \tau)$ is definably homeomorphic to a point or an interval with either the euclidean, discrete, right half-open interval or left half-open interval topology. 
\end{corollary}

Since the topologies $\tau_r$, $\tau_l$ and $\tauc$ are all finer that $\tau_e$, it follows by o-minimality that any definable function $(X,\tau)\rightarrow (R,\tau_e)$, where $\dim X\leq 1$ and $\tau$ is a Hausdorff definable topology, is cell-wise continuous.

We end this section with a statement noting that, for spaces in the line, having an interval subspace with any one of the euclidean, discrete or half-open interval topologies is a definable topological invariant. 
This is an easy consequence of the monotonicity theorem of o-minimality, and in particular the observation that the push-forward of an interval with the $\tau_r$ or $\tau_l$ topology by an \econtinuous strictly monotone definable function is an interval with either the $\tau_r$ or $\tau_l$ topology. 
This holds in weakly o-minimal structures too, since these have a form of monotonicity (see~\cite{arefiev97}).

\begin{lemma}\label{remark_classification_spaces_line}
If $(X,\tau)$ and $(Y,\mu)$, where $X,Y\subseteq R$, are definable topological spaces, and $f:(X,\tau)\rightarrow (Y,\mu)$ is a definable homeomorphism, then 
\begin{enumerate}[(i)]
\item \label{itm:tau-copies-1} if $(X,\tau)$ contains an interval subspace with the discrete topology, then $(Y,\mu)$ contains an interval subspace with the discrete topology;
\item if $(X,\tau)$ contains an interval subspace with the right half-open or left half-open interval topology, then  $(Y,\mu)$ contains an interval subspace with the right half-open or left half-open interval topology;
\item \label{itm:tau-copies-3} if $(X,\tau)$ contains an interval subspace with the euclidean topology, then $(Y,\mu)$ contains an interval subspace with the euclidean topology. 
\end{enumerate}
\end{lemma}
Note that (\ref{itm:tau-copies-1}) and (\ref{itm:tau-copies-3}) hold for spaces of all dimensions if we substitute \textquotedblleft interval subspace" with \textquotedblleft definable subspace of dimension $n$".  

Hence definable topological spaces in the line can be classified up to definable homeomorphism according to whether or not they contain interval subspaces with the euclidean, discrete or half-open interval topologies. Moreover, by Proposition~\ref{them_main}, every infinite $T_1$ space will fall into at least one of these categories.

\section{Hausdorff regular ($T_3$) spaces. Decomposition in terms of the $\taulex$ and $\taualex$ topologies.}\label{section: universal spaces}

In this section, we study Hausdorff regular (i.e. $T_3$) spaces in the line. The main result is Theorem~\ref{them_ADC_or_DOTS}, which states that any such space can be partitioned into a finite set and two definable open subspaces, one of which definably embeds into a space with the lexicographic order topology, and the other into a space which we label the definable Alexandrov $n$-line. In the next section, we use this result and its proof, as well as Theorem~\ref{them 1.5}, to address universality questions in our setting. In Section~\ref{section:compactifications}, we extend Theorem~\ref{them_ADC_or_DOTS} to show that any $T_3$ definable topological space in the line has a definable compactification. We also combine Theorem~\ref{them_ADC_or_DOTS} and its proof with our affineness result (Theorem~\ref{them_main_2}) in Section~\ref{section: affine} to address questions of Fremlin and Gruenhage on perfectly normal, compact, Hausdorff spaces (see Subsection~\ref{subsection: Fremlin}). 

We start by introducing the relevant topologies. 
Given $X\subseteq R^n$, we denote by $\lex$ the lexicographic order on $X$ and by $(X,\taulex)$ the topological space induced by $\lex$ on $X$. Clearly this space is definable whenever $X$ is.

\begin{definition}[Definable $n$-split interval, Appendix~\ref{section:examples}, Example~\ref{example:n-split}] \label{dfn:lex}

We call the space $(R\times\{0,\ldots, n-1\}, \taulex)$ the \emph{definable $n$-split interval}. This space has the property that all the points in $R\times\{i\}$, for $0<i<n-1$, are isolated. In the case that $n>1$, for any $x\in R$ a basis of open neighbourhoods of $\al x,0\ar$ is given by sets of the form 
\[
\al x,0\ar \cup ((y,x)\times\{0,\ldots, n-1\}) \text{ for } y<x,
\]
and a basis of open neighbourhoods of $\al x, n-1\ar$ is given by sets of the form 
\[
\al x,n\ar \cup ((x,y)\times\{0,\ldots, n-1\}) \text{ for } y>x.
\]
If $n=1$, then $(R\times\{0\}, \taulex)=(R\times \{0\}, \tau_e)$.
\end{definition}

\begin{definition}[Definable Alexandrov $n$-line, Appendix~\ref{section:examples}, Example~\ref{example_n_line}]\label{example_n_line_0}
Let $\taualex$ be the topology on $R^2$ where all points in $R^2\setminus R\times\{0\}$ are isolated and, for any $x\in R$, a basis of open neighbourhoods of $\al x, 0\ar$ is given by sets of the form 
\[
\{\al x,0 \ar\} \cup (((z,y)\setminus \{x\}) \times R), \text{ for } z<x<y.  
\]
Then, for any $n>0$, we call the space $(R\times\{0,\ldots, n-1\},\taualex)$ the \emph{definable Alexandrov $n$-line}. 
\end{definition}

Note, in particular, that $(R\times\{0\},\taulex)=(R\times\{0\},\taualex)=(R\times\{0\},\tau_e)$. 

We may now state the main theorem of this section.

\begin{theorem}\label{them_ADC_or_DOTS}
Let $(X,\tau)$, $X\subseteq R$, be a regular and Hausdorff definable topological space. Then there exist disjoint definable open sets $Y, Z \subseteq X$ with \mbox{$X\setminus (Y\cup Z)$} finite, and $n_Y, n_Z>0$, such that the following holds. 
\begin{enumerate}
\item \mbox{There exists a definable embedding $h_Y:(Y,\tau)\hookrightarrow (R \times \{0,\ldots, n_Y\}, \taulex)$.}
\item \mbox{There exists a definable embedding $h_Z:(Z,\tau)\hookrightarrow(R \times \{0,\ldots, n_Z\},\taualex)$.} 
\end{enumerate}
\end{theorem}   

Lemmas~\ref{remark_for_two_theorems} and~\ref{lemma_proof_two_thems} are the bulk of the proof of Theorem~\ref{them_ADC_or_DOTS}. They are also used in Section~\ref{section:compactifications} to prove that all regular Hausdorff definable topological spaces in the line can be definably Hausdorff compactified (Theorem~\ref{them_compactification}). 
In Lemma~\ref{remark_for_two_theorems}, we construct a finite family $\Xinf$ of pairwise disjoint definable open subsets of $X$ of a very special form such that $X\setminus \bigcup\Xinf$ is finite. 
In Lemma~\ref{lemma_proof_two_thems}, we construct, for every $A\in\Xinf$, a set $A^*$ of the form $I_A\times\{0,\ldots, n_A\}$, for some interval $I_A$ and natural number $n_A$, and a definable embedding $h_A:(A,\tau)\hookrightarrow (A^*,\tau_A)$, where $\tau_A$ is either $\taulex$ or $\taualex$. The construction will be such that $I_{A}\cap I_{A'}=\emptyset$ for distinct $A, A'\in\Xinf$. Then $Z$ will be the union of all the sets $A$ in $\Xinf$ such that $(A^*,\tau_A)=(A^*,\taualex)$, and $h_Z$ the union of the respective embeddings $h_A$. The set $Y$ and embedding $h_Y$ are constructed similarly from the remaining sets in $\Xinf$.
 
Until the end of the proof of Theorem~\ref{them_ADC_or_DOTS}, we fix a definable topological space $(X,\tau)$ with $X\subseteq R$.   

We introduce an equivalence relation on $X$ induced by the topology $\tau$ that is defined as follows. Given $x,y\in X$, we write $x\ssim y$ when one of the following holds:
\begin{enumerate}[(i)]
\item $x=y$; 
\item there exists some $z\in X$ such that $\{x,y\}\cap \E_z \neq \emptyset$ and, for all $z\in X, \, x\in \E_z \Leftrightarrow y\in \E_z$.
\end{enumerate}
This relation is clearly reflexive and symmetric, and one easily checks that it is transitive. Moreover, by Proposition~\ref{basic_facts_P_x}(\ref{itm: basic_facts_2}), it is definable. For any $x\in X$, we denote by $[x]$ the equivalence class $\{y\in X: y\ssim x\}$. We prove some preliminary facts regarding this relation. 

\begin{lemma}\label{lemma_finite_classes}
If $(X,\tau)$ is $T_1$, then every equivalence class of $\ssim$ is finite. 
\end{lemma}
\begin{proof}
Let $x\in X$. If $x\in X\setminus \bigcup_{y\in X} \E_y$, then $[x]=\{x\}$. If there is some $y\in X$ such that $x\in \E_y$, then, from the definition of $\ssim$, it follows that $[x]\subseteq \E_y$. If $(X,\tau)$ is $T_1$, then, by Lemma~\ref{lemma_2}, the set $\E_y$ is finite, so $[x]$ is finite.  
\end{proof}

\begin{lemma}\label{properties_cofnite_subset_1.5}
If $(X,\tau)$ is regular and Hausdorff, then there exists a cofinite set $\Xgood\subseteq X$ with the following properties.
\begin{enumerate}[(a)]
\item \label{itm:properties_cofnite_subset_1.5_a} For any $x\in X$, either $[x] \subseteq \Xgood$ or $[x]\cap \Xgood= \emptyset$ (i.e. $X'$ is compatible with $\ssim$). 
\item  \label{itm:properties_cofnite_subset_1.5_b}For any $x\in \Xgood$ and $y\in X$, if $x\in \E_y$, then $y\in [x]$. 
\item \label{itm:properties_cofnite_subset_1.5_c} For any $x\in \Xgood$, $\E_x \subseteq [x]$. In particular, if $\E_x$ is non-empty, then $E_x=[x]$.  
\end{enumerate} 
\end{lemma}
\begin{proof} 
Let $H=\{(x,y)\in X^2 : y\in \E_x \text{ and } x\not\sim_\tau y\}$. Let $H_1$ and $H_2$ be the projections of $H$ onto the first and second coordinate respectively. We start by showing that these sets are finite. If $H_2$ is finite then, by Hausdorffness (see Proposition~\ref{lemma:basic_facts_P_x_2}(\ref{itm2:basic_facts_P_x_2})), there exist at most two $x \in X$ such that $y \in E_x$, hence $H_1$ is finite. Towards a contradiction, we suppose that $H_2$ is infinite. Let $g:H_2 \rightarrow H_1$ be the function given by $y \mapsto \min\{x : y\in \E_x \text{ and } y\not\sim_\tau x\}$. By Hausdorffness, this function is well defined and, by Lemma~\ref{lemma_2}, it cannot be constant on an interval, so, by o-minimality, there exists an interval $I\subseteq H_2$ such that $g|_I$ is strictly monotonic and $e$-continuous. But then, since $(X,\tau)$ is regular, by Lemma~\ref{lemma_f_Rx_Lx}, for any $y\in I$, it holds that $y \ssim g(y)$, a contradiction. 

Set $\Xgood:=X\setminus (\bigcup_{x\in H_1 \cup H_2} [x])$. By Lemma~\ref{lemma_finite_classes} and the finiteness of $H_1\cup H_2$, this set is cofinite in $X$. By definition of $H$, it follows that $X'$ satisfies \eqref{itm:properties_cofnite_subset_1.5_a}-\eqref{itm:properties_cofnite_subset_1.5_c}.
\end{proof}
The next lemma strengthens Lemma~\ref{properties_cofnite_subset_1.5} and is the core construction in the proofs of Theorems~\ref{them_ADC_or_DOTS} and~\ref{them_compactification}.

\begin{lemma}\label{remark_for_two_theorems}
Suppose that $(X,\tau)$ is regular and Hausdorff. Let $X'\subseteq X$ be as given by Lemma~\ref{properties_cofnite_subset_1.5}. There exists a finite partition $\XX$ of $X$ into singletons $\Xsing$ and infinite definable $\tau$-open sets $\Xinf$, the latter being subsets of $X'$, with the following properties.

For every $A\in\Xinf$, there exists $n>0$, an interval $I\subseteq A$, and definable \econtinuous strictly monotonic functions $f_0, f_1,\ldots, f_{n-1}:I\rightarrow A$ such that, for every $x\in I$, $[x]=\{ f_i(x) : 0\leq i < n\}$. In particular, $f_0$ is the identity map. 
Further, the intervals $f_i(I)$, for $0\leq i < n$, are pairwise disjoint and $A=\bigcup_{0\leq i <n} f_i(I)=\bigcup_{x\in I} [x]$. Additionally, for every $x\in I$ and $0<i<n-1$, the point $f_i(x)$ is $\tau$-isolated. 

Moreover, if we set $[x]^E=\{ y\in [x] : \E_y\neq \emptyset\}$ for every $x\in I$, then exactly one of the following conditions holds:
\[ 
\begin{split}
\forall x\in I, \quad & [x]=\{x\} \text{ and } \E_x=\emptyset \text{, so }[x]^E=\emptyset \text{ (and $A=I$ contains only} \\ 
& \text{$\tau$-isolated points);} \\
\forall x\in I, \quad & [x]^E=\{x\}; \\
\forall x\in I, \quad & [x]^E=\{x, f_{n-1}(x)\} \text{ (this case only applies if $n>1$)}. 
\end{split}
\] 

Finally, in each of the latter two cases, exactly one of the following conditions is also satisfied:
\[ 
\begin{split}
\forall x\in I,\quad &x\in L_x \setminus R_x; \\
\forall x\in I,\quad &x\in R_x \setminus L_x;\\
\forall x\in I,\quad &x\in R_x \cap L_x.
\end{split}
\] 
\end{lemma}
\begin{proof}
We construct $\XX$ by describing the family $\Xinf$ of $\tau$-open subsets of $X'$, while making sure that $\bigcup\Xinf$ is cofinite in $X$. Since the set $X'$ is cofinite in $X$ (see Lemma~\ref{properties_cofnite_subset_1.5}) it suffices to ensure that $\bigcup\Xinf$ is cofinite in $X'$. In particular, we consider a finite number of definable sets that partition $X'$ and, for each such set $S$, describe a partition of a cofinite subset of $S$, which becomes the collection of subsets of $S$ in $\Xinf$. 

Let $\Xisol=\Xgood \setminus \bigcup_{y\in X} \E_y$. Note that $[x]=\{x\}$ for every $x\in \Xisol$. By Lemma~\ref{properties_cofnite_subset_1.5}(c) and definition of $\Xisol$, for all $x\in \Xisol$, $\E_x =\emptyset$, and so, by Lemma~\ref{lemma_explaining_neighbourhoods}, these points are $\tau$-isolated. If $\Xisol$ is infinite, let $\XXisol$ be a finite family of disjoint intervals whose union is cofinite in $\Xisol$. Otherwise let $\XXisol=\emptyset$. 

We now consider $X' \setminus \Xisol$. By Lemma~\ref{lemma_finite_classes} and uniform finiteness, there exists $n'\geq 1$ such that, for every $x\in X$, $|[x]|\leq n'$. For every $1\leq n \leq n'$, set $X_n:=\{ x\in \Xgood\setminus \Xisol: |[x]|=n \}$. These sets are definable and partition $\Xgood\setminus \Xisol$. We fix $1\leq n\leq n'$. If $X_n$ is finite, let $\XX_n=\emptyset$. If $X_n$ is infinite, then the following describes a finite partition $\XX_n$ of a cofinite subset of $X_n$ into definable $\tau$-open sets as desired.

Recall the notation $[x]^E=\{y\in [x]: \E_y \neq \emptyset\}$ for $x\in X$. Since $\Xgood\setminus \Xisol\subseteq \bigcup_{y\in X} \E_y$ then, by Lemma~\ref{properties_cofnite_subset_1.5}(\ref{itm:properties_cofnite_subset_1.5_b}), for every $x\in X' \setminus \Xisol$, and hence for every $x\in X_n$, it holds that $|[x]^E|\geq 1$. Moreover, by Lemma~\ref{properties_cofnite_subset_1.5}(\ref{itm:properties_cofnite_subset_1.5_c}) and Hausdorffness (see Lemma~\ref{lemma:basic_facts_P_x_2}(\ref{itmb:lemma_basic_facts_P_x_2})), for every $x\in X' \setminus \Xisol$, and hence for every $x\in X_n$, we have that $|[x]^E|\leq 2$. Let $\Xone=\{ x\in X_n : |[x]^E|=1 \}$ and $\Xtwo=\{ x\in X_n : |[x]^E|=2 \}$. These sets partition $X_n$. 

Set $Dom(\Xone):=\bigcup\{[x]^E : x\in \Xone\}$, $Dom(\Xtwo):=\{\min[x]^E : x \in \Xtwo\}$ and $Dom(X_n):=Dom(\Xone)\cup Dom(\Xtwo)$. Note that, for every $x\in X_n$, it holds that $|Dom(X_n)\cap [x]|=1$.

First, let $f_0$ denote the identity map on $Dom(X_n)$. Then, for every $1\leq i < n$, let $f_i:Dom(X_n)\rightarrow X$ be the function defined as follows. For every $x\in Dom(\Xone)$, $f_i(x)$ is the $i$-th smallest element in $[x]\setminus\{x\}$. For every $x\in Dom(\Xtwo)$, if $1\leq i < n-1$, then $f_i(x)$ is the $i$-th smallest element in $[x]\setminus [x]^E$, while $f_{n-1}(x)=\max [x]^E$. By construction, for every $y\in X_n$ there exists a unique $x\in Dom(X_n)$ and unique $0\leq i < n$ such that $f_i(x)=y$. In particular, for every $x \in Dom(X_n)$, $[x]=\{f_i(x) : 0\leq i < n\}$, all functions $f_i$ are injective and the family of images $\{f_i(Dom(X_n)) : 0\leq i < n\}$ is pairwise disjoint and covers $X_n$. Moreover, by construction, for every $x\in Dom(X_n)$ and $0<i<n-1$, $\E_{f_i(x)}=\emptyset$, so, by Lemma~\ref{lemma_explaining_neighbourhoods}, $f_i(x)$ is $\tau$-isolated.  

Since, by definition of $Dom(X_n)$, $E_x \neq \emptyset$ for all $x\in Dom(X_n)$, by Lemma~\ref{properties_cofnite_subset_1.5}(\ref{itm:properties_cofnite_subset_1.5_c}) and definition of $\sim_\tau$ it follows that $x \in E_x$ for all $x\in Dom(X_n)$. So, by Lemma~\ref{prop_basic_facts_P_x_2}(\ref{itm2:basic_facts_P_x_2}), we may further partition $Dom(X_n)$ into three definable sets as follows.
\[
\begin{split}
Dom(X_n)'=\{ x\in Dom(X_n) : x\in L_x \setminus R_x\}, \\
Dom(X_n)''=\{ x\in Dom(X_n) : x\in R_x \setminus L_x\},\\
Dom(X_n)'''=\{ x\in Dom(X_n) : x\in R_x \cap L_x\}.
\end{split}
\]

By o-minimality, there exists a finite partition $\mathcal{D}_n$ of $Dom(X_n)$, compatible with $\{Dom(\Xone),Dom(\Xtwo), Dom(X_n)', Dom(X_n)'', Dom(X_n)'''\}$, which contains only singletons and intervals and is such that, for every interval $I\in\mathcal{D}_n$, $f_i|_I$ is \econtinuous and strictly monotonic, for every $0\leq i < n$. The family of those sets in $\Xinf$ which are subsets of $X_n$ is then defined by
\[
 \XX_n=\left\{\bigcup_{0\leq i < n} f_i(I) : I\in \mathcal{D}_n, I \text{ is an interval} \right\}. 
\]
Note that, by construction, $\bigcup_{0\leq i < n} f_i(I)= \bigcup_{x\in I} [x]$, for any interval $I\in \mathcal{D}_n$. Moreover, since the functions $f_i$ are $e$-continuous, the sets $f_i(I)$ are all intervals. It easily follows that any $A \in \XX_n$  is $\tau$-open, by observing that, for every $x \in A$, we have $E_x \subseteq [x]$ by Lemma~\ref{properties_cofnite_subset_1.5}(\ref{itm:properties_cofnite_subset_1.5_c}), and using the form of a basis of $\tau$-neighbourhoods for $x$ given by Lemma~\ref{lemma_explaining_neighbourhoods}.

Finally, set $\Xinf:=\XXisol \cup \XX_1 \cup \cdots \cup \XX_{n'}$, and let $\Xsing$ denote the collection of singletons given by the points in $X\setminus \bigcup\Xinf$. By construction, this partition satisfies the properties stated in the lemma. 
To check this, for any $A\in \Xinf$, if $A\subseteq X_n$ for some $n$, let $I$ and $f_i$, for $0\leq i <n$, be as described above. If $A\subseteq \Xisol$, then simply consider $I=A$ with $f_0$ denoting the identity map on $I$. 
\end{proof}

Continuing with the construction in Lemma~\ref{remark_for_two_theorems}, the next lemma describes how each of the sets $A\in \Xinf$ definably embeds into a space with either the lexicographic or Alexandrov $n$-line topology. 

We first require a definition extending the notion of $e$-convergence from the right or from the left (Remark~\ref{remark_side_convergence}), and a remark on how this convergence relates to convergence with respect to the topologies $\taulex$ and $\taualex$.
\begin{definition}\label{remark_side_convergence_1}
Given a definable set $\tilde{X} \subseteq R \times \{0,1,\ldots\}$ we say that a definable curve in $\tilde{X}$ \emph{$e$-converges to $\al x,i \ar\in \tilde{X}$ from the right} (respectively \emph{left}) if it $e$-converges to $\al x,i \ar$ and its projection to the first coordinate, namely $\pi \circ \gamma$, $e$-converges to $x$ from the right (respectively left). 
\end{definition}

\begin{remark} \label{remark:lex_side_convergence}
Consider a definable set $\tilde{X}=I\times\{0,\ldots, n\}$, with $I\subseteq R$ an interval, and an injective definable curve $\gamma$ in $\tilde{X}$. For $x\in I$, note that, by o-minimality, if $n>0$, then $\gamma$ converges in $(\tilde{X},\taulex)$ to $\al x,0 \ar$ if and only if it $e$-converges to $\al x,i \ar$ from the left, for some $0\leq i \leq n$ (recall the basis of open neighbourhoods for $\al x, 0 \ar$ described in Definition~\ref{dfn:lex}). Similarly, maintaining the assumption that $n>0$, we have that $\gamma$ converges in $(\tilde{X},\taulex)$ to $\al x,n \ar$ if and only if it $e$-converges to $\al x,i \ar$ from the right, for some $0\leq i \leq n$ (likewise consider the basis of open neighbourhoods for $\al x, n \ar$ given earlier). Moreover (recalling Definition~\ref{example_n_line_0}), by o-minimality, $\gamma$ converges to $\al x,0 \ar$ in $(\tilde{X},\taualex)$ if and only if it $e$-converges to $\al x,i \ar$, for some $0\leq i \leq n$ (from the right or from the left).  
\end{remark}

\begin{lemma}\label{lemma_proof_two_thems}
Suppose that $(X,\tau)$ is regular and Hausdorff. Let $\mathcal{X}$ be a finite partition of $X$ as given by Lemma~\ref{remark_for_two_theorems}. For each $A=\bigcup_{0\leq i <n} f_i(I)\in \Xinf$, there exists a definable set $A^*\subseteq R^2$ and a definable injection \mbox{$h_A:A\rightarrow A^*$} such that the following hold. 
\begin{enumerate}[(1)]
\item  \label{itm1:lemma_proof_two_theorems} $A^*=I\times \{0,\ldots, m\}$, for some $m\in\{n-1, n, 2\}$. In particular, for every distinct pair $A_0, A_1$ in $\Xinf$ we have that $A_0^* \cap A_1^*=\emptyset$. 

\item \label{itm1.5:lemma_proof_two_theorems} For every $0\leq i <n$ and $x\in f_i(I)$, it holds that $(\pi\circ h_A)(x)= f_i^{-1}(x)$.

\item \label{itm2:lemma_proof_two_theorems} The map $h_A:(A,\tau)\hookrightarrow (A^*,\tau_A)$ is an embedding, where either 
\begin{enumerate}
\item $(A^*,\tau_A)=(A^*,\taulex)$ or
\item $(A^*,\tau_A)=(A^*,\taualex)$.
\end{enumerate}
\end{enumerate}
\end{lemma}
\begin{proof}

Following Lemma~\ref{remark_for_two_theorems}, we distinguish different cases of possible sets $A=\bigcup_{0\leq i <n} f_i(I)$ in $\Xinf$, based on properties of $I$. In each case, we define $A^*$ and $h_A$ (which for simplicity we denote by $h$) so that~(\ref{itm1:lemma_proof_two_theorems}), (\ref{itm1.5:lemma_proof_two_theorems}) and~(\ref{itm2:lemma_proof_two_theorems}) hold.

\textbf{Case 0:} $I=A$ is a set of isolated points in $(X,\tau)$.  
Specifically, $[x]=\{x\}$ and $\E_x=\emptyset$, for every $x\in I$. 

Let $A^*=A\times \{0,1,2\}$ and let $h:A\rightarrow A^*$ be given by $x\mapsto \al x,1 \ar$. Let $\tau_A$ be the topology induced by the lexicographic order on $A^*$. With this topology all the points in $A\times \{1\}$ are isolated, so $h$ is an embedding.

For the remaining cases we will make use of Proposition~\ref{prop_cont_lim} to prove~(\ref{itm2:lemma_proof_two_theorems}), that is, after specifying what $A^*$, $\tau_A$ and $h$ will be in each case in order to satisfy the other requirements, we will fix a point $x\in X$ and prove that an injective definable curve $\gamma$ converges in $(A,\tau)$ to $x$ if and only if $h\circ \gamma$ converges in $(A^*,\tau_A)$ to $h(x)$. To do this we will use Remarks~\ref{remark_side_convergence} and~\ref{remark:lex_side_convergence} extensively.

\textbf{Case 1:} $[x]^E=\{x\}$ and $x\in L_x \setminus R_x$, for every $x\in I$.

In this case, for every $x \in I$, we have the following. 
Since $[x]^E=\{x\}$, by Lemma~\ref{properties_cofnite_subset_1.5}\eqref{itm:properties_cofnite_subset_1.5_c} it holds that $\E_x = [x]$, hence $E_x = \{f_i(x) : 0\leq i <n\}$. Furthermore, any point $y\in [x]\setminus \{x\}$ satisfies that $\E_y=\emptyset$ and so, by Lemma~\ref{lemma_explaining_neighbourhoods}, it is $\tau$-isolated. 
Moreover, we will make use of the following observation, which follows from applying Lemma~\ref{lemma_f_Rx_Lx} and the fact that $x\in L_x \setminus R_x$, for all $x \in I$.
For every $0\leq i <n$ and $x \in I$,
\begin{equation}\label{eqn_cases_1} 
\begin{split}
f_i(x)\in L_x \Leftrightarrow f_i \text{ is increasing}; \\
f_i(x) \in R_x \Leftrightarrow f_i \text{ is decreasing}.
\end{split}  
\end{equation}

Let $A^*= I \times \{0,\ldots, n\}$ 
and let $h:A\rightarrow A^*$ be the definable injection given by $h(f_i(x))=\al x,i \ar$,
for every $x\in I$ and $0\leq i <n$.  

Note that, for each $0\leq i <n$, $h|_{f_i(I)}$ is given by $x\mapsto \al f_i^{-1}(x),i \ar$, and so $h$ is an $e$-embedding. 
We show (using Proposition~\ref{prop_cont_lim}) that $h$ is an embedding $(A,\tau)\hookrightarrow (A^*,\taulex)$. 

Fix $x\in A$. If $x\in f_i(I)$ with $i>0$, then both $x$ and $h(x)$ are isolated, and hence not the limit of any injective definable curve, so we may assume that $x\in I$. 
Let $\gamma$ be an injective definable curve in $A$ that $\tau$-converges to $x$. 
By o-minimality, $\gamma$ is $e$-convergent. By the fact that $\E_x=\{f_i(x) : 0\leq i <n\}$ and Proposition~\ref{basic_facts_P_x_1}\eqref{itm: basic_facts_3}, we have that $\gamma$ necessarily $e$-converges to some $f_i(x)$, from either the right or the left (see Remark~\ref{remark_side_convergence}). Since $h$ is an $e$-homeomorphism, 
$h\circ\gamma$ $e$-converges to $\al x,i \ar$. If $\gamma$ $e$-converges from the left, then, by Remark~\ref{remark_side_convergence}, it must be that $f_i(x) \in L_x$.  
In this case it follows from~\eqref{eqn_cases_1} that $f_i$ is increasing, and consequently $h \circ \gamma$ also $e$-converges to $\al x,i \ar$ from the left, and therefore (by Remark~\ref{remark:lex_side_convergence}), $h \circ \gamma$ converges in $(A^*,\taulex)$ to $\al x,0 \ar=h(x)$. Similarly, if $\gamma$ $e$-converges from the right, then it must be that $f_i(x)\in R_x$, and so, by~\eqref{eqn_cases_1}, $f_i$ is decreasing, and thus $h \circ \gamma$ again $e$-converges to $\al x,i \ar$ from the left, so again it converges in $(A^*,\taulex)$ to $\al x,0 \ar=h(x)$.

Conversely, let $\gamma'\subseteq h(A)$ be an injective definable curve converging in $(A^*,\taulex)$ to $\al x,0 \ar$, in which case it must $e$-converge from the left to some $\al x,i \ar$ (see Remark~\ref{remark:lex_side_convergence}). We may assume that $\gamma'\subseteq I\times\{i\}$, for some $0 \leq i < n$ (see Remark~\ref{remark_assumptions_curves}). Note that $h^{-1}\circ \gamma' = f_i \circ \pi \circ \gamma'$. If $f_i$ is increasing, then, by~\eqref{eqn_cases_1}, $f_i(x)\in L_x$, and moreover $h^{-1}\circ \gamma'$ $e$-converges to $f_i(x)$ from the left, so it $\tau$-converges to $x$ (see Remark \ref{remark_side_convergence}). Similarly, if $f_i$ is decreasing, then, by~\eqref{eqn_cases_1}, $f_i(x)\in R_x$, and moreover $h^{-1}\circ \gamma$ $e$-converges to $f_i(x)$ from the right, so again it $\tau$-converges to $x$. 

\textbf{Case 2:} $[x]^E=\{x\}$ and $x\in R_x \setminus L_x$ for every $x\in I$. 

This case is very similar to Case 1, so we only indicate here the key details and requisite changes. 
As in Case 1, it holds that $\E_x=[x]=\{f_i(x) : 0\leq i <n\}$ for every $x\in I$, and every point $y\in A\setminus I$ is $\tau$-isolated. On the other hand, in this case we have the following two equivalences arising from Lemma~\ref{lemma_f_Rx_Lx}. For every $0\leq i <n$ and $x\in I$,
\begin{equation}\label{eqn_cases_2}
\begin{split}
f_i(x)\in R_{x} \Leftrightarrow f_i \text{ is increasing}, \\
f_i(x) \in  L_{x} \Leftrightarrow f_i \text{ is decreasing}.
\end{split}  
\end{equation}

We let $A^*=I\times \{0,\ldots, n\}$ and let $h:A\rightarrow A^*$ be the definable injection given by 
$
h(f_i(x))=\al x,n-i \ar,
$
for $x\in I$ and $0\leq i <n$. Again, this is clearly an $e$-embedding. It can be shown that $h:(A,\tau)\hookrightarrow (A^*,\taulex)$ is an embedding by analogy to Case 1.
 That is, we fix $x\in I$ and a definable curve $\gamma$ in $A$ that $\tau$-converges to $x$. Then we observe that $\gamma$ $e$-converges to some $f_i(x)$ and $h\circ \gamma$ $e$-converges to $h(f_i(x))=\al x, n-i \ar$. Finally, we consider separately the cases where $\gamma$ $e$-converges to $f_i(x)$ from the left and from the right and show, using~\eqref{eqn_cases_2} together with Remarks~\ref{remark_side_convergence} and~\ref{remark:lex_side_convergence}, that in both cases $h\circ \gamma$ converges in $(A^*,\taulex)$ to $h(x)$. The converse, i.e. the case of a definable curve $\gamma'$ in $h(A)$, may be similarly argued by analogy to Case 1.

\textbf{Case 3:} $n>1$, $[x]^E=\{x, f_{n-1}(x)\}$ and $x\in L_x \setminus R_x$ for every $x\in I$. 

In this case we have that, by Lemma~\ref{properties_cofnite_subset_1.5}\eqref{itm:properties_cofnite_subset_1.5_c}, for every $x \in I$ it holds that $\E_x = \E_{f_{n-1}(x)}=[x]=\{f_i(x) : 0\leq i <n\}$. Moreover, every point $y\in A\setminus (I \cup f_{n-1}(I))$ satisfies that $\E_y=\emptyset$ and so (by Lemma~\ref{lemma_explaining_neighbourhoods}) is $\tau$-isolated. Furthermore, by Hausdorffness of $\tau$ (see Proposition~\ref{prop_basic_facts_P_x_2}\eqref{itm3:basic_facts_P_x_2}), it holds that, for every $0\leq i <n$ and $x\in I$, exactly one of the following two possibilities is satisfied: 
\begin{equation}\label{eqn_cases_3}
\begin{split}
f_i(x)\in R_x \cap L_{f_{n-1}(x)};\\ 
f_i(x)\in R_{f_{n-1}(x)} \cap L_x. 
\end{split}
\end{equation}

Let $A^*=I\times\{0,\ldots, n-1\}$, and let $h:A\rightarrow A^*$ be defined in a similar manner to Case $1$, namely by $h(f_i(x))=\al x,i \ar$, for every $x\in I$ and $0\leq i <n$. In this case, $h:A\rightarrow A^*$ is a bijection. We show that $h$ is a homeomorphism $(A,\tau)\rightarrow (A^*,\taulex)$, by showing that an injective definable curve $\gamma$ in $A$ $\tau$-converges to a point $y\in A$ if and only if $h \circ \gamma$ $\taulex$-converges to $h(y)$. 

We fix $y\in A$. The case $y\in f_i(I)$, for $0<i<n-1$, is as usual trivial, since in this case both $y$ and $h(y)$ are isolated in their respective spaces. If $y\in I$, then the required statement in that case follows from the corresponding argument in Case 1. 

Therefore, suppose that $y\in f_{n-1}(I)$.  
Let $x \in I$ be such that $y=f_{n-1}(x)$. 
Let $\gamma$ be an injective definable curve in $A$ $\tau$-converging to $y$. Using that $\E_y=\{f_i(x) : 0\leq i <n\}$ and following the arguments in Case 1, we observe that $\gamma$ $e$-converges to some $f_i(x)$. Since $h$ is an \ehomeomorphismp, $h\circ\gamma$ $e$-converges to $\al x,i \ar$, a priori from either the right or from the left (see Remark~\ref{remark_side_convergence}). However, if $h \circ \gamma$ $e$-converges to $\al x, i \ar$ from the left -- converging thus in $(A^*,\taulex)$ to $\al x,0 \ar$ (by Remark~\ref{remark:lex_side_convergence}) -- then, by continuity of $h^{-1}$ at $h(x)=\al x,0 \ar$ (i.e. using that the statement is already established for all points in $I$), it must be that $h^{-1}\circ h\circ \gamma = \gamma$ $\tau$-converges to $x$, a contradiction (by Hausdorffness of $\tau$). So $h\circ \gamma$ must $e$-converge to $\al x, i \ar$ from the right, meaning that it converges in $(A^*,\taulex)$ to $\al x,n-1 \ar=h(y)$ (see Remark~\ref{remark:lex_side_convergence}). 

Conversely, let $\gamma'$ be an injective definable curve converging in $(A^*,\taulex)$ to $h(y)=\al x,n-1 \ar$. Then it $e$-converges to some $\al x,i \ar$, meaning that $h^{-1}\circ \gamma'$ $e$-converges to $f_i(x)$. If $h^{-1}\circ\gamma'$ does not $\tau$-converge to $y$ then, by~\eqref{eqn_cases_3} and Remark~\ref{remark_side_convergence}, it $\tau$-converges to $x$, but then, by continuity of $h$ at $x$, it follows that $\gamma'$ converges in $(A^*,\taulex)$ to $h(x)=\al x,0 \ar$, a contradiction. So $h^{-1}\circ \gamma'$ $\tau$-converges to $y$.

\textbf{Case 4:} $n>1$, $[x]^E=\{x, f_{n-1}(x)\}$ and $x\in R_x \setminus L_x$ for every $x\in I$. 

Again let $A^*=I\times\{0,\ldots, n-1\}$, and now let $h:A\rightarrow A^*$ be given in a similar manner to Case 2, namely by $h(f_i(x))=\al x,n-1-i \ar$. Note that $h$ is again a bijection. Moreover we may again in this case show that $h:(A,\tau)\rightarrow (A^*,\taulex)$ is a homeomorphism. 

The proof of (\ref{itm2:lemma_proof_two_theorems}) here follows from the proofs of the other cases. The argument in Case $2$ shows that, for any $x\in A\setminus f_{n-1}(I)$, both $h$ and $h^{-1}$ are continuous at $x$ and $h(x)$ respectively. Then, for the points in $f_{n-1}(I)$, one may use an argument analogous to the one in Case 3. 
 
\textbf{Case 5:}  $x\in R_x \cap L_x$ for every $x\in I$.

Set $A^*= I \times \{0,\ldots,n-1\}$ and let $h:A\rightarrow A^*$ be given by $h(f_i(x))=\al x,i \ar $. This map is clearly bijective. We show that it is a homeomorphism $(A,\tau)\rightarrow (A^*,\taualex)$.

If there exists $0 < i < n$ and $x \in I$ such that $E_{f_{i}(x)} \neq \emptyset$, we must have, by Lemma~\ref{properties_cofnite_subset_1.5}(\ref{itm:properties_cofnite_subset_1.5_c}), that $E_{f_{i}(x)} = [f_i(x)] = [x]$, and in particular $x \in R_{f_i(x)} \cup L_{f_i(x)}$, which contradicts Hausdorffness (Proposition~\ref{prop_basic_facts_P_x_2}\eqref{itm3:basic_facts_P_x_2}) and the fact that $x\in R_x \cap L_x$. Thus, for every $x\in I$ and $0<i<n$, we have $E_{f_{i}(x)} = \emptyset$, meaning that the point $f_i(x)$ is $\tau$-isolated. Moreover, by definition of the topology $\taualex$, the points $\al x,i \ar$ for $x\in I$ and $0<i<n$ are isolated in $(A^*,\taualex)$.

Applying Lemma~\ref{lemma_f_Rx_Lx}, we note that, for every $x\in I$, it follows from $x\in R_x \cap L_x$ that $f_i(x)\in R_x \cap L_x$, for every $0\leq i <n$. Furthermore, by Lemma~\ref{properties_cofnite_subset_1.5}(\ref{itm:properties_cofnite_subset_1.5_c}), we have that $R_x=L_x=\{ f_i(x) : 0\leq i < n\}$. Consequently, for any $x\in I$, any injective definable curve converges in $(A,\tau)$ to $x$ if and only if it $e$-converges to some $f_i(x)$ (see Remark~\ref{remark_side_convergence}). Similarly, any injective definable curve converges in $(A^*,\taualex)$ to $\al x,0 \ar$ if and only if it $e$-converges to some $\al x,i \ar$ (see Remark~\ref{remark:lex_side_convergence}). Thus the result follows from the fact that $h$ is a \ehomeomorphismp, which is clear from the definition.

This covers all possible cases for $A$, and thus completes the proof of the lemma. 
\end{proof}

We may now prove Theorem~\ref{them_ADC_or_DOTS}.

\begin{proof}[Proof of Theorem~\ref{them_ADC_or_DOTS}]
Let $\XX$ be a partition of $X$ as given by Lemma~\ref{remark_for_two_theorems} and, for each $A\in \Xinf$, let $A^*$, $h_A$ and $\tau_A$ be as given by Lemma~\ref{lemma_proof_two_thems}. 

Set $h:=\bigcup\{h_A : A\in\Xinf\}$. By Lemma~\ref{lemma_proof_two_thems}(\ref{itm1:lemma_proof_two_theorems}), $h$ is an injection $\bigcup\Xinf \rightarrow \bigcup_{A\in\Xinf} A^*$. Let 
\[Y=\bigcup\{A\in\Xinf : (A^*,\tau_A)=(A^*,\taulex)\neq (A^*,\taualex)\}\]
and 
\[Z=\bigcup\{ A\in\Xinf : (A^*,\tau_A)=(A^*,\taualex)\}.\]
By construction, these sets are disjoint, $\tau$-open and definable, and $X\setminus (Y \cup Z)$ is finite. Set $Y^*:=\bigcup\{ A^* : A\in \Xinf,\, A\subseteq Y\}$ and $Z^*:=\bigcup\{ A^* : A\in \Xinf,\, A\subseteq Z \}$. We claim that $h|_Y:(Y,\tau)\rightarrow (Y^*,\taulex)$ and $h|_Z:(Z,\tau)\rightarrow (Z^*,\taualex)$ are embeddings. We show that this claim holds by noting that we may decompose the maps into embeddings between open subspaces of their domains and codomains. 

Recall that, by Lemma~\ref{lemma_proof_two_thems}(\ref{itm1:lemma_proof_two_theorems}), for each $A\in \Xinf$, the set $A^*$ is of the form $I\times \{0,\ldots, m\}$, for some $m$ and interval $I\subseteq A$. It follows that, for each $A\in \Xinf$, if $A\subseteq Z$, then $A^*$ is open in $(Z^*,\taualex)$. Similarly, if $A \subseteq Y$, then $A^*$ is open in $(Y^*,\taulex)$, and the subspace topology on $A^*$ is precisely the topology induced by the lexicographic order on $A^*$. The claim then follows from Lemma~\ref{lemma_proof_two_thems}(\ref{itm2:lemma_proof_two_theorems}).   

Finally, if $Y\neq \emptyset$, then let $n_Y=\max\{n : (R\times\{n\})\cap Y^* \neq \emptyset\}$. It is easy to see that the map on $Y^*$ given by $\al x,m \ar\mapsto \al x,n_Y \ar$, if $m=\max\{m' : \al x, m'\ar \in Y^*\}$, and the identity otherwise is a definable embedding $(Y^*,\taulex)\hookrightarrow (R\times\{0,\dots,n_Y\},\taulex)$. Hence we conclude that $(Y,\tau)$ embeds definably into $(R\times\{0,\dots,n_Y\},\taulex)$. Furthermore, we may define $n_Z$ analagously to $n_Y$, and then $h|_Z$ is clearly a definable embedding $(Z,\tau) \hookrightarrow (R\times\{0,\dots,n_Z\},\taualex)$.
\end{proof}

\begin{remark}~\label{remark_general_them_ADC_or DOTS} 
If $\RR$ expands an ordered field, then, by Remark~\ref{remark_assumption_X}, Theorem~\ref{them_ADC_or_DOTS} applies to all $T_3$ one-dimensional spaces. Otherwise, the theorem and its proof may be rewritten to apply to one-dimensional spaces as follows. 

Let $(X,\tau)$ be a $T_3$ one-dimensional definable topological space. In the context of Remark~\ref{remark_generalizing_Rx_Lx}, one may prove analogues of Lemmas~\ref{lemma_finite_classes} and~\ref{properties_cofnite_subset_1.5}, and reach a partition $\XX=\Xinf \cup \Xsing$ of $X$ analogous to the one described in Lemma~\ref{remark_for_two_theorems}. Then, analogously to the proof of Lemma~\ref{lemma_proof_two_thems}, it is possible to show that, for any $A\in \Xinf$, there exists some $m \geq 0$, an interval $I$, and a definable embedding $h_A:(A,\tau)\hookrightarrow (A^*,\tau_A)$, where $A^*=I \times \{0, \ldots, m\}$, and $\tau_A\in\{\taulex, \taualex\}$. We may then derive that $X$ has a cofinite subset that is definably homeomorphic to the disjoint union of finitely many spaces of the form $(R \times \{0, \ldots, m\},\taulex)$ or $(R \times \{0, \ldots, m\},\taualex)$, for various $m$. 
\end{remark}

In~\cite{ram13}, Ramakrishnan shows that, if $\RR$ has definable choice and defines an order-reversing injection (e.g. if $\RR$ expands an ordered group), then every definable linear order definably embeds into $(R^n,\lex)$, for some $n$. In particular, under these assumptions, for any definable order topological space one may assume that, up to definable homeomorphism, the topology is induced by the lexicographic order. Theorem~\ref{them_ADC_or_DOTS} adds to the understanding of $T_3$ definable spaces in the line by describing how the $\taualex$ topology also plays a role describing them. We complete this picture with the next proposition.

\begin{proposition}\label{prop_n_line_is_not_order_space}
For any interval $I$ and any $n>0$, the space $(I\times\{0,\ldots, n\},\taualex)$ does not definably embed into a definable order topological space.
\end{proposition}
\begin{proof}
It suffices to prove the propostion for $n=1$. Towards a contradiction, assume that there exists an interval $I$ and one such embedding from $(I\times\{0,1\},\taualex)$ into a definable topological space $(X,\tau)$, where $\tau$ is given by a definable linear order $\precc$. Let $Y$ denote the image of the aforementioned embedding. 

We begin by observing that any clopen definable subset of $(I\times\{0,1\},\taualex)$ is either finite or cofinite. We then complete the proof by contradiction by showing that $(Y,\tau)$ contains an infinite coinfinite definable clopen subset.

Let $Z\subseteq I\times\{0,1\}$ denote a $\taualex$-clopen definable subset. Note that, by o-minimality, the subspace $(I\times\{0\}, \taualex)=(I\times\{0\},\tau_e)$ is definably connected, and so it must be that either $Z \cap (I\times \{0\}) = \emptyset$ or $Z \cap (I\times \{0\}) = I\times \{0\}$. We show that in the first case $Z$ is finite and in the second it is cofinite. By passing from $Z$ to $(I\times\{0,1\})\setminus Z$ in the first case it suffices to prove the second case. By o-minimality, in order to prove that $Z$ is cofinite it suffices to show that it contains a point in every set of the form $J\times \{i\}$ for any interval $J\subseteq I$ and $i\in \{0,1\}$. However, this follows directly from the condition $Z \cap (I\times \{0\}) = I\times \{0\}$, and the fact that $Z$ is open in $(I\times\{0,1\}, \taualex)$. 

We now show that $(Y,\tau)$ contains an infinite coinfinite definable clopen subset.
Let $Y_1=\{ x\in Y : (x,+\infty)_{\precc}\cap Y \text{ is finite}\}$. Since the intervals $(x,+\infty)_{\precc}$ are nested note that, by uniform finiteness, the set $Y_1$ is finite. Similarly the set $Y_2=\{x\in Y : (-\infty, x)_{\precc}\cap Y \text{ is finite}\}$ is also finite. 

Since $(Y,\tau)$ is homeomorphic to $(I\times\{0,1\},\taualex)$ and all the points in $I\times\{1\}$ are $\taualex$-isolated, it follows that $(Y,\tau)$ has infinitely many isolated points. Let us fix $x\in Y$, an isolated point in $(Y,\tau)$ that does not belong in $Y_1\cup Y_2$. There must exist $y \prec z$ in $X$ such that 
\begin{equation*}
(y,z)_{\precc} \cap Y=\{x\}.
\end{equation*} 
It follows that the set $(x,+\infty)_{\precc}\cap Y= [z,+\infty)_{\precc} \cap Y$ is clopen in $(Y,\tau)$. 
Since $x\notin Y_1\cup Y_2$, then this set is also infinite and coinfinite in $Y$, the desired contradiction.  
\end{proof}

Finally, the next corollary of Theorem~\ref{them_ADC_or_DOTS} implies that, for any $T_3$ definably separable definable topological space $(X,\tau)$, where $X\subseteq R$, there exists a cofinite subset $Y\subseteq X$ such that $\tau|_Y$ is induced by a definable linear order. 

\begin{corollary}\label{cor_T3_separable_embedding}
Let $(X,\tau)$, $X\subseteq R$, be a regular Hausdorff definable topological space. If $(X,\tau)$ is definably separable, then there exists a cofinite subset $Y\subseteq X$, a definable set $Y^*\subseteq R\times\{0,1\}$ and a definable embedding $(Y,\tau)\rightarrow (Y^*,\taulex)$. 
\end{corollary}
\begin{proof}
Let $\XX$ be a finite partition of $X$ as given by Lemma~\ref{remark_for_two_theorems} and, for each $A \in \Xinf$, let $A^*$ and $h_A$ be as given by Lemma~\ref{lemma_proof_two_thems}. Our aim is to show that, under the additional assumption that $(X,\tau)$ is definably separable, the construction in these lemmas yields that, for every $A\in\Xinf$, we have that $A^* \subseteq R\times \{0,1\}$ and $(A^*,\tau_A)=(A^*,\taulex)$. We will then let $Y=\bigcup \Xinf$ and $Y^*=\bigcup\{A^* : A\in\Xinf\}$. Following the proof of Theorem~\ref{them_ADC_or_DOTS}, $Y$ is cofinite and $h=\bigcup \{ h_A : A\in\Xinf\}$ is an embedding $(Y,\tau)\hookrightarrow (Y^*,\taulex)$.

If $(X,\tau)$ is definably separable, then it can have only finitely many isolated points. Following Lemma~\ref{remark_for_two_theorems}, fix $A=\bigcup_{0\leq i <n} f_i(I) \in \Xinf$. For every $x\in I$ and $0<i<n-1$, the point $f_i(x)$ is $\tau$-isolated. It follows that we must have $n\leq 2$. Similarly, for every $x\in I$, it must be that $[x]^E=\{x, f_{n-1}(x)\}$, since otherwise $\E_{y}=\emptyset$ (i.e. $y$ is $\tau$-isolated), for every $y\in f_{n-1}(I)$.

If $n=1$, this is considered in Cases 1, 2 and 5 in the proof of Lemma~\ref{lemma_proof_two_thems}. In Cases 1 and 2, $A^*=I\times \{0,1\}$ and $(A^*,\tau_A)=(A^*,\taulex)$. In Case 5, we have that $A^*=I\times \{0\}$ and $(A^*,\tau_A)=(A^*,\taualex)$, and in this case $(A^*,\taualex)=(A^*,\tau_e)=(A^*,\taulex)$.  
If $n=2$, this is considered in Cases 3 and 4 in the proof of Lemma~\ref{lemma_proof_two_thems}. In both of these cases, we have that $A^*=I\times \{0,1\}$ and $(A^*,\tau_A)=(A^*,\taulex)$. 
\end{proof}

\section{Some universality results} \label{section:universal_spaces_2}

In this section, we consider universality questions in the o-minimal definable setting, building on results from the previous two sections. This will also lead us to introduce several key concepts that will be important for the results in subsequent sections.

The main type of question that we examine here is the following. Given a certain class of definable topological spaces $\mathcal{C}$, is there a topological space $(X, \tau)$ that is universal for $\mathcal{C}$ in a definable sense, i.e. such that every space in $\mathcal{C}$ embeds definably into $(X,\tau)$? Moreover, is there such a universal space that lies within the given class $\mathcal{C}$ itself? We also consider a closely related property, which is natural to consider in the context of o-minimality, that we call `almost definable universality'. Our analysis is very much in the spirit of universality questions considered in the classical study of Banach spaces, where there is an extensive literature on universal spaces for classes of separable Banach spaces, going back to the classical Banach--Mazur Theorem \cite{banach} (see for example \cite{szlenk_68}, \cite{bourgain_80}, \cite{bossard_02},  \cite{odell_schlumprecht_06}, \cite{brech_kosz_12}). This question has also classically been studied in the context of topological spaces, where results of this kind include the Menger--N\"obeling Theorem \cite[Theorem V.2]{hur_wall41}, which states that $(\mathbb{R}^{2n+1},\tau_e)$ is universal for the class of all $n$-dimensional compact metric spaces.

We begin with the main definitions.

\begin{definition} \label{dfn:(almost)universal}
Let $\CC$ be a class of definable topological spaces and let $(X,\tau)$ be a definable topological space. 
We say that $(X,\tau)$ is \emph{definably universal for $\CC$} if every space $(Y,\mu)\in\CC$ embeds definably into $(X,\tau)$. 

We say that $(X,\tau)$ is \emph{almost definably universal for $\CC$} if, for every $(Y,\mu)\in\CC$, there exists a definable subset $Z\subseteq Y$ with $\dim (Y\setminus Z) < \dim Y$ such that $(Z,\mu)$ embeds definably into $(X,\tau)$. 
\end{definition}

Note that, if a space is definably universal for a class, then in particular it is almost definably universal. 

We begin by observing how o-minimality implies that, if $\RR$ expands an ordered field, then $R^n$ is almost definably universal for the class of euclidean spaces of dimension at most $n$. An analogous result can be proved for the class of bounded euclidean spaces of dimension at most $n$ when $\RR$ expands an ordered group, using the observations in Remark~\ref{remark_assumption_X}. Moreover, a slight weakening of the result can also be obtained in the most general case, which we discuss below in Remark~\ref{remark:universality_euclidean_spaces}.

\begin{proposition}\label{prop:universal_euclidean_space}
Suppose that $\RR$ expands an ordered field and $n>0$. Then $(R^n,\tau_e)$ is almost definably universal for the class of euclidean spaces of dimension less than or equal to $n$.
\end{proposition}
\begin{proof}

Let $(X,\tau_e)$ be a euclidean space. We prove the case $\dim X = n$. The case $\dim X = m$, with $1\leq m <n$, then follows by the easy fact that, for any such $m$, the space $(R^m,\tau_e)$ embeds definably into $(R^n, \tau_e)$, while in the case $\dim X = 0$, i.e. where $X$ is finite, $(X,\tau_e)$ clearly embeds definably into $(R^n, \tau_e)$.

Applying Remark~\ref{remark_assumption_X}, let $(Y,\mu)$ denote the push-forward of $(X,\tau_e)$ into $R^n$ by some definable injection $f:X\rightarrow R^n$. We prove that $(Y,\mu)$ contains a definable subspace $Z\subseteq Y$, with $\dim (Y\setminus Z)<\dim Y$, where the subspace topology is euclidean.

Observe that, by o-minimal cell decomposition applied to $f$ and $f^{-1}$, the function $f$ is a finite union of definable $e$-homeomorphisms. 
It follows that $Y$ can be partitioned into finitely many cells $\mathcal{D}$ where the subspace topology is euclidean (this is the property discussed in Section~\ref{section: affine} of being `cell-wise euclidean'; see Definition~\ref{dfn:cell-wise_euclidean}). Let \mbox{$Z=\bigcup\{int_\mu D : D\in\mathcal{D},\, \dim D=n\}$}. Note that, since $(Y,\mu)$ is the push-forward of a euclidean space, it has the frontier dimension inequality. It follows that $\dim (Y\setminus Z) < \dim Y$. Note that, since any cell of dimension $n$ is $e$-open, we have that, for any $D\in\mathcal{D}$ with $\dim D=n$, the set $int_\mu D$ is $e$-open, as well as $\mu$-open, in $Z$. Moreover the subspace topology on $int_\mu D$ is euclidean. We conclude that the subspace topology on $Z$ is euclidean. 
\end{proof}

\begin{remark}\label{remark:universality_euclidean_spaces}
In the general case where $\RR$ does not necessarily expand an ordered field one may still adapt the proof of Proposition~\ref{prop:universal_euclidean_space} to show that, if $(X,\tau)$ is a euclidean space of dimension $n$ with cell decomposition $\mathcal{D}$, then the union $Z$ of the interiors in $X$ of cells in $\mathcal{D}$ of dimension $n$ embeds definably into finitely many disjoint copies of $R^n$ (one for each cell). It follows that the space $(R^{n+1}, \tau_e)$ is almost definably universal for the class of euclidean spaces of dimension at most $n$.  
\end{remark}

The question of definable universality (as opposed to almost definable universality) for euclidean spaces is less straightforward. Walsberg has announced to the authors (through private correspondence) that he and C. Miller have obtained a definable version of the classical Menger--N\"obeling Theorem which implies that, whenever $\RR$ expands an ordered field, any euclidean space of dimension $n$ embeds definably into ($R^{2n+1}, \tau_e)$. 

We now illustrate how results from previous sections can be used to derive o-minimal definable universality results.
To begin, we show that Theorem~\ref{them_ADC_or_DOTS} may be framed in terms of existence of an almost definably universal space as follows.

\begin{corollary}\label{cor_universal_space}
The disjoint union of $(R\times [0,1],\taulex)$ and $(R\times[0,\infty), \taualex)$ is Hausdorff, regular and almost definably universal for the class of Hausdorff regular definable topological spaces $(X,\tau)$, where $X\subseteq R$.  
\end{corollary}
\begin{proof}
It is easy to observe that the spaces $(R\times[0,1], \taulex)$ and $(R\times[0,\infty),\taualex)$ are Hausdorff and regular, from where it follows that their disjoint union is too.  

Let $(X,\tau)$, where $X\subseteq R$, be a regular Hausdorff definable topological space. By Theorem~\ref{them_ADC_or_DOTS}, there exist definable disjoint open sets $Y,Z\subseteq X$ and $n_Y>0$ such that $X\setminus (Y\cup Z)$ is finite and there are definable embeddings $(Y,\tau)\hookrightarrow (R\times\{0,\ldots, n_Y\},\taulex)$ and $(Z,\tau)\hookrightarrow (R\times [0,\infty),\taualex)$. Hence it suffices to show that, for any $n>0$, there exists a definable embedding $(R\times\{0,\dots,n\},\taulex) \hookrightarrow (R\times [0,1],\taulex)$. Fix parameters $0=a_0<a_1<\cdots<a_{n}=1$. Then the map given by $\al x,i \ar \mapsto \al x, a_i \ar$ does the job.  
\end{proof}

In the specific case of spaces in the line that embed definably into definable order topological spaces we may refine the above corollary as follows.

\begin{corollary} \label{cor_universal_space_embed_ordertopspace}
    The disjoint union of $(R\times [0,1], \taulex)$ and $(R,\tau_e)$, with topology given by the lexicographic order, is almost definably universal for the class of definable topological spaces $(X,\tau)$, with $X\subseteq R$, that embed definably into a definable order topological space. 
\end{corollary}

\begin{proof}
Let $(X,\tau)$, with $X\subseteq R$, be a definable topological space that embeds definably into a definable order topological space. Observe first that, since order topological spaces are Hausdorff and regular and these properties are hereditary, then $(X,\tau)$ is Hausdorff and regular.
Now note that, in this case, by Proposition~\ref{prop_n_line_is_not_order_space}, the partition $\XX$ of $X$ described in Lemmas~\ref{remark_for_two_theorems} and~\ref{lemma_proof_two_thems} is such that, for every $A\in \Xinf$, the space $(A^*,\tau_A)$ described in Lemma~\ref{lemma_proof_two_thems} satisfies $\tau_A=\taualex$ if and only if $A^*=I\times \{0\}$, for some interval $I\subseteq X$ (see the cases in the proof of Lemma~\ref{lemma_proof_two_thems}). Following the proof of Theorem~\ref{them_ADC_or_DOTS}, we derive that there exist disjoint definable $\tau$-open sets $Y, Z \subseteq X$, whose union is cofinite in $X$, and some $n_Y$ such that $(Y,\tau)$ embeds definably into $(R\times \{0,\ldots, n_Y\}, \taulex)$ and $(Z,\tau)$ embeds definably into $(R\times \{0\}, \taualex)=(R \times\{0\},\tau_e)$. The proof then concludes in a similar manner to the proof of Corollary~\ref{cor_universal_space}.
\end{proof}

From now on, we let $\CCone$ denote the class of one-dimensional regular Hausdorff definable topological spaces.

\begin{remark}\label{remark:universality_T_3_spaces}
If $\RR$ expands an ordered field then, by Remark~\ref{remark_assumption_X}, the space described in Corollary~\ref{cor_universal_space} is almost definably universal for $\CCone$. 
 
In general, recall that Remark~\ref{remark_general_them_ADC_or DOTS} states that any space in $\CCone$ can be partitioned into finitely many points and open subsets each definably homeomorphic to some set $R\times\{0,\ldots, n-1\}$ with either the $\taulex$ or $\taualex$ topology, for various $n$. Consequently, following the arguments in the proof of Corollary~\ref{cor_universal_space}, one may show that a space given by infinitely many copies of $(R\times [0,1],\taulex)$ and $(R\times[0,\infty), \taualex)$ is almost definably universal for $\CCone$. Such a space exists as a three-dimensional space. That is, consider $(X,\tau)$, where $X=((-\infty, 0) \times R\times[0,1]) \cup ([0,\infty) \times R \times [0,\infty))$. Then let $\tau$ be the topology such that, for every $t\in R$, the fiber of $X$, which is given by either $\{t\}\times R \times [0,1]$ or $\{t\}\times R\times[0,\infty)$, is open and its projection to the last two coordinates is a homeomorphism onto $(R\times [0,1],\taulex)$ or $(R\times[0,\infty), \taualex)$ respectively.

Using Corollary~\ref{cor_universal_space_embed_ordertopspace}, as well as Remark~\ref{remark_general_them_ADC_or DOTS} and Proposition~\ref{prop_n_line_is_not_order_space}, one may analogously identify an almost definably universal space specifically for the class of one-dimensional definable topological spaces that embed definably into a definable order topological space. Such a space can be found that is two-dimensional if $\RR$ expands an ordered field, and three-dimensional in general.
\end{remark}

Note that Proposition~\ref{prop:universal_euclidean_space} states that, when $\RR$ expands an ordered field, the class of euclidean spaces of dimension at most $n$ contains an almost definably universal space for itself (namely $(R^n,\tau_e)$). On the other hand, the space described in Corollary~\ref{cor_universal_space}, which, by Remark~\ref{remark:universality_T_3_spaces}, is almost definably universal for $\CCone$ whenever $\RR$ expands an ordered field, is two-dimensional, and so it does not belong in $\CCone$. In light of these results, and in the context of the classical literature on universal Banach spaces mentioned at the start of the section,  it is natural to ask if, in the case that $\RR$ expands an ordered field, there exists a space $(X,\tau)\in\CCone$ that is almost definably universal for $\CCone$, as well as to ask more generally which classes of spaces admit an almost definably universal space, and when does the space belong in the class. 

We first answer the question regarding the class $\CCone$ negatively and then derive, from Theorem~\ref{them 1.5} and Corollary~\ref{cor_T3_separable_embedding}, positive answers for two other classes of one-dimensional spaces.

\begin{proposition}\label{prop_nonexistence_universal_space}
There does not exist a $T_1$ one-dimensional definable topological space $(X,\tau)$ that is almost definably universal for $\CCone$. 
\end{proposition}

In order to prove the proposition we first introduce a notion of equivalence for curves and some preliminary results. 
\begin{lemma} \label{lem:curve-equiv}
Let $\gamma:(a,b)\rightarrow R^n$ and $\gamma':(a',b')\rightarrow R^n$ be two definable curves, with convergence endpoints $c\in \{a,b\}$ and $c'\in\{a',b'\}$ respectively. For any $s, t \in \exR$, let $I(s,t)$ denote the interval with endpoints $s$ and $t$. The following are equivalent.
\begin{enumerate}[(a)]
     \item \label{itm0:curve-equiv} For any $a < t < b$ and $a' < t' < b'$, it holds that $\gamma(I(c,t)) \cap \gamma(I(c',t')) \neq \emptyset$. 

    \item \label{itm1:curve-equiv} For any $a < t < b$, there exists some $a' < t' < b'$ such that $\gamma'(I(c',t'))\subseteq \gamma(I(c,t))$.

    \item \label{itm2:curve-equiv} For any $a' < t' < b'$, there exists some $a < t < b$ such that $\gamma(I(c,t))\subseteq \gamma'(I(c',t'))$. 

    \item \label{itm3:curve-equiv} For any definable topological space $(X,\tau)$ with $\gamma[(a,b)] \cup \gamma'[(a',b')] \subseteq X\subseteq R^n$, and any $x\in X$, it holds that $\gamma$ $\tau$-converges to $x$ if and only if $\gamma'$ $\tau$-converges to $x$. 
\end{enumerate}
\end{lemma}
\begin{proof}
It is easy to see that $(\ref{itm1:curve-equiv}) \vee (\ref{itm2:curve-equiv}) \Rightarrow (\ref{itm0:curve-equiv})$.
Similarly, one may easily check, using the definition of curve convergence (Definition~\ref{dfn:curve}), that $(\ref{itm1:curve-equiv}) \wedge (\ref{itm2:curve-equiv}) \Rightarrow (\ref{itm3:curve-equiv})$. We show that $(\ref{itm0:curve-equiv})\Rightarrow (\ref{itm1:curve-equiv}) \wedge (\ref{itm2:curve-equiv})$ and $(\ref{itm3:curve-equiv})\Rightarrow(\ref{itm0:curve-equiv})$.

\textbf{Proof of $(\ref{itm0:curve-equiv})\Rightarrow (\ref{itm1:curve-equiv}) \wedge (\ref{itm2:curve-equiv})$.}

By symmetry of the statements (\ref{itm1:curve-equiv}) and (\ref{itm2:curve-equiv}) it suffices to show that $(\ref{itm0:curve-equiv})\Rightarrow(\ref{itm1:curve-equiv})$. Hence suppose that (\ref{itm0:curve-equiv}) holds, and let us fix some $a<t<b$. Consider the definable set $J$ of all $a' < s' < b'$ such that $\gamma'(s')\in \gamma(I(c,t))$. By o-minimality, there exists some $a' < t' < b'$ such that either $I(c',t') \subseteq J$ or $I(c',t') \cap J = \emptyset$. In the second case however we have that $\gamma'(I(c',t')) \cap \gamma(I(c,t)) = \emptyset$, contradicting (\ref{itm0:curve-equiv}). So $I(c',t') \subseteq J$, meaning that $\gamma'(I(c',t'))\subseteq \gamma(I(c,t))$.

\textbf{Proof of $(\ref{itm3:curve-equiv})\Rightarrow(\ref{itm0:curve-equiv})$.}

By contraposition, assume that there exists some $a < t < b$ and $a' < t' < b'$ such that $\gamma(I(c,t)) \cap \gamma'(I(c',t')) = \emptyset$. In particular, we have that $\gamma'(I(c',t')) \subsetneq R^n$. Fix any point $x\in R^n\setminus \gamma'(I(c',t'))$. Consider the definable topology $\tau$ on $R^n$ where every point in $R^n \setminus \{x\}$ is isolated, and a basis of open neighbourhoods for $x$ is given by the definable family of sets $\{ \{x\}\cup\gamma(I(c,s)) : a < s < b\}$. Clearly $\gamma$ $\tau$-converges to $x$ but $\gamma'$ does not, i.e. we reach the negation of (\ref{itm3:curve-equiv}). 
\end{proof}

\begin{definition}\label{dfn:curve_equiv}
Let $(X,\tau)$ be a definable topological space. We say that two definable curves $\gamma:(a,b)\rightarrow X$ and $\gamma':(a',b')\rightarrow X$, with fixed convergence endpoints $c\in\{a',b'\}$ and $c'\in\{a,b\}$ respectively, are \emph{equivalent} if any of the equivalent conditions in Lemma~\ref{lem:curve-equiv} holds. 
\end{definition}

It is easy to check that two injective definable curves in $R$ are equivalent if and only if they $e$-converge to the same point in $\exR$ from the same side, i.e. from left or right (see Remark~\ref{remark_side_convergence}).

\begin{definition}
Let $(X,\tau)$ be a definable topological space.
For any $x\in X$, let $\num(x,X,\tau)$ ($\num(x)$ for short, when the underlying topological space $(X,\tau)$ is clear from the context) denote the maximum cardinality of a  set of non-equivalent definable curves (together with fixed convergence endpoints) in $X$ $\tau$-converging to $x$.

Let $\num(X,\tau)$ ($\num(X)$ for short) be defined to be $\sup \{ \num(x) : x\in X\}$.
\end{definition}

\begin{exmpl}\label{ex:cross}
For any $n$ and $1 \leq i\leq n$, let $L(i)$ denote the line 
\[
\{0\}\times \overset{i-1}{\cdots} \times \{0\}\times R \times \{0\}\times \overset{n-i}{\cdots} \times \{0\},
\]
where $R$ is in the $i$-th coordinate position. Consider the euclidean space $X=\bigcup_{1\leq i \leq n} L(i)\subseteq R^n$. For any $s,t\in \exR$, let $I(s,t)$ denote the interval with endpoints $s$ and $t$. Observe that, by o-minimality, a definable curve $\gamma:(a,b)\rightarrow X$, with convergence endpoint $c\in \{a,b\}$, $e$-converges to the point $\al 0,\ldots, 0 \ar \in R^n$ if and only if there exists some $1\leq i \leq n$ and $a<t<b$ such that $\gamma(I(c,t))\subseteq L(i)$, and moreover the composition $\pi_i \circ \gamma$, where $\pi_i$ denotes the projection to the $i$-th coordinate, $e$-converges (as $t$ tends to $c$) to $0$, from either the left or the right.  
One may derive from this that $\num(\al 0,\ldots, 0 \ar, X, \tau_e)=2n$.  
\end{exmpl}

\begin{remark}\label{rem:finite-n(x)}
If $\dim X \leq 1$ and $(X,\tau)$ is $T_1$, then one may easily check that $\num(x)=1+|R_x|+|L_x|$, for every $x\in X$, using Remarks~\ref{remark_side_convergence} and~\ref{remark_generalizing_Rx_Lx}, as well as the fact that, by $T_1$-ness, for any $x \in X$, there is only one eventually constant definable curve (up to equivalence) that $\tau$-converges to $x$. By Lemma~\ref{lemma_2}, Proposition~\ref{prop_basic_facts_P_x_2}~(\ref{itm1:basic_facts_P_x_2}) and~(\ref{itm2:basic_facts_P_x_2}), and uniform finiteness, it follows that $\num(X) <\omega$.
\end{remark}

\begin{lemma}\label{lemma_universality_index}
Let $f:(X,\tau)\rightarrow (Y,\mu)$, be a continuous injective definable map between definable topological spaces. For any $x\in X$, it holds that $\num(x) \leq \num(f(x))$. In particular, $\num(X) \leq \num (Y)$.   
\end{lemma}
\begin{proof}
By Proposition~\ref{prop_cont_lim}, if $\gamma$ and $\gamma'$ are two definable curves $\tau$-converging to $x$ (with convergence endpoints $c$ and $c'$, respectively) then $f\circ \gamma$ and $f \circ \gamma'$ (with the same respective convergence endpoints, namely $c$ and $c'$) $\mu$-converge to $f(x)$. It therefore suffices to show that if $\gamma$ and $\gamma'$ are non-equivalent then $f \circ \gamma$ and $f\circ \gamma'$ (taken with the corresponding convergence endpoints) are non-equivalent. This however follows easily from the characterization of curve equivalence given by Lemma~\ref{lem:curve-equiv}~(\ref{itm0:curve-equiv}) and the injectivity of $f$.     
\end{proof}

\begin{remark}
Some observations can be derived from Lemma~\ref{lemma_universality_index} regarding the existence of definably universal spaces. Recall that in Example~\ref{ex:cross} we define, for every $n$, a one-dimensional set $X\subseteq R^n$ such that $\num(\al 0,\ldots, 0 \ar, X, \tau_e)=2n$. 
By Lemma~\ref{lemma_universality_index} and Remark~\ref{rem:finite-n(x)}, it follows that there does not exist a one-dimensional $T_1$ definable topological space that is definably universal for the class of all one-dimensional euclidean spaces, and in particular neither for $\CCone$. 
\end{remark}

We may now prove Proposition~\ref{prop_nonexistence_universal_space}.

\begin{proof}[Proof of Proposition~\ref{prop_nonexistence_universal_space}]
Let $Y=R\times\{0,\ldots, n-1\}$ and consider the space $(Y,\taulex)$, which belongs to $\CCone$. If, for every $0\leq i <n$, we identify the subspace $R\times\{i\}$ with $R$ through the projection to the first coordinate (see Remark~\ref{remark_generalizing_Rx_Lx}) then, for any $x\in R$, it holds that $L_{\al x,0\ar}=\{\al x,0 \ar,\ldots, \al x,n-1\ar\}$ and $R_{\al x,0\ar}=\emptyset$. So $\num(\al x,0\ar)=n$, and in fact it can easily be shown that $\num(Y)=n$. Moreover note that, for any cofinite subset $Y'\subseteq Y$, it still holds that $\num(Y')=n$, since we may always find an interval $I\subseteq R$ such that $I\times\{0,\ldots, n-1\}\subseteq Y'$.

Suppose that $(X,\tau)$ is a $T_1$ one-dimensional definable topological space that is almost definably universal for $\CCone$. Since $(X,\tau)$ is $T_1$, we have $\num(X) < \omega$ (see Remark~\ref{rem:finite-n(x)}). However, by Lemma~\ref{lemma_universality_index} and the above observation, we have that $\num(X)\geq n$ for every $n$, a contradiction.
\end{proof}

The same proof would still have worked if we had considered the space $(R\times \{0,\ldots, n-1\}, \taualex)$ in place of $(R\times \{0,\ldots, n-1\}, \taulex)$. Ultimately, one may show that there exists no one-dimensional definable topological space that is almost definably universal for either of the following two classes: all one-dimensional spaces with the $\taulex$ topology and all one-dimensional spaces with the $\taualex$ topology.

It is our belief that Proposition~\ref{prop_nonexistence_universal_space} can likely be improved by dropping the condition of being $T_1$. In other words, we believe that the following question has a negative answer:
\begin{question}
Is there a one-dimensional definable topological space (which is necessarily not $T_1$ by Proposition~\ref{prop_nonexistence_universal_space}) that is almost definably universal for $\CCone$?
\end{question}

As a counterpoint to Proposition~\ref{prop_nonexistence_universal_space}, Theorem~\ref{them_ADC_or_DOTS} -- more specifically, Corollary~\ref{cor_T3_separable_embedding} -- and Theorem~\ref{them 1.5} do  yield the existence of two  classes of one-dimensional spaces, each of which contains a space that is almost definably universal for itself, as shown by the following corollaries (see also Remark~\ref{rem:2dim-univ}). 

Let $\CCCsep$ denote the class of Hausdorff regular definably separable one-dimensional spaces.

\begin{corollary}\label{cor:universality_CCCsep}
The disjoint union of $(R,\tau_e)$ and $(R\times\{0,1\}, \taulex)$ is Hausdorff, regular, definably separable and almost definably universal for the class of Hausdorff, regular, definably separable spaces $(X,\tau)$, where $X\subseteq R$.

It follows that, whenever $\RR$ expands an ordered field, the class $\CCCsep$ contains an almost definably universal space.
\end{corollary}
\begin{proof}
The second paragraph of the corollary follows from the first by direct application of Remark~\ref{remark_assumption_X}. We prove the first paragraph.

Since $(R,\tau_e)$ and $(R\times\{0,1\}, \taulex)$ are regular, Hausdorff and definably separable (see Lemma~\ref{lem:T1-sep} and Proposition~\ref{prop:def-sep_eucl_disc_llt}(\ref{itm1:sep-spaces-examples})) their disjoint union is too. 

By Corollary~\ref{cor_T3_separable_embedding}, it suffices to show that, for any definable set $X\subseteq R\times\{0,1\}$, there exists a cofinite subspace $Y$ of $(X,\taulex)$ that embeds definably into the disjoint union of $(R,\tau_e)$ and $(R\times\{0,1\}, \taulex)$.

We partition $X\subseteq R\times\{0,1\}$ as follows. Let $X_1=\{ \al x, i \ar\in X : \al x, 1-i \ar\notin X \}$ and $X_2=X\setminus X_1$. By o-minimality, there exists a partition $\XX$ of a cofinite subset of $X$ with the following properties. For every $A\in\XX$, there exists an interval $I \subseteq R$ such that either $A=I\times \{i\}$, for some $i\in\{0,1\}$, and $A\subseteq X_1$, or $A=I\times\{0,1\}$ (and so $A\subseteq X_2$). Let $\XX_1=\{ A\in\XX: A\subseteq X_1\}$ and $\XX_2=\XX\setminus \XX_1$. 

Note that every $A\in\XX$ is open in $(X,\taulex)$, and that the subspace topology on $A$ corresponds precisely to the lexicographic order topology on $A$. If $A\subseteq X_1$, then the projection $\al x,i\ar\mapsto x$ is an open embedding $(A, \taulex)\hookrightarrow (R,\tau_e)$, and otherwise the identity is an open embedding $(A,\taulex)\hookrightarrow (R\times\{0,1\},\taulex)$.  Hence the projection to the first coordinate is an open embedding $(\bigcup \XX_1,\taulex)\hookrightarrow (R,\tau_e)$ and the identity is an open embedding $(\bigcup\XX_2,\taulex)\rightarrow (R\times\{0,1\},\taulex)$, which completes the proof. 
\end{proof}

Finally, we consider the class of one-dimensional Hausdorff definable topological spaces satisfying the frontier dimension inequality (\textbf{fdi}; see Definition~\ref{dfn:fdi}), which we denote $\CCCfdi$. By Proposition~\ref{prop_T2_frontier_ineq_regular}, we know that these spaces are regular. We show that, whenever $\RR$ expands an ordered field, the class contains an almost definably universal space. This is a corollary of Theorem~\ref{them 1.5}.

\begin{corollary}\label{cor_universal_space_for_T2_fdi}
The disjoint union of the spaces $(R,\tau_e)$, $(R,\tau_r)$, $(R,\tau_l)$ and $(R,\tauc)$ is Hausdorff, satisfies the frontier dimension inequality, and is almost definably universal for the class of Hausdorff definable topological spaces $(X,\tau)$ with the frontier dimension inequality, where $X\subseteq R$. 

It follows that, whenever $\RR$ expands an ordered field, the class $\CCCfdi$ contains an almost definably universal space.
\end{corollary}
\begin{proof}
As in the proof of Corollary~\ref{cor:universality_CCCsep}, the second paragraph of the corollary follows from the first by direct application of Remark~\ref{remark_assumption_X}. We prove the first paragraph.

Since the spaces $(R,\tau_e)$, $(R,\tau_r)$, $(R,\tau_l)$ and $(R,\tauc)$ are Hausdorff and satisfy the \textbf{fdi}, their disjoint union has these properties too. 

We fix $(X,\tau)$, $X\subseteq R$, a Hausdorff definable space which satisfies the \textbf{fdi}. By Theorem~\ref{them 1.5}, there exists a finite partition $\XX$ of $X$ into points and intervals such that, for each $I\in \XX$, the subspace topology $\tau|_I$ is one of $\tau_e$, $\tau_r$, $\tau_l$ or $\tauc$. Let $\XX'$ be the subfamily of intervals in $\XX$ and let $X'=\bigcup\{ int_\tau I : I\in \XX'\}$. By the \textbf{fdi}, the set $X'$ is cofinite in $X$. Let $X_e=\{ x\in X' : x\in I\in\XX,\, (I,\tau)=(I,\tau_e)\}$. Then the identity is an open embedding $(X_e,\tau)\hookrightarrow (R, \tau_e)$. By repeating this argument with the topologies $\tau_r$, $\tau_l$ and $\tauc$, we may conclude that $X'$ can be partitioned into four definable open subspaces on which the identity is an embedding into one of $(R,\tau_e)$, $(R,\tau_r)$, $(R,\tau_l)$ or $(R,\tauc)$. The corollary follows.
\end{proof}

\begin{remark}\label{rem:2dim-univ}
Following Remarks~\ref{remark:universality_euclidean_spaces} and~\ref{remark:universality_T_3_spaces}, if $\RR$ does not expand an ordered field, then one may adapt the proof of Corollaries~\ref{cor:universality_CCCsep} and~\ref{cor_universal_space_for_T2_fdi} to show the existence of two-dimensional almost definably universal spaces for each of the classes $\CCCsep$ and $\CCCfdi$. More precisely, any space containing infinitely many disjoint copies of $(R,\tau_e)$ and $(R\times\{0,1\}, \taulex)$ is almost definably universal for $\CCCsep$, and any space with infinitely many disjoint copies of the spaces $(R,\tau_e)$, $(R,\tau_r)$, $(R,\tau_l)$ and $(R,\tauc)$ is almost definably universal for $\CCCfdi$.  
\end{remark}

Note that, by Lemma~\ref{lem:T1-sep} and o-minimality, any definable subspace of $(R,\tau_e)$ or $(R\times\{0,1\}, \taulex)$ is definably separable. 
By Corollary~\ref{cor:universality_CCCsep} and Remark~\ref{rem:2dim-univ}, it follows that any definable subspace of a space in $\CCCsep$ is also definably separable. 
In other words, definable separability is a hereditary property for $T_3$ one-dimensional spaces. Since being $T_3$ is also a hereditary property, we have that the class $\CCCsep$ is closed under passing to one-dimensional definable subspaces.

\section{Definable Hausdorff compactifications}\label{section:compactifications}

In this section, we use the decomposition of $T_3$ spaces described in Section~\ref{section: universal spaces} to address the question of which definable topological spaces can be Hausdorff compactified in a definable sense. Our main result (Theorem~\ref{them_compactification}) shows that a Hausdorff definable topological space in the line has a definable Hausdorff compactification of dimension at most one if and only if it is $T_3$. This result then generalizes to all one-dimensional definable topological spaces (Remark~\ref{rmk_compactification}). 

Recall that a definable topological space is definably compact (Definition~\ref{dfn:compact}) if and only if every definable curve in it converges. 
We present the following definition, and prove that it characterizes the one-dimensional $T_3$ spaces that can be definably one-point Hausdorff compactified.

\begin{definition} \label{dfn:near-compact}
A definable topological space $(X,\tau)$, $\dim X\leq 1$, is \emph{definably near-compact} if, up to equivalence, there are only finitely many non-convergent definable curves in $(X,\tau)$.  

Clearly definable compactness implies definable near-compactness. We say that a definable topological space $(X^*,\tau^*)$, $X^*\subseteq R$, is a \emph{definable near-compactification} of $(X,\tau)$ if $(X^*, \tau^*)$ is definably near-compact and there exists a definable embedding $(X,\tau)\hookrightarrow(X^*,\tau^*)$.   
\end{definition}

We have the following characterization of definably near-compact definable topological spaces in the line, 
which is a direct consequence of Remark~\ref{remark_side_convergence}, and can be seen as an analogue of Lemma~\ref{lem:RL-compact}, a statement about definably compact spaces in the line. 

\begin{lemma} \label{lemma_near-compact_characterization}
A definable topological space $(X,\tau)$, where $X\subseteq R$, is definably near-compact if and only if the set $(\bigcup_{x\in X} R_x) \cap (\bigcup_{x\in X} L_x)$ is cofinite in $cl_e X$.
\end{lemma}

\begin{remark}\label{remark:ADC_or_DOTS_compactification}
Note that, for any $n$ and interval $I$, the spaces $(I\times\{0,\ldots,n\},\taulex)$ and $(I\times\{0,\ldots,n\},\taualex)$ are definably near-compact, and furthermore they are definably compact if and only if $I$ is a closed interval. 
It follows that, if $(X, \tau)$, $X \subseteq R$, is a regular Hausdorff definable topological space, then the embedding $h :(Y\cup Z, \tau) \hookrightarrow (Y^*\cup Z^*, \taulex|_{Y^*}\cup \taualex|_{Z^*})$ described in the proof of Theorem~\ref{them_ADC_or_DOTS} is a definable near-compactfication of a cofinite open subspace of $X$.
\end{remark}

We extract the following observation from the proof of Lemma~\ref{lemma_proof_two_thems}, which can be seen as an improvement of Theorem~\ref{them_ADC_or_DOTS} for $T_3$ definably near-compact spaces in the line. 
We will use this result in Section~\ref{section: affine} (Corollary~\ref{thm:compact_affine}).

\begin{lemma} \label{prop_proof_two_thems}
Let $(X,\tau)$ be a regular Hausdorff definable topological space with $X\subseteq R$. Let $\mathcal{X}$ be a finite partition of $X$ as given by Lemma~\ref{remark_for_two_theorems}. For each $A\in \Xinf$, let $A^*$ and $h_A:A\rightarrow A^*$ be as given by Lemma~\ref{lemma_proof_two_thems}. 

If $(X,\tau)$ is definably near-compact, then, for any $A\in\Xinf$, the map $h_A$ is a bijection. In particular, by Lemma~\ref{lemma_proof_two_thems}(\ref{itm2:lemma_proof_two_theorems}), it is a definable homeomorphism $(A,\tau)\rightarrow (A^*,\mu)$, where $\mu$ is one of $\taulex$ or $\taualex$.
\end{lemma}
\begin{proof}
Recall the proof by cases of Lemma~\ref{lemma_proof_two_thems}. Let $A = \bigcup_{0 \leq i < n}f_{i}(I)$. Observe that, by Lemma~\ref{lemma_near-compact_characterization}, if $(X,\tau)$ is definably near-compact then, for all but finitely many $x\in I$, it holds that $x\in R_{y} \cap L_{z}$, where, by Lemma~\ref{properties_cofnite_subset_1.5}(\ref{itm:properties_cofnite_subset_1.5_b}), $y,z\in [x]^E\subseteq \{x, f_{n-1}(x)\}$. It follows that, if $(X,\tau)$ is definably near-compact, Cases 0, 1 and 2  in the proof of Lemma~\ref{lemma_proof_two_thems} are not possible. In the remaining cases, the function $h_A$ defined is a bijection.
\end{proof}

In the context of Remark~\ref{remark_general_them_ADC_or DOTS}, Lemma~\ref{prop_proof_two_thems} generalizes in the natural way to all $T_3$ one-dimensional spaces. 

We now show that if a $T_3$ one-dimensional space is definably near-compact then it can be definably one-point Hausdorff compactified.
The converse implication can also be derived, from the observation that, given a Hausdorff definably compact definable topological space $(\Xopc,\tauopc)$, with $\dim \Xopc = 1$, and any point $x\in \Xopc$, by Lemmas~\ref{lemma_2} and~\ref{lemma_near-compact_characterization} and Remark~\ref{remark_generalizing_Rx_Lx} the subspace $(\Xopc\setminus\{x\}, \tauopc)$ is definably near-compact. 
We will use the following result to show that $T_3$ spaces in the line can be definably Hausdorff compactified (Theorem~\ref{them_compactification}).

\begin{proposition}\label{lemma_one_point_compact}
Let $(X,\tau)$, $\dim X \leq 1$, be a regular Hausdorff definably near-compact definable topological space. 
Then there exists a Hausdorff definably compact definable topological space $(\Xopc, \tauopc)$ and a definable embedding $h:(X,\tau)\hookrightarrow (\Xopc,\tauopc)$, where $\Xopc \setminus h(X)$ is a singleton. 
\end{proposition}
If $\RR$ expands the field of reals then we leave it to the reader to check that $(\Xopc,\tauopc)$ is the classical one-point compactification of $(X,\tau)$.  
\begin{proof}[Proof of Proposition~\ref{lemma_one_point_compact}]
We prove the lemma in the case where $X\subseteq R$. Given the assumptions in Remark~\ref{remark_generalizing_Rx_Lx}, the proof adapts to a proof of the general case. 

Let $c=\al 0,1 \ar \in R^2$ and let $\Xopc=(X\times \{0\}) \cup \{c\}$. Let $h:X \to \Xopc$ be given by $x\mapsto \al x,0 \ar$,  and let $\tau_h$ be the push-forward topology of $\tau$ by $h$ (see Definition~\ref{dfn:push-forward}). We will define $\tauopc$ as an extension of $\tau_h$ to a topology on $\Xopc$. If $\tau$ is definably compact then it clearly it suffices to let $\tauopc = \tau_h \cup \{c\}$, and so we assume otherwise.

Set $R_c:=\{ x\in \exR\setminus \bigcup_{x\in X} R_x : \exists \, y>x \; (x,y)\subseteq X\}$ and $L_c:=\{ x\in \exR\setminus \bigcup_{x\in X} L_x : \exists \, y<x \; (y,x)\subseteq X\}$. Set $\E_c:= R_c \cup L_c$. Since $(X,\tau)$ is definably near-compact, $\E_c$ is finite, by Lemma~\ref{lemma_near-compact_characterization}.  
Since, by assumption, $(X,\tau)$ is not definably compact, we also have that $\E_c\neq \emptyset$. Let $R_c=\{y_1,\ldots, y_n\}$ and $L_c=\{z_1,\ldots, z_m\}$, and let $\UU(c)$ be the family of sets 
\[
U=\bigcup_{1\leq i \leq n} (y_i,y'_i) \cup \bigcup_{1 \leq j \leq m} (z'_j, z_j)
\] 
definable uniformly over those parameters $(y'_1, \ldots, y'_n, z'_1,\ldots, z'_m)\in R^{n+m}$ for which $y_i<y'_i$, $z'_j<z_j$ and moreover $U\subseteq X$. 

Let $\tauopc$ be the definable topology with basis $\{ (int_\tau U \times \{0\}) \cup \{c\}: U\in\UU(c)\} \cup \tau_h$. It is routine to check that this is a well-defined topology and that $h:(X,\tau)\hookrightarrow (\Xopc, \tauopc)$ is an embedding. Since $(X,\tau)$ is Hausdorff, by definition of $\UU(c)$ and Lemma~\ref{lemma_explaining_neighbourhoods} it is immediate that $(\Xopc,\tauopc)$ is also Hausdorff. It remains to prove that it is definably compact. 

Let $\gamma'$ be a definable curve in $(\Xopc,\tauopc)$. We may assume that $\gamma'$ is injective and hence lies in $X\times\{0\}$. Let $\gamma= h^{-1}\circ \gamma'$. Let $x_0\in \exR$ denote the limit of $\gamma$ in the euclidean topology. Since the remaining case is analogous, we consider only the case where $\gamma$ $e$-converges to $x_0$ from the right. Then clearly there must exist $y>x_0$ such that $(x_0,y)\subseteq X$. If $x_0\notin R_c$, then $x_0\in\bigcup_{x\in X} R_x$ and so, by Remark~\ref{remark_side_convergence}, $\gamma$ $\tau$-converges to some $x\in X$, and it follows that $\gamma'$ $\tauopc$-converges to $h(x)$. 

Therefore, it remains to consider the case $x_0\in R_c$. We will show that $\gamma'$ $\tauopc$-converges to $c$. To prove this, it suffices to show that, for every $U\in\UU(c)$, there is $x_U>x_0$ such that $(x_0,x_U)\subseteq int_\tau U$. 

Towards a contradiction, suppose otherwise. Then, by o-minimality, there exists $U_1\in\UU(c)$ and $x_1>x_0$ such that $(x_0,x_1)\cap int_\tau U_1=\emptyset$. By definition of $\UU(c)$, we may moreover take $x_1$ close enough to $x_0$ to satisfy that $(x_0,x_1)\subseteq U_1$. For every $x_0<x<x_1$, we have that $x\in \partial_\tau (X\setminus U_1)$, and so, from Proposition~\ref{basic_facts_P_x_1}(\ref{itm: basic_facts_4}), it follows that $\E_x \setminus U_1 \neq \emptyset$, as $U_1$ is $e$-open. Let $f:(x_0,x_1)\rightarrow \exR$ be the definable map given by $x\mapsto \min \E_x \setminus U_1$. By Hausdorffness (Lemma~\ref{lemma:basic_facts_P_x_2}(\ref{itmb:lemma_basic_facts_P_x_2})) and o-minimality, this function is \econtinuous and strictly monotone on some subinterval $(x_0,x_2)\subseteq (x_0,x_1)$. Let $y_0=\elim_{x\rightarrow x_0} f(x)$. If $f$ is increasing on $(x_0,x_2)$, then, by construction of $\UU(c)$ and the fact that $f$ maps into $\exR\setminus U_1$, it cannot be that $y_0\in R_c$. However, there clearly exists $y'\in R$ such that $(y_0,y')\subseteq X$. So there exists $y\in X$ such that $y_0\in R_y$ and, by Lemma~\ref{lemma_f_Rx_Lx}(\ref{itmc:lemma_f_Rx_Lx}) and regularity, it follows that $x_0\in R_y$, a contradiction since $x_0\in R_c$. The case where $f$ is decreasing is analogous.
\end{proof}

We now present the main result of this section. 

\begin{theorem}\label{them_compactification}
Let $(X,\tau)$, $X\subseteq R$, be a Hausdorff definable topological space. Then $(X,\tau)$ is regular if and only if there exists a definably compact Hausdorff definable topological space $(\Xcp,\taucp)$, with $\dim \Xcp \leq 1$, and a definable embedding $(X,\tau)\hookrightarrow (\Xcp,\taucp)$. 
\end{theorem}

The \textquotedblleft if" direction of the theorem is proven largely by the following lemma. 

\begin{lemma}\label{lemma_Hausdorff_compact_is_regular}
Let $(X,\tau)$, $\dim X\leq 1$, be a definably compact Hausdorff definable topological space. Then $(X,\tau)$ is regular. 
\end{lemma}
\begin{proof}
Let $(X,\tau)$ be as in the lemma and towards a contradiction suppose that it is not regular. Let $x\in X$ and let $C\subseteq X$ be a $\tau$-closed set such that $x\notin C$ and $cl_\tau (A)\cap C\neq \emptyset$, for every $\tau$-neighbourhood $A$ of $x$. By passing to a larger set if necessary (i.e. by passing if necessary to the complement of a definable $\tau$-neighbourhood of $x$ contained in $X\setminus C$), we may assume that $C$ is definable. Let $\UU$ denote a definable basis of $\tau$-neighbourhoods of $x$. Note that $\{cl_\tau (U) \cap C : U\in \UU\}$ is a definable downward directed family of non-empty sets of dimension at most one, so, by Remark~\ref{remark_curve_selection}, there exists a definable curve $\gamma: (a,b) \rightarrow C$ that is cofinal for this family. By definable compactness, $\gamma$ $\tau$-converges to some point $y\in C\subseteq X\setminus\{x\}$. By definition of $\gamma$, it holds that $y\in cl_\tau (A)$ for every $\tau$-neighbourhood $A$ of $x$. So $x$ and $y$ cannot be separated by $\tau$-neighbourhoods, which contradicts that $(X,\tau)$ is Hausdorff. 
\end{proof}

In light of the above lemma, the question arises of whether or not definably compact Hausdorff spaces of any dimension are regular. Using~\cite[Corollary 25]{atw1}, one may show that the answer is positive whenever $\RR$ expands an ordered field. On the other hand, it is easy to prove that if a Hausdorff definable topological space is definably compact in the sense of condition~\eqref{dfn:directed-compact} (i.e. every downward directed definable family of non-empty closed sets has non-empty intersection), then it is regular (see~\cite[Lemma 5.4.7]{andujar_thesis} for a proof). 
Furthermore, as noted earlier in Remark~\ref{rem:compactness}, the first author proved in~\cite{ag_FTT} that condition~\eqref{dfn:directed-compact} is equivalent to definable compactness for all Hausdorff definable topological spaces. Consequently, Lemma~\ref{lemma_Hausdorff_compact_is_regular} can be generalized to spaces of any dimension. 

We now prove Theorem~\ref{them_compactification}. 
\begin{proof}[Proof Theorem~\ref{them_compactification}]
Let $(X,\tau)$ be a Hausdorff definable topological space with $X \subseteq R$. Since the finite case is trivial, we assume that $\dim X=1$. The \textquotedblleft if" implication of the theorem follows directly from Lemma~\ref{lemma_Hausdorff_compact_is_regular} using the observation that regularity is a hereditary property. We prove the  \textquotedblleft only if" implication. 
 
Assume that $(X,\tau)$ is regular.  
We will make use of Lemmas~\ref{remark_for_two_theorems} and~\ref{lemma_proof_two_thems} and Proposition~\ref{lemma_one_point_compact} to construct a one-dimensional definable Hausdorff compactification for $(X,\tau)$. 

Recall the embedding $(Y\cup Z, \tau) \hookrightarrow (Y^*\cup Z^*, \taulex|_{Y^*}\cup \taualex|_{Z^*})$ described in the proof of Theorem~\ref{them_ADC_or_DOTS}. As noted in Remark~\ref{remark:ADC_or_DOTS_compactification}, this embedding is a definable near-compactfication of a cofinite $\tau$-open subspace of $X$. The idea of the current proof is to extend this embedding to an embedding of $(X,\tau)$ into a regular Hausdorff definably near-compact space $(X^*,\tau^*)$, with $X^*\subseteq R\times\{0,1,\ldots\}$. Applying Proposition~\ref{lemma_one_point_compact} then completes the proof.

Let $\XX=\Xinf\cup \Xsing$ be a finite partition of $X$ as given by Lemma~\ref{remark_for_two_theorems}. For each $A\in\Xinf$, let $A^*$, $\tau_A$ and $h_A$ be as given by Lemma~\ref{lemma_proof_two_thems}. In particular, recall that each set $A^*$ is of the form $I\times\{0,\ldots,n\}$, for some interval $I\subseteq A$ and some $n$, and that $\tau_A$ is either the $\taulex$ or the $\taualex$ topology on $A^*$. Moreover, $h_A:(A,\tau)\hookrightarrow (A^*,\tau_A)$ is a definable embedding.

Let $X^*=\bigcup\{A^* : A\in\Xinf\} \cup \{\al x,0\ar : \{x\}\in \Xsing\}$ and $h=\bigcup\{ h_A : A\in \Xinf\} \cup h'$, where $h'$ is the map with domain $\Asing:= \bigcup\Xsing$ given by $x\mapsto \al x,0 \ar$. Note that $h$ is injective. 

We construct a regular Hausdorff topology $\tau^*$ on $X^*$ such that, for every $A\in\Xinf$, $A^*$ is $\tau^*$-open and $(A^*,\tau^*)=(A^*,\tau_A)$. Since every space $(A^*,\tau_A)$ is definably near-compact, it follows that $(X^*,\tau^*)$ is definably near-compact. We then prove that $h:(X,\tau)\hookrightarrow (X^*,\tau^*)$ is an embedding. 

Let $s=|\Xinf|$. We define $\tau^*$ as follows. For every $x\in \Asing$ and $A\in\Xinf$, we first construct a downward directed definable family $\BB_A(x)$ of $\tau_A$-open subsets of $A^*$. Then we use this to define, for each $x \in \Asing$,
\[
\BB(x):=\{ \{\al x,0 \ar\}\cup V_1\cup\cdots \cup V_s : (V_1,\ldots, V_s)\in \prod_{A\in\Xinf} \BB_A(x) \}.
\]  
It is then routine to check that the family $\bigcup\{ \BB(x) : x\in\Asing\} \cup \bigcup\{ \tau_A : A\in \Xinf\}$ is a basis for a topology $\tau^*$ on $X^*$, which will clearly satisfy that $\tau^*|_{A^*}=\tau_A$, for every $A\in\Xinf$. Since $\Asing$ is finite and the topologies $\tau_A$ are definable, $\tau^*$ is also definable. It will then remain to check that $(X^*,\tau^*)$ is Hausdorff and regular, and that $h \colon (X,\tau) \hookrightarrow (X^*,\tau^*)$ is an embedding.

We fix $x\in \Asing$ and $A\in\Xinf$ and describe $\BB_A(x)$. Recall the notation $A=\bigcup_{0\leq i < n} f_i(I)$ from Lemma~\ref{remark_for_two_theorems}, and let $I=(a,b)$. 
We define families of sets $V_a(x,y)$ and $V_b(x,y)$, as $y \in I$ varies, as follows, according to whether or not $a \in R_x$ and whether or not $b \in L_x$:
\begin{equation*}
V_a(x,y)=
\begin{cases}
((a, y)\times R) \cap A^* &\text{for all } y \in I, \text{ if } a\in R_x,\\
\emptyset &\text{for all } y \in I, \text{ if } a\notin R_x,
\end{cases}
\end{equation*}
\begin{equation*}
V_b(x,y)=
\begin{cases}
((y, b)\times R) \cap A^* &\text{for all } y \in I, \text{ if } b\in L_x,\\
\emptyset &\text{for all } y \in I, \text{ if } b\notin L_x.
\end{cases}
\end{equation*}

Note that these sets are always open in $(A^*,\tau_A)$. Then $\BB_A(x)$ is defined to be the family of sets 
\[
V_A(x,y,z)=V_a(x,y) \cup V_b(x,z), 
\]
definable uniformly in $y$ and $z$ with $a<y<z<b$.
Clearly, $\BB_A(x)$ is a definable downward directed family of $\tau_A$-open subsets of $A^*$. 

We now consider the induced topology $\tau^*$ as described above. Since $\bigcap_{a<y<z<b} V_A(x,y,z)=\emptyset$, for each $A \in \Xinf$, it is immediate from the definition that $\tau^*$ is $T_1$. 
We now show that $(X^*,\tau^*)$ is regular. It will follow (since in a $T_1$ topological space singletons are closed) that it is also Hausdorff.

Consider the sets $V_A(x,y,z)$, as $x \in \Asing$ varies, for some fixed $A \in \Xinf$ and fixed $a<y<z<b$. By Hausdorffness of $\tau$, for any two distinct $x, x'\in \Asing$, if $a\in R_x$, then $a\notin R_{x'}$ and, if $b\in L_x$, then $b\notin L_{x'}$, so $V_A(x,y,z)\cap V_A(x',y,z)=\emptyset$, from where it follows that, for all $x \in \Asing$, $cl_{\tau^*} V_A(x,y,z)\subseteq \{\al x,0\ar\}\cup A^*$, and consequently $cl_{\tau^*} V_A(x,y,z)=\{\al x,0\ar\}\cup cl_{\tau_A} V_A(x,y,z)$. Moreover, note that, since $\tau_A$ is one of $\taulex$ or $\taualex$, it holds, for all $x \in \Asing$ and $y,z \in I$, that $cl_{\tau_A} V_a(x,y)\subseteq ((a,y]\times R)\cap A^*$ and $cl_{\tau_A} V_b(x,z)\subseteq ([z,b)\times R) \cap A^*$. So, for any $x \in \Asing$, $a<y'<y$ and $z<z'<b$, we have that $cl_{\tau_A} V_A(x,y',z') \subseteq V_A(x,y,z)$. It follows that, for all $x \in \Asing$ and $U\in\BB(x)$, there exists $U'\in \BB(x)$ such that $cl_{\tau^*} U' \subseteq U$. Since each $A^*$ is $\tau^*$-open with $(A^*,\tau^*)=(A^*,\tau_A)$, it is easy to check that the same property holds among $\tau^*$-neighbourhoods of points in $A^*$. 
It follows that the topology $\tau^*$ is regular. 

It remains to show that $h:(X,\tau)\hookrightarrow (X^*,\tau^*)$ is an embedding. We fix $x\in X$ and show continuity of $h$ and $h^{-1}$ at $x$ and $h(x)$ respectively. Since, for every $A\in\Xinf$, $h_A:(A,\tau)\hookrightarrow (A^*,\tau_A)$ is an embedding between open subsets of $(X, \tau)$ and $(X^*, \tau^*)$ respectively, where $\tau^*|_{A^*}=\tau_A$, this holds whenever $x\in A$, for some $A\in\Xinf$, so we may assume that $x\in \Asing$. We make use of Proposition~\ref{prop_cont_lim}. 

Fix $\gamma$, an injective definable curve in $X$, and set $\gamma':=h\circ \gamma$. 
We may assume that there is some fixed $A=\bigcup_{0\leq i<n} f_i(I) \in  \Xinf$ and $0\leq j <n$ such that $\gamma$ is contained in the interval $f_j(I)$. Recall from Lemma~\ref{remark_for_two_theorems} that, for every $0\leq i <n$, $f_i : I \rightarrow R$ is a definable \econtinuous strictly monotonic function. We will prove the case where $f_j$ is increasing. The decreasing case is analogous. Let $I=(a,b)$ and  $f_j(I)=(a_j,b_j)$.  We require the following simple fact that follows, for each $x \in \Asing$, from the definition of the definable families $\BB_A(x)$ (recall that $\pi:R^2\rightarrow R$ denotes the projection to the first coordinate). 
\begin{equation}\label{fact_proof_compactification_theorem_alt}
\parbox{0.9\textwidth}{The curve $\gamma'$ $\tau^*$-converges to $h(x)=\al x,0 \ar$ if and only if either $a\in R_x$ and $\pi\circ\gamma'$ $e$-converges to $a$, or $b\in L_x$ and $\pi\circ \gamma'$ $e$-converges to $b$.}
\end{equation} 

Suppose that $\gamma$ $\tau$-converges to $x$. 
Recall that, by Lemmas~\ref{properties_cofnite_subset_1.5}(\ref{itm:properties_cofnite_subset_1.5_b}) and~\ref{remark_for_two_theorems}, if $A\cap E_x \neq \emptyset$, then $x\in A$, a contradiction. So, by o-minimality and Proposition~\ref{basic_facts_P_x_1}(\ref{itm: basic_facts_3}), $\gamma$ must $e$-converge to either $a_j$ or $b_j$. Suppose that $\gamma$ $e$-converges to $a_j$, in which case $a_j\in R_x$ (see Remark~\ref{remark_side_convergence}). Since $f_j$ is increasing, we have that $f_j^{-1}\circ\gamma$ $e$-converges to $a$. By regularity of $\tau$ and Lemmas~\ref{lemma_f_Rx_Lx} and~\ref{remark_for_two_theorems}, it follows that $a \in R_x$. Now note that, by Lemma~\ref{lemma_proof_two_thems}(\ref{itm1.5:lemma_proof_two_theorems}), $\pi\circ \gamma'= \pi \circ h \circ \gamma =f^{-1}_j \circ \gamma$. Hence $\pi\circ\gamma'$ $e$-converges to $a$. So, by (\ref{fact_proof_compactification_theorem_alt}), we conclude that $\gamma'$ $\tau^*$-converges to $h(x)=\al x,0\ar$. Analogously, if $\gamma$ $e$-converges to $b_j$, then $b_j\in L_x$ and, again by regularity of $\tau$ and Lemmas~\ref{lemma_f_Rx_Lx} and~\ref{remark_for_two_theorems}, $b\in L_x$. Moreover, again by Lemma~\ref{lemma_proof_two_thems}(\ref{itm1.5:lemma_proof_two_theorems}), $\pi \circ \gamma'$ $e$-converges to $b$ and so, by (\ref{fact_proof_compactification_theorem_alt}), $\gamma'$ $\tau^*$-converges to $h(x)=\al x,0 \ar$ in this case as well. Hence, by Proposition~\ref{prop_cont_lim}, $h$ is continuous at $x$.

Now suppose that $\gamma'$ $\tau^*$-converges to $h(x)$. By (\ref{fact_proof_compactification_theorem_alt}), there are two cases to consider: either $a \in R_x$ and $\pi \circ \gamma'$ $e$-converges to $a$, or $b \in L_x$ and $\pi \circ \gamma'$ $e$-converges to $b$. In the first case, since $f_j$ is increasing, by Lemmas~\ref{lemma_f_Rx_Lx} and~\ref{remark_for_two_theorems} we have that $a_j\in R_x$. Moreover, by Lemma~\ref{lemma_proof_two_thems}(\ref{itm1.5:lemma_proof_two_theorems}) and since $f_j$ is increasing, $\gamma=h^{-1} \circ\gamma'=f_j \circ \pi \circ \gamma'$ $e$-converges to $a_j$. We conclude that $h^{-1}\circ \gamma'$ $\tau$-converges to $x$, by Remark~\ref{remark_side_convergence}. 
In the other case it can analogously be shown that $h^{-1}\circ \gamma'$ $\tau$-converges to $x$. Hence $h^{-1}$ is continuous at $h(x)$. This completes the proof of the theorem.
\end{proof}

\begin{remark} \label{rmk_compactification} Theorem~\ref{them_compactification} may be generalized to all Hausdorff definable topological spaces of dimension at most one. In particular, using Remark~\ref{remark_general_them_ADC_or DOTS}, one may note that every $T_3$ one-dimensional space $(X,\tau)$ has a cofinite subspace that embeds definably into a finite disjoint union of spaces of the form $(I\times\{0,\ldots\}, \taulex)$ or $(I\times\{0,\ldots\}, \taualex)$, where $I\subseteq R$ is an interval and $m\geq 0$. Then, using a construction similar to the one in the proof of Theorem~\ref{them_compactification}, one may expand this disjoint union, by adding finitely many points, to a definably near-compact space $(X^*,\tau^*)$, and extend the embedding of a cofinite subspace of $(X,\tau)$ into an embedding $(X,\tau) \hookrightarrow (X^*,\tau^*)$. Finally, by applying Proposition~\ref{lemma_one_point_compact} to $(X^*,\tau^*)$, one reaches a definable compactification of $(X,\tau)$.
\end{remark}

\section{Affine topologies}\label{section: affine}

Unless stated otherwise, throughout this section we assume that $\RR$ expands an ordered field. 
We restrict to this setting to be able to assume, using Remark~\ref{remark_assumption_X}, that, up to definable homeomorphism, every definable topological space $(X,\tau)$ satisfies that $X$ is a bounded set. We will also use in this section (in the proof of Lemma~\ref{lemma_walsberg_homeomorphism}) facts about o-minimal expansions of ordered fields from~\cite[Chapter 10]{dries98}.

We wish to classify the one-dimensional definable topologies that are, up to definable homeomorphism, euclidean. We will refer to these topologies as \emph{affine}\footnote{In the o-minimal context, euclidean and affine are used by different authors to refer to the same canonical topology on definable sets. For clarity, we use only the latter to denote these topologies up to definable homeomorphism.}. Our main result is the following, which will utilize results of the previous sections, in particular the decomposition of Hausdorff one-dimensional definable topological spaces (Corollary~\ref{cor 1.5}) and the Hausdorff compactification of certain one-dimensional $T_3$ definable topological spaces (Proposition~\ref{lemma_one_point_compact}). 
Later in the section we will use all of these results to address in our setting several questions of Fremlin and Gruenhage on the nature of perfectly normal, compact, Hausdorff spaces (see Subsection~\ref{subsection: Fremlin}).

\begin{theorem}\label{them_main_2}
Suppose that $\RR$ expands an ordered field. Let $(X,\tau)$, $\dim X\leq 1$, be a Hausdorff definable topological space. Exactly one of the following holds.
\begin{enumerate}[(1)]
\item $(X,\tau)$ contains a subspace definably homeomorphic to an interval with either the discrete or the right half-open interval topology. \label{itm:them_main_2_1}
\item $(X,\tau)$ is definably homeomorphic to a euclidean space (i.e. $(X,\tau)$ is affine). \label{itm:them_main_2_2}
\end{enumerate}
\end{theorem}

Note that, since the map $x\mapsto -x$ is a homeomorphism $(R,\tau_r)\rightarrow (R,\tau_l)$, Theorem~\ref{them_main_2} still holds if we replace the right half-open interval topology by the left half-open interval topology in~(\ref{itm:them_main_2_1}).

The main theorem (Theorem 7.1) in~\cite{walsberg15} states that if a definable metric space contains no infinite definable discrete subspace (equivalently, by Lemma~\ref{remark_def_sep}, if it is definably separable), then it is affine. Hence, taking into account this result and the fact that the euclidean topology $\tau_e$ is definably separable (Proposition~\ref{prop:def-sep_eucl_disc_llt}(\ref{itm1:sep-spaces-examples})) and definably metrizable, statement~(\ref{itm:them_main_2_2}) in Theorem~\ref{them_main_2} can be changed to \textquotedblleft $(X,\tau)$ is definably separable and definably metrizable".

\begin{remark}\label{rem:(3)implies(1)}
Recall that, for any interval $I\subseteq R$, the space $(I,\mu)$, where $\mu\in \{\tau_r, \tauc\}$, is totally definably disconnected (i.e. singletons are the only definably connected non-empty subspaces). On the other hand, by o-minimality, every euclidean space has finitely many definably connected components. Hence Theorem~\ref{them_main_2} implies that a Hausdorff one-dimensional definable topological space $(X,\tau)$ is affine if and only if every definable topological subspace has finitely many definably connected components (or is not totally definably disconnected). 
\end{remark}

A definition, a remark and a lemma precede the proof of Theorem~\ref{them_main_2}.

\begin{definition}~\label{dfn:cell-wise_euclidean}
We say that a definable topological space $(X,\tau)$ is \emph{cell-wise euclidean} if there is a finite partition $\XX$ of $X$ into cells such that, for each $C\in \XX$, $(C,\tau)=(C,\tau_e)$.
\end{definition}

By o-minimal cell decomposition, we can clearly relax the requirement in Definition~\ref{dfn:cell-wise_euclidean} that the sets in $\XX$ be cells to just that they be any definable subset of $X$. 
Furthermore, since it also follows by o-minimal cell decomposition that every definable bijection is a finite union of disjoint definable $e$-homeomorphisms, the property of being cell-wise euclidean is maintained by definable homeomorphism, and is equivalent to being cell-wise affine (i.e. the property that $X$ admits a finite partition into cells on which the subspace topology is affine). 
In particular, any affine space is cell-wise euclidean. Theorem~\ref{them_main_2} implies that the converse holds for one-dimensional Hausdorff definable topological spaces. 
This statement cannot be generalized to definable topological spaces of all dimensions, as illustrated by Appendix~\ref{section:examples}, Example~\ref{example_space_cell-wise_euclidean_not_metrizable}, which describes a two-dimensional definable topological space that is 
$T_3$ and cell-wise euclidean but not even definably metrizable (and hence not affine). Moreover, in Section~\ref{section:examples}, Example~\ref{example: T_1 fdi not regular}, we produce a cell-wise euclidean one-dimensional definable topological space that is not Hausdorff, and in particular not affine.

\begin{remark}\label{remark_cell-wise_euclidean}
By Corollary~\ref{cor 1.5} (and the fact that the map $x\mapsto -x$ is a definable homeomorphism $(R,\tau_l)\rightarrow (R,\tau_r)$), if a Hausdorff definable topological space $(X,\tau)$, with $\dim X \leq 1$, does not have a subspace that is definably homeomorphic to an interval with either the $\tau_r$ or the $\tauc$ topology, then $(X,\tau)$ is cell-wise euclidean.  
\end{remark}

By the above remark, in order to prove that the negation of \eqref{itm:them_main_2_1} implies \eqref{itm:them_main_2_2} in Theorem~\ref{them_main_2}, it suffices to show that if a Hausdorff one-dimensional space is cell-wise euclidean then it is affine. 

The following lemma is essentially Lemma 5.7 in \cite{walsberg15}, proved for definable metric spaces, which we extend to one-dimensional spaces using Lemma~\ref{lemma:compact_homeomorphism} and results of van den Dries \cite{dries98}, one of which requires the setting of an o-minimal expansion of an ordered field. Recall that a euclidean space is definably compact if and only if it is closed and bounded. 

\begin{lemma}\label{lemma_walsberg_homeomorphism}
Suppose that $\RR$ expands an ordered field. 
Let $(X,\tau)$ be a definably compact Hausdorff definable topological space. Let $(Y,\tau_e)$ be a definably compact euclidean space that admits a definable continuous surjection $f:(Y,\tau_e) \rightarrow (X,\tau)$. Then there exists a definable set $Z$ and a definable homeomorphism $(Z,\tau_e)\rightarrow (X,\tau)$. 
\end{lemma}
\begin{proof}
Let $E$ be the kernel of $f$, namely $E=\{\al x,y \ar\in Y^2 : f(x)=f(y)\}$. By continuity of $f$, $E$ is closed in $Y^2$, and so definably compact. By~\cite[Chapter 10, Corollary 2.16]{dries98}, there exists a definable set $Z$ and a definable quotient map $g:(Y,\tau_e)\rightarrow (Z,\tau_e)$ of $E$, i.e. $g$ has kernel $E$, is surjective, continuous and, for every $C\subseteq Z$, if $g^{-1}(C)$ is closed in $(Y,\tau_e)$, then $C$ is closed in $(Z,\tau_e)$. Moreover, by~\cite[Chapter 6, Proposition 1.10]{dries98}, the space $(Z,\tau_e)$ is definably compact.

The definable map $h:(Z,\tau_e)\rightarrow (X,\tau)$ given by $h(g(x))= f(x)$ is well-defined, and clearly continuous and bijective. By Lemma~\ref{lemma:compact_homeomorphism}, it is a homeomorphism. 
\end{proof}

We may now prove Theorem~\ref{them_main_2}.

\begin{proof}[Proof of Theorem~\ref{them_main_2}]
Let $(X,\tau)$ be a Hausdorff definable topological space. The case where $\dim(X)=0$ is trivial (as~(\ref{itm:them_main_2_2}) trivially holds in this case) and so we assume that $\dim(X)=1$. By Remark~\ref{remark_assumption_X}, we may assume that $X\subseteq R$ and is bounded. Note that, by Remark~\ref{remark_no_homeomorphism}, $(X,\tau)$ cannot be both cell-wise euclidean and have a definable copy of an interval with either the $\tau_r$ or the $\tauc$ topology, so (\ref{itm:them_main_2_1}) and (\ref{itm:them_main_2_2}) in the statement of the theorem are mutually exclusive. Applying Remark~\ref{remark_cell-wise_euclidean}, we assume that $(X,\tau)$ is cell-wise euclidean and derive that it is affine.

Since $(X,\tau)$ is cell-wise euclidean, it is also definably near-compact so, by Proposition~\ref{lemma_one_point_compact}, by passing to $(\Xopc,\tauopc)$ if necessary, we may assume that $(X,\tau)$ is definably compact. 

Let $\XX$ be a partition of $X$ into points and intervals such that, for each $C \in \XX$, the subspace $(C,\tau)$ is euclidean. We define, for each $C\in \XX$, a continuous function $f_C: (cl_e C,\tau_e)  \rightarrow (cl_\tau C, \tau)$ extending the identity on $C$. 
Once we have defined these functions, we complete the proof as follows.

Let $num:\XX \rightarrow \omega$ be an enumeration of the elements in $\XX$ and let $Y=\bigcup_{C\in\XX} (cl_{e} C  \times \{num(C)\})$ be the disjoint union of the euclidean closures of the sets in $\XX$. Clearly, $(Y,\tau_e)$ is definably compact. Let $f: (Y,\tau_e) \rightarrow (X,\tau)$ be the function given by $f(x,num(C))=f_C(x)$, where $x\in cl_e C$. This function is clearly definable, surjective and continuous. The result then follows by Lemma~\ref{lemma_walsberg_homeomorphism}. 

It remains to define, for each $C\in \XX$, the function $f_C$.   
If $C\in \XX$ is a singleton, let $f_C$ simply be the identity. Now let us fix an interval $C=I=(a_I,b_I)\in \XX$. By Hausdorffness (Proposition~\ref{prop_basic_facts_P_x_2}(\ref{itm3:basic_facts_P_x_2})) and definable compactness (Lemma~\ref{lem:RL-compact}), there exists a unique point $x_I\in X$ such that $a_I\in R_{x_I}$ and similarly a unique point $y_I\in X$ such that $b_I\in L_{y_I}$. Note that, since $(I,\tau)=(I,\tau_e)$, the points $x_I$ and $y_I$ do not belong in $I$. Let $f_I$ be defined as 
\[
f_I|_I = id, \, f(a_I)=x_I \text{ and } f(b_I)=y_I.
\]
It is routine to check that $f_I$ is continuous as a map $([a_I,b_I], \tau_e) \rightarrow (cl_\tau I,\tau)$ 
\end{proof}

In Appendix \ref{section:examples}, Example~\ref{example_line-wise_euclidean_not_euclidean} we describe a Hausdorff definable topological space of dimension two that has no definable copy of an interval with either the $\tauc$ or the $\tau_r$ topology but fails to be cell-wise euclidean. In Appendix \ref{section:examples}, Example~\ref{example_space_cell-wise_euclidean_not_metrizable}, we describe, as already noted, a $T_3$ definable topological space of dimension two that is cell-wise euclidean but not affine. 
Hence, although equivalent for one-dimensional definable topological spaces, the following three implications are strict in general. 
\begin{center}
\begin{tabular}{ c c c c c }
 Affine & $\Rightarrow$  & Hausdorff and       & $\Rightarrow$ & Hausdorff and does not contain \\ 
        &                & cell-wise euclidean &               & a definable copy of an interval \\  
        &                &                     &               & with either the $\tau_r$ or $\tauc$ topology    
\end{tabular}
\end{center}
In fact it is not even the case that being Hausdorff and cell-wise euclidean implies being definably metrizable, since, by o-minimality, a cell-wise euclidean space cannot contain an infinite definable discrete subspace, and so, by~\cite[Theorem 7.1]{walsberg15}, every cell-wise euclidean definably metrizable space is affine.

This complicates the task of generalising Theorem~\ref{them_main_2} to spaces of all dimensions. The next corollary however offers a possibility.  

\begin{corollary}\label{cor_them_2}
Suppose that $\RR$ expands an ordered field. Let $(X,\tau)$, \mbox{$\dim X\leq 1$}, be a definably compact Hausdorff definable topological space. The following are equivalent. 
\begin{enumerate}[(1)]
\item $(X,\tau)$ satisfies the frontier dimension inequality. \label{itm:cor_them_2_1}
\item $(X,\tau)$ is definably metrizable. \label{itm:cor_them_2_2}
\item $(X,\tau)$ is affine. \label{itm:cor_them_2_3}
\end{enumerate} 
\end{corollary}
\begin{proof}
We fix $(X,\tau)$, $\dim X\leq 1$, a definably compact Hausdorff definable topological space. 

$(\ref{itm:cor_them_2_3})\Rightarrow (\ref{itm:cor_them_2_2})$ is trivial. If $(X,\tau)$ is definably metrizable, then, by~\cite[Lemma 7.15]{walsberg15}, it satisfies the \textbf{fdi}, i.e. $(\ref{itm:cor_them_2_2})\Rightarrow (\ref{itm:cor_them_2_1})$. We complete the proof by showing $(\ref{itm:cor_them_2_1})\Rightarrow (\ref{itm:cor_them_2_3})$, that is, if $(X,\tau)$ satisfies the \textbf{fdi}, then it is affine. 

We prove the contrapositive. Suppose that $(X, \tau)$ is not affine. Then, by Theorem~\ref{them_main_2}, there exists an interval with either the $\tau_r$ or the $\tauc$ topology that definably embeds into $(X, \tau)$. We prove that $(X, \tau)$ does not have the \textbf{fdi}.

By Remark~\ref{remark_assumption_X}, we may assume that $X\subseteq R$. By Lemma~\ref{remark_classification_spaces_line}, there exists an interval $I\subseteq X$ such that $\tau|_I\in\{\tau_r, \tau_l, \tauc\}$. Considering the push-forward of $(X,\tau)$ by $x\mapsto -x$ if necessary, we may assume that $\tau|_I\in\{\tau_r, \tauc\}$. By definable compactness and Hausdorffness, for every $y\in I$, there exists a unique $x \in X$ such that $y\in L_x$ (see Proposition~\ref{prop_basic_facts_P_x_2}(\ref{itm3:basic_facts_P_x_2}) and Lemma~\ref{lem:RL-compact}). Since $\tau|_I\in\{\tau_r, \tauc\}$, we have that this $x$ must not belong in $I$, and in particular $x\in \partial_\tau I$. By Lemma~\ref{lemma_2}, it follows that $\partial_\tau I$ is infinite, and so $(X,\tau)$ does not satsify the \textbf{fdi}.    
\end{proof}

Using our work in Sections~\ref{section: universal spaces} and~\ref{section:compactifications}, we present a second refinement of Theorem~\ref{them_main_2} for definably compact spaces. We do not know whether or not this result generalizes to spaces of any dimension.

\begin{corollary} \label{thm:compact_affine}
Suppose that $\RR$ expands an ordered field. Let $(X,\tau)$, $\dim X\leq 1$, be a definably compact Hausdorff definable topological space. Exactly one of the following holds. 
\begin{enumerate}[(1)]
    \item \label{itm1:compact-affine} There exists an interval $I\subseteq R$, some $n>0$, and a definable open embedding $(I\times\{0,\ldots, n\}, \mu)\hookrightarrow (X,\tau)$, where $\mu\in \{\taulex, \taualex\}$.
    \item \label{itm2:compact-affine} $(X,\tau)$ is affine. 
\end{enumerate}
Furthermore, if $(X,\tau)$ is definably separable, then we can take $n$ to be $1$ and $\mu$ to be $\taulex$ in (1). 
\end{corollary}
\begin{proof}
Applying Remark~\ref{remark_assumption_X}, we may assume that $X\subseteq R$. On the one hand, if $(X,\tau)$ is affine, then in particular it is cell-wise euclidean (see comments after Definition~\ref{dfn:cell-wise_euclidean}). On the other hand, note that, for any interval $I$ and $n>0$, the space $(I\times\{0,\ldots, n\}, \mu)$, where $\mu$ is $\taulex$ (respectively $\taualex$), satisfies that the subspace $I\times\{n\}$ is homeomorphic, taking the projection to the first coordinate, to the interval $I$ with the $\tau_r$ (respectively $\tauc$) topology. Hence, if (\ref{itm1:compact-affine}) in the corollary holds, then $(X,\tau)$ contains a subspace definably homeomorphic to an interval with the $\tau_r$ or $\tauc$ topology, and so, by Remark~\ref{remark_no_homeomorphism}, (\ref{itm1:compact-affine}) and (\ref{itm2:compact-affine}) are mutually exclusive. 

By Lemma~\ref{lemma_Hausdorff_compact_is_regular}, the space $(X,\tau)$ is regular. Let $\XX=\Xinf \cup \Xsing$ be a partition of $X$ as defined by Lemma~\ref{remark_for_two_theorems} and, for every $A\in \Xinf$, let $A^*$, $\tau_A$ and $h_A$ be as given by Lemma~\ref{lemma_proof_two_thems}.
Since $(X,\tau)$ is definably compact, by Lemma~\ref{prop_proof_two_thems} we have that, for every $A\in \Xinf$, $h_A:(A,\tau)\rightarrow (A^*,\tau_A)$ is a homeomorphism. Furthermore, since every set in $\Xinf$ is $\tau$-open, we have that, for every $A\in\Xinf$, the map $h_A^{-1}:(A^*,\tau_A)\hookrightarrow (X,\tau)$ is an open embedding. Now recall from Lemma~\ref{lemma_proof_two_thems} that, for every $A\in \Xinf$, the set $A^*$ is of the form $I\times\{0,\ldots, n\}$, for some interval $I\subseteq R$ and some $n=n_A\geq 0$, and moreover $\tau_A$ is one of $\taulex$ or $\taualex$. On the one hand, if $n_A>0$ for some $A\in \Xinf$, then we clearly have condition (\ref{itm1:compact-affine}) in the corollary. On the other hand, if $n_A=0$ for every $A\in \Xinf$ then, using the fact that, for any interval $I$ it holds that $(I\times\{0\}, \taulex)=(I\times\{0\}, \taualex)=(I\times\{0\}, \tau_e)$, we derive that $(X,\tau)$ is cell-wise euclidean (see comments after Definition~\ref{dfn:cell-wise_euclidean}), and so, by the proof of Theorem~\ref{them_main_2}, it is affine.   

It remains to show that, if $(X,\tau)$ is definably separable, then we may take $n$ to be $1$ and $\mu$ to be $\taulex$ in (\ref{itm1:compact-affine}). By the above paragraph, it suffices to show that, if $(X,\tau)$ is definably separable, then $n_A\leq 1$ and $\tau_A=\taulex$ for every $A\in\Xinf$. This however is shown in the proof of Corollary~\ref{cor_T3_separable_embedding}.
\end{proof}

\subsection{Fremlin's conjecture}\label{subsection: Fremlin}

As indicated in Subsection~\ref{subsection: 3-el_basis_conj}, the 3-element basis conjecture (which asserts that statement~\eqref{star} is consistent with ZFC) arose as a result of interest in various questions about the nature of perfectly normal, compact, Hausdorff spaces. (Recall that a topological space is perfectly normal if and only if it is normal and any open set is a countable union of closed sets.) In particular, this conjecture is closely related to a conjecture of Fremlin, which posits that the following statement is consistent with ZFC (see Question 1 in \cite{gm07}): 
\begin{equation} \label{Fremlin}
    \parbox{0.9\textwidth}{Every perfectly normal, compact, Hausdorff space admits a continuous, at most $2$-to-$1$ map onto a compact metric space.}
\end{equation}
The underlying question is whether or not the Split Interval, which in our setting is $([0,1] \times \{0, 1\}, \taulex)$ defined in $(\mathbb{R},<)$ (see Example~\ref{example:split_interval}), is in essence the only example in ZFC of a non-metrizable, perfectly normal, compact, Hausdorff space, and how to make that notion precise. Fremlin's initial conjecture (see~\cite{gm07}) had been that every perfectly normal, compact, Hausdorff space in ZFC is a continuous image of $([0,1] \times \{0, 1\}, \taulex) \times ([0,1],\tau_e)$ (where likewise, in our setting, this space should be understood as being defined in $(\mathbb{R},<)$), but this does not hold (see~\cite{watsonweiss88}). However, $([0,1] \times \{0, 1\}, \taulex)$ certainly admits a continuous, at most 2-to-1 map onto a compact metric space, as does the counterexample of~\cite{watsonweiss88}. 

In considering Fremlin's conjecture, Gruenhage put forward the following question that is in the same spirit in~\cite{gru88} (see also Question 2.2 in~\cite{gru90}). He asked if the following statement is consistent with ZFC:
\begin{equation} \label{Fremlin2.2}
    \parbox{0.9\textwidth}{Every non-metrizable, perfectly normal, compact, Hausdorff space contains a copy of $(A \times \{0,1\}, \taulex)$, for some uncountable $A \subseteq [0,1]$.}
\end{equation}
It is indicated in~\cite{gru90} and~\cite{gm07} that, under PFA, statement~\eqref{Fremlin} is equivalent to the existence of a 3-element basis, as posited in~\eqref{star}, for the class of subspaces of perfectly normal, compact, Hausdorff spaces, and furthermore that, under PFA, both of these statements imply~\eqref{Fremlin2.2}. It is also indicated that, even without assuming PFA, statement~\eqref{star} implies both~\eqref{Fremlin} and~\eqref{Fremlin2.2}. 

Here, we consider the questions above and show, as a corollary to Theorems~\ref{them_ADC_or_DOTS} and~\ref{them_main_2}, that statement~\eqref{Fremlin} holds in a definable sense for any regular, Hausdorff one-dimensional topological space which is either perfectly normal or separable and which is definable in any o-minimal expansion of $(\mathbb{R}, +, \cdot, <)$ (Corollary~\ref{cor_Fremlin_2-2}). Moreover, we address statement~\eqref{Fremlin2.2} by showing that any  one-dimensional, perfectly normal, compact, Hausdorff space definable in any o-minimal expansion of $(\mathbb{R}, +, \cdot, <)$ is either affine or there exists an interval $I\subseteq \mathbb{R}$ and a definable open embedding $(I\times\{0,1\}, \taulex) \hookrightarrow (X,\tau)$ (Corollary~\ref{cor:Fremlin-Q2.2}). We also address definable generalizations of some of these statements for o-minimal expansions of ordered fields.

We begin by proving a definable result in our setting that is closely related to statement~\eqref{Fremlin}. We then specifically address statement~\eqref{Fremlin} in our setting, whenever $\RR$ expands the field of reals, through a subsequent lemma and corollary.

\begin{corollary} \label{cor_Fremlin_2-1}
Suppose that $\RR$ expands an ordered field. Let $(X,\tau)$, $\dim \leq 1$, be a regular and Hausdorff definable topological space. If $(X,\tau)$ is definably separable, then there exists a definable continuous map \mbox{$f:(X,\tau)\rightarrow (R^n,\tau_e)$}, for some $n$, where $f$ is at most $2$-to-$1$ (i.e. $|f^{-1}(z)|\leq 2$ for every $z\in R^n$). 
\end{corollary}
\begin{proof}
Applying Remark~\ref{remark_assumption_X} we may assume that $X\subseteq R$. 
We will use the construction in the proof of Theorem~\ref{them_compactification} of a (definably near-compact) space $(X^*,\tau^*)$ and of a definable embedding $h:(X,\tau) \hookrightarrow (X^*,\tau^*)$. We show that, if $(X,\tau)$ is definably separable, then $(X^*,\tau^*)$ admits a definable at most $2$-to-$1$ continuous map into a Hausdorff cell-wise euclidean space $(Y,\mu)$, with $Y\subseteq R$. Since $(Y,\mu)$ is Hausdorff and cell-wise euclidean it follows from Remark~\ref{remark_no_homeomorphism} and Theorem~\ref{them_main_2} that it is affine, completing the proof.   

Let $\XX=\Xinf \cup \Xsing$ be a partition of $X$ as defined by Lemma~\ref{remark_for_two_theorems} and, for every $A\in \Xinf$, let $A^*$ and $\tau_A$ be as given by Lemma~\ref{lemma_proof_two_thems}. Set $\Asing:=\bigcup \Xsing$, and let $(X^*, \tau^*)$ be as given by the proof of Theorem~\ref{them_compactification}. Recall that $X^*=\bigcup\{A^* : A\in \Xinf\} \cup \{\al x, 0\ar : x\in\Asing\}$. Now, in the proof of Corollary~\ref{cor_T3_separable_embedding} it is shown that, if $(X,\tau)$ is definably separable, then, for every $A\in \Xinf$, it holds that $A^*\subseteq R\times\{0,1\}$ and $(A^*,\tau_A)=(A^*,\taulex)$. In particular, we have that $X^*\subseteq R\times\{0,1\}$. 
Consider the projection $\pi:X^*\rightarrow R$ to the first coordinate. Since $X^*\subseteq R\times\{0,1\}$, then $\pi$ is at most $2$-to-$1$. Let $Y=\pi(X^*)$. We define a Hausdorff cell-wise euclidean topology $\mu$ on $Y$ such that $\pi:(X^*,\tau^*)\rightarrow (Y,\mu)$ is continuous. 

For every $A\in \Xinf$, recall that the set $A^*$ is of the form $I\times\{0,1\}$ for some interval $I\subseteq A$. We denote this interval by $(a_A,b_A)$. In particular, for every $A\in \Xinf$, we have that $\pi(A^*)=(a_A,b_A)$, and $\Xsing \cup \{ (a_A, b_A) : A\in \Xinf\}$ is a partition of $Y$. 
We define $\mu$ as follows. For every $A\in \Xinf$, the interval $(a_A, b_A)$ is $\mu$-open and its subspace topology is euclidean. Furthermore, for every $x\in \Asing$, we say that a set $V\subseteq Y$ is a $\mu$-neighbourhood of $x$ if and only if $\pi^{-1}(V)$ is a $\tau^*$-neighbourhood of $\al x, 0\ar$.

We give a precise description of a basis of neighbourhoods in $(Y,\mu)$ for each point in $\Asing$. Let $x\in \Asing$ and let $\BB(x)$ denote the definable basis of neighbourhoods of the point $\al x, 0 \ar$ in $(X^*,\tau^*)$ described in the proof of Theorem~\ref{them_compactification}. Then a definable basis of neighbourhoods of $x$ in $(Y,\mu)$ is given by the family $\UU(x)=$ $\{\pi(V) : V\in \BB(x)\}$. We can describe $\UU(x)$ more precisely using the definition of $\BB(x)$ as follows. 
Let $\Xinf^{R_x} = \{ A \in \Xinf : a_A \in R_x\}$ and $\Xinf^{L_x} = \{ A \in \Xinf : b_A \in L_x\}$. Then $\UU(x)$ consists of every set of the form 
\[
\{x\} \cup \bigcup_{A \in \Xinf^{R_x}} (a_A,y_A) \cup \bigcup_{A \in \Xinf^{L_x}} (z_A, b_A),
\]
for parameters $\{y_A : A \in \Xinf^{R_x}\}$ and $\{ z_A : A \in \Xinf^{L_x}\}$ satisfying that, for every $A \in \Xinf^{R_x}$, we have $y_A \in (a_A, b_A)$ and, for every $A \in \Xinf^{L_x}$, we have $z_A \in (a_A, b_A)$.
It is thus easy to check that the topology $\mu$ is well defined and, by Proposition~\ref{prop_basic_facts_P_x_2}(\ref{itm3:basic_facts_P_x_2}), Hausdorff. Furthermore it is cell-wise euclidean by definition. 

It remains to show that $\pi:(X^*,\tau^*)\rightarrow (Y,\mu)$ is continuous; in other words, that, for every point $\al x, i\ar \in X^*$ and every $\mu$-neighbourhood $V$ of $\pi(\al x,i \ar)=x$, the set $\pi^{-1}(V)$ is a $\tau^*$-neighbourhood of $\al x, i \ar$. 
We observe that this follows easily from the definitions of $(X^*,\tau^*)$ and $(Y,\mu)$. Specifically, if $x\in \Asing$ then, by definition of $X^*$, we have that $i=0$, and the observation is explicit in the definition of $\mu$. For the remaining case, suppose that $x\in (a_A,b_A)$ for some $A\in \Xinf$. By definition of $\taulex$ note that every interval $J\subseteq (a_A,b_A)$ satisfies that $\pi^{-1}(J)$ is open in $(A^*,\taulex)$. Hence the observation follows from the facts that $((a_A,b_A),\mu)=((a_A,b_A),\tau_e)$, and $A^*$ is $\tau^*$-open with $(A^*,\tau^*)=(A^*,\taulex)$.
\end{proof}

Through the following lemma and corollary we prove a stronger form of statement~\eqref{Fremlin} in our setting, where we do not require the assumption that the space be compact, and furthermore we may replace perfect normality by separability. Since the proof of Corollary~\ref{cor_Fremlin_2-2} is similar to the proof of Corollary~\ref{cor_Fremlin_2-1} above, we only sketch its proof. 

\begin{lemma}\label{lem:T5}
Let $I\subseteq \mathbb{R}$ be an interval and $0 \leq n \leq m$. The subspace $I\times\{0,\ldots,n\}$ of the space $(I\times\{0,\ldots,m\}, \taulex)$ is perfectly normal if and only if either $n=0$ or $n=m=1$. 
The space $(I\times\{0,\ldots, n\}, \taualex)$ is perfectly normal if and only if $n=0$. 
\end{lemma}
\begin{proof}
The fact that the euclidean space $(I\times\{0\},\taulex)=(I\times\{0\},\taualex)=(I\times\{0\},\tau_e)$ and the space $(I\times\{0,1\},\taulex)$ are both perfectly normal is classical and so we omit the proofs. Similarly, it is also classical that the space $(I,\tau_l)$ is perfectly normal, and so, for any $m>0$, the subspace $I\times\{0\}$ of the space $(I\times\{0,\ldots, m\},\taulex)$ is perfectly normal, since this subspace corresponds to the push-forward of $(I,\tau_l)$ by the map $x\mapsto\al x, 0 \ar$.

Let $(I\times\{0,\ldots,n\},\mu)$ be a space in any of the remaining cases, i.e. such that $n>0$, and either with $\mu=\taualex$ or otherwise satisfying that there is $m>n$ such that $\mu$ is the subspace topology inherited from the space $(I\times\{0,\ldots,m\},\taulex)$. Note that, in all of these cases, there exists some $0<i \leq n$ such that the subset $(I\times \{i\}, \mu)$ only contains $\mu$-isolated points, and in particular it is $\mu$-open. We fix such an $i$. 

Recall that a topological space is perfectly normal if and only if it is normal and any open set is a countable union of closed sets. Hence to show that $(I\times\{0,\ldots,n\},\mu)$ is not perfectly normal it suffices to show that $I\times\{i\}$ is not a countable union of $\mu$-closed sets. Towards a contradiction, suppose that it were, in which case there would exist an uncountable $\mu$-closed subset $C$ of $I\times\{i\}$.
Since such a $C$ is uncountable, the projection $\pi(C)$ of $C$ to the first coordinate would satisfy that there exists a point $x\in \pi(C)$ such that, for any $a<x<b$, it holds that $(a,x)\cap \pi(C)\neq \emptyset$ and $(x,b)\cap \pi(C) \neq \emptyset$ (see Lemma~\ref{lemma_limit_points} for a generalization of this classical fact of the reals). By definition of the $\taulex$ and $\taualex$ topologies, it follows that $\al x, 0\ar$ is in the $\mu$-closure of $C$, a contradiction, since $C$ is a $\mu$-closed subset of $I\times \{i\}$ with $i>0$. 
\end{proof}

\begin{corollary} \label{cor_Fremlin_2-2}
Suppose that $\RR$ expands the ordered field of reals \mbox{$(\mathbb{R},+,\cdot,<)$.} Let $(X,\tau)$, $\dim X \leq 1$, be a regular and Hausdorff definable topological space. If $(X,\tau)$ is either separable or perfectly normal, then there exists a definable continuous map \mbox{$f:(X,\tau)\rightarrow (R^n,\tau_e)$}, for some $n$, where $f$ is at most $2$-to-$1$.
\end{corollary}
\begin{proof}
If $(X,\tau)$ is separable then, by Proposition~\ref{prop:sep-equiv}, it is also definably separable, and the result follows from Corollary~\ref{cor_Fremlin_2-1}. Suppose that $(X,\tau)$ is perfectly normal. 

Applying Remark~\ref{remark_assumption_X} we may assume that $X\subseteq R$. Let $\XX=\Xinf \cup \Xsing$ be a partition of $X$ as defined by Lemma~\ref{remark_for_two_theorems} and, for every $A\in \Xinf$, let $h_A$ be as given by Lemma~\ref{lemma_proof_two_thems}. Since perfect normality is a hereditary property we have that, for every set $A\in \Xinf$, the subspace $(A,\tau)$ is perfectly normal.
Applying Lemma~\ref{lem:T5} to the construction by cases in the proof of Lemma~\ref{lemma_proof_two_thems}, we observe that, for every $A\in\Xinf$, it holds that $h_A(A)\subseteq R\times\{0,1\}$. The rest of the proof is now analogous to the proof of Corollary~\ref{cor_Fremlin_2-1}, taking $(h(X),\tau^*)$ in place of $(X^*,\tau^*)$ in said proof. 
\end{proof}

Observe that the proofs of Corollaries~\ref{cor_Fremlin_2-1} and~\ref{cor_Fremlin_2-2} can be adapted to yield a different version of Corollary~\ref{cor_Fremlin_2-1}, which holds for any expansion $\RR$ of an ordered field, where the condition that $(X, \tau)$ is definably separable is replaced by the weaker condition that $(X, \tau)$ does not have a subspace definably homeomorphic to any space that is the $\RR$-definable analogue of one of the spaces that Lemma~\ref{lem:T5} shows are not perfectly normal.

We now consider statement~\eqref{Fremlin2.2} in our setting. We present a definable version of this statement, in which the metrizability condition is replaced by the stronger property of being affine. It is an immediate consequence of Lemma~\ref{lem:T5} applied to the dichotomy described in Corollary~\ref{thm:compact_affine}.

\begin{corollary}\label{cor:Fremlin-Q2.2}
Suppose that $\RR$ expands the ordered field of reals \mbox{$(\mathbb{R},+,\cdot,<)$.} Let $(X,\tau)$, $\dim X\leq 1$, be a definably compact (equivalently, by Remark~\ref{rem:compactness}, compact) and Hausdorff definable topological space. If $(X,\tau)$ is perfectly normal then it is either affine or otherwise there exists an interval $I\subseteq \mathbb{R}$ and a definable open embedding $(I\times\{0,1\}, \taulex) \hookrightarrow (X,\tau)$.
\end{corollary}

One may in fact replace the assumption of perfect normality in Corollary~\ref{cor:Fremlin-Q2.2} by the weaker condition that there does not exist an interval $I\subseteq \mathbb{R}$ and a definable embedding $(I\times\{0,1\},\taualex)\hookrightarrow (X,\tau)$. By Corollary~\ref{thm:compact_affine} alone, this condition already implies that $(X,\tau)$ is affine or otherwise there exists an interval $I\subseteq \mathbb{R}$ and a definable embedding $(I\times\{0,1\}, \taulex) \hookrightarrow (X,\tau)$. (To see this observe that, for any interval $I\subseteq \mathbb{R}$ and $n>0$, the subspace $I\times \{0, n\}$ of the space $(I\times \{0,\ldots, n\}, \taulex)$ is definably homeomorphic to $(I\times \{0,1\}, \taulex)$.) Additionally, since Corollary~\ref{thm:compact_affine} holds in the setting of an arbitrary o-minimal expansion $\RR$ of an ordered field, this version of Corollary~\ref{cor:Fremlin-Q2.2} moreover holds in this more general setting. 

\section{Definable metrizability}\label{section:metrizability}

In this section, we use our affineness characterization of the previous section (Theorem~\ref{them_main_2}) to derive the equivalence of metrizability and definable metrizability for one-dimensional definable topological spaces in certain o-minimal expansions of ordered fields. 

Throughout, we continue to assume that $\RR$ expands an ordered field. Recall that in our setting \textquotedblleft metric" refers to an $\RR$-metric (Definition~\ref{dfn:R-metric}), including those instances when it appears implicitly in notions such as metrizability and metric space. 

We will frequently utilize two classical topological notions in stating and proving results throughout this section. The first of these is the weight $w_\tau(X)$ of a topological space $(X,\tau)$, i.e. the minimum cardinality of a basis for $\tau$ (this was discussed earlier in Subsection \ref{subsection: Cantor_space}). Note that it follows from Theorem~\ref{them_main_2}, Lemma~\ref{lem:weight-taur} and Proposition~\ref{prop:weight}(\ref{itm_weight_2}) that, whenever $\RR$ expands an ordered field, every one-dimensional Hausdorff definable topological space $(X,\tau)$ satisfies that $w_\tau(X)\in\{w_e(R), |R|\}$. The other classical notion that we will use is the density of $(X,\tau)$, i.e. the minimum cardinality of a $\tau$-dense subset of $X$, which we denote by $den_\tau(X)$. 

Our main result, Theorem~\ref{thm:metrizability}, shows that, whenever $\RR$ is an o-minimal expansion of an ordered field satisfying that $den_e(R)<|R|$ (e.g. whenever $\RR$ expands the field of reals), every metrizable one-dimensional space is in fact definably metrizable.

We begin by deriving the following as a straightforward corollary of Theorem~\ref{them_main_2}.
\begin{corollary}\label{thm:metrizability_2}
Let $(X,\tau)$, $\dim(X)\leq 1$, be a definable topological space that is metrizable and separable. Suppose that either of the following two conditions holds. 
\begin{enumerate}[(a)]
\item $\RR$ expands the field of reals.
\item $(X,\tau)$ is compact.
\end{enumerate}
Then $(X,\tau)$ is affine. In particular it is definably metrizable. 
\end{corollary}
\begin{proof}
The case where $X$ is finite is trivial, so we assume that $\dim(X)=1$. Recall from Remark~\ref{remark_compact_implies_reals} that, if there exists a compact infinite $T_1$ topological space definable in $\RR$, then $\RR$ expands the field of reals. Since a metrizable space is $T_1$, the second case here therefore reduces to the first, and we may consequently assume that $\RR$ expands the field of reals, i.e. $R=\mathbb{R}$.
 
By Remark~\ref{remark_assumption_X}, we assume that $X\subseteq \mathbb{R}$. By Theorem~\ref{them_main_2}, it is enough to show that $(X,\tau)$ does not have a definable copy of an interval with either the discrete or the right half-open interval topology. This follows from the fact that $(X,\tau)$ is a separable metric space, hence second countable, so $w_{\tau}(X)<|\mathbb{R}|$, while the topological weight of an interval with the Sorgenfrey Line or discrete topology is $|\mathbb{R}|$ (see Lemma~\ref{lem:weight-taur}). 
\end{proof}  

We now state the main theorem of this section, which improves the metrization part of Corollary~\ref{thm:metrizability_2}. 

\begin{theorem}\label{thm:metrizability}
Suppose that $\RR$ expands an ordered field and satisfies that \mbox{$den_e(R)<|R|$}. Let $(X,\tau)$, $\dim X\leq 1$, be a definable topological space. Then $(X,\tau)$ is metrizable if and only if it is definably metrizable. 
\end{theorem}

In the case where $den_e(R)=|R|$, we point the reader towards~\cite{diep17} for a proof that the space $(R,\tau_r)$, which is shown in Proposition~\ref{prop:tau_r not metrizable} not to be definably metrizable, is metrizable whenever $\RR$ expands a countable densely ordered group. Classically, the Sorgenfrey Line $(\mathbb{R},\tau_r)$ is separable but not second countable, and so it is not ($\mathbb{R}$-)metrizable. 
We show below (Remark~\ref{remark_not_metrizable}) that the analogous statement holds whenever $den_e(R) < |R|$, i.e. in this case, $(X, \tau_r)$ is not ($\RR$)-metrizable when $X \subseteq R$ is any infinite set.

In order to prove Theorem~\ref{thm:metrizability}, we require two simple lemmas, whose aim is to generalise basic results in metric topology and the topology of the real line to our setting. In what follows, recall that, since $\RR$ expands an ordered field, any two intervals are definably $e$-homeomorphic and in particular, for any interval $I\subseteq R$, we have that $|I|=|R|$ and $den_e(I)=den_e(R)$.

\begin{lemma}\label{lemma_tech_2nd_countable}
Let $(X,d)$ be a metric space. Let $\tau:=\tau_d$, and let $A,C\subseteq X$ satisfy that $A\subseteq cl_\tau C$ (i.e. $C$ is $\tau$-dense in $A$). Let $D$ be an $e$-dense subset of $(0,+\infty)$. Consider the family of $d$-balls $\BB=\{B_d(y,\delta) : y\in C, \delta\in D\}$. Then, for every $x\in A$, there exists a subfamily $\BB_x$ of $\BB$ that is a basis of open $\tau$-neighbourhoods of $x$.

In particular $w_\tau(X)\leq den_\tau(X)den_e(R)$.   
\end{lemma}
\begin{proof}
Let $x\in A$ and $\varepsilon>0$. We must show that there exists $y\in C$ and $\delta\in D$ such that $x\in B_d(y, \delta)\subseteq B_d(x,\varepsilon)$. 

Let $\delta\in D$ be such that $0<\delta< \varepsilon/2$. Since $A\subseteq cl_\tau C$, there exists $y\in C$ such that $d(x,y)<\delta$. Consider the ball $B_d(y,\delta)$. Clearly $x\in B_d(y,\delta)$ and, if $z\in B_d(y,\delta)$, then, by the triangle inequality, $d(x,z) \leq d(x,y) + d(y,z) \leq \delta+\delta < \varepsilon$. Hence $x\in B_d(y,\delta) \subseteq B_d(x,\varepsilon)$, which completes the proof of the first part of the lemma.

For the second part, suppose that $C$ is a dense subset of $(X,\tau)$ of cardinality $den_\tau(X)$ and $D$ is an $e$-dense subset of $(0,+\infty)$  of cardinality $den_e(R)$. Then, by the above, the family $\{ B_d(y,\delta) : y\in C, \delta\in D\}$ is a basis for $\tau$ which has cardinality bounded by $den_\tau(X)den_e(R)$.
\end{proof}

\begin{remark}\label{remark_not_metrizable}
Let $X\subseteq R$ be an infinite definable set. Recall that, by Lemma~\ref{lem:weight-taur}, because $\RR$ expands an ordered field, the space $(X,\tau_*)$, where $\tau_*\in\{\tau_l, \tau_r\}$, has weight $|R|$. Moreover, clearly the density of $(X,\tau_*)$ is equal to $den_e(R)$. 
From Lemma~\ref{lemma_tech_2nd_countable}, it follows that, if $(X,\tau_*)$ is metrizable, then $|R|=w_{\tau_*}(X)\leq den_{\tau_*}(X) den_e(R)= den_e(R)^2=den_e(R)\leq |R|$, i.e. $den_e(R)=|R|$. It follows that, if in fact we have $den_e(R)<|R|$, then $(X,\tau_*)$ is not metrizable.   
\end{remark}

For the next lemma recall that, for any set $X\subseteq R$, a right (respectively left) limit point of $X$ is a point $x\in R$ satisfying that, for every $y>x$ (respectively $y<x$), $(x,y)\cap X\neq \emptyset$ (respectively $(y,x)\cap X\neq \emptyset$).

\begin{lemma}\label{lemma_limit_points}
Suppose that $den_e(R)<|R|$ and let $X\subseteq R$ be a subset of cardinality $|R|$. Then there exist $|R|$-many elements of $X$ that are both right and left limit points of $X$.
\end{lemma}
\begin{proof}
We show that all but at most $den_e(R)$-many points in $X$ are right limit points of $X$. The same holds for left limit points. The result then follows from the fact that $den_e(R)<|R|=|X|$. 

Let $Y$ be the set of points in $X$ that are not right limit points of $X$. For every $x\in Y$, there is some $x'>x$ such that the interval $(x,x')=I_x$ is disjoint from $X$. The family $\{I_x : x\in Y\}$ has cardinality $|Y|$ and contains only non-empty pairwise disjoint intervals. It follows that $|Y|\leq den_e(R)$. The proof for the set of left limit points is analogous.   
\end{proof}

We may now prove Theorem~\ref{thm:metrizability}.

\begin{proof}[Proof of Theorem~\ref{thm:metrizability}]
Clearly any definably metrizable topological space is metrizable.  Fix $(X,\tau)$, with $\dim X \leq 1$, a definable topological space whose topology is induced by a metric $d$. We prove that $(X,\tau)$ is definably metrizable by describing a definable metric $\hat{d}$ that induces $\tau$. Since every finite metric space is discrete, we may assume that $\dim X=1$. Let $D$ be a dense subset of $R$ of cardinality $den_e(R)$. By Remark~\ref{remark_assumption_X}, we assume that $X$ is a bounded subset of $R$. 

Consider the definable set $S=\{x\in X : \E_x \setminus \{x\} \neq \emptyset\}$. We begin by proving the following claim. 
\begin{claim} \label{claim:S-finite}
$S$ is finite.
\end{claim}
\begin{proof}[Proof of claim]
Towards a contradiction suppose that $S$ is infinite. Let $f:S\rightarrow \exR$ be the map given by $x\mapsto \min \E_x\setminus \{x\}$, which, by Proposition~\ref{basic_facts_P_x}(\ref{itm: basic_facts_2}) and Lemma~\ref{lemma_2}, is definable. By Hausdorffness (Lemma~\ref{lemma:basic_facts_P_x_2}(\ref{itmb:lemma_basic_facts_P_x_2})) and o-minimality, there exists an interval $I\subseteq S$ on which $f$ is \econtinuous and strictly monotonic. Note (see Lemma~\ref{lemma_explaining_neighbourhoods}) that $I$ is in the $\tau$-closure of $D \cap f(I)$. Consider the family of $d$-balls $\BB=\{B_d(q,\delta) : q \in D \cap f(I), \delta \in D \cap (0,+\infty) \}$. This family has cardinality bounded by $den_e(R)$ and, by Lemma~\ref{lemma_tech_2nd_countable}, contains, for every $x\in I$, a subfamily that is a basis of open $\tau$-neighbourhoods of $x$.

Now, let $h : I \to \BB$ be a function with the property that, for every $x \in I$, $h(x) \in \BB$ is a $\tau$-neighbourhood of $x$ such that $f(x) \notin h(x)$. Such a function can be defined since, for every $x \in I$, $f(x) \neq x$ and $\tau$ is $T_1$. 
Since $|I|=|R|$ and $|\BB|\leq den_e(R)<|R|$, there must exist, by the pigeonhole principle, some $d$-ball $B\in h(I)$ such that the set $h^{-1}(B)$ has cardinality $|R|$. By Lemma~\ref{lemma_limit_points}, there exists $x\in h^{-1}(B)$ that is both a right and left limit point of $h^{-1}(B)$. Recall that $f(x)\in\E_x$. Suppose that $f(x)\in R_x$. Then, since $B=h(x)$ is a $\tau$-neighbourhood of $x$, there is some $z>f(x)$ such that $(f(x),z)\subseteq B$. If $f$ is increasing then, by $e$-continuity, there is some $y>x$ with $(x,y)\subseteq I$ such that $f[(x,y)]\subseteq (f(x),z)$. 
Hence, for every $x'\in (x,y)$, it holds that $f(x') \in B$ and so, by definition of $h$, that $h(x')\neq B$. However, this contradicts that $x$ is a right limit point of $h^{-1}(B)$. Similarly, if $f$ is decreasing, there is some $y<x$ with $(y,x)\subseteq I$ such that $(y,x)\cap h^{-1}(B)=\emptyset$, contradicting that $x$ is a left limit point of $h^{-1}(B)$. The argument in the case where $f(x)\in L_x$ is analogous. This completes the proof of the claim. 
\renewcommand{\qedsymbol}{$\square$ (claim)}
\end{proof}

We now continue the proof of Theorem~\ref{thm:metrizability}. By Theorem~\ref{them 1.5} and Remark~\ref{remark_not_metrizable}, there exists a partition $\XX$ of $X$ into finitely many points and intervals where each interval subspace in $\XX$ has either the euclidean or the discrete topology. Let $E_S=\bigcup_{x\in S} \E_x$. By Claim~\ref{claim:S-finite} and Lemma~\ref{lemma_2}, both $S$ and $E_S$ are finite sets. By passing to a finer partition if necessary, we may require that $\XX$ has the following two properties. 
\begin{enumerate}[(i)]
\item \label{itm1: metrizability} The elements in $S$ and in $E_S$ do not belong in any interval in $\XX$. 
\item \label{itm2: metrizability} For any interval $(a,b)\in \XX$ with the discrete subspace topology, it holds that, if $a\in \bigcup_{x\in X} R_x$, then $b\notin \bigcup_{x\in X} L_x$ and, if $b\in \bigcup_{x\in X} L_x$, then $a\notin \bigcup_{x\in X} R_x$.  
\end{enumerate} 

Note that (\ref{itm2: metrizability}) can be arranged since any discrete interval subspace $I$ that is disjoint from $E_S$ is also disjoint from $\bigcup_{x\in X} \E_x$, and so any proper subinterval of $I$ has the desired property.

First note that, by (\ref{itm1: metrizability}), for any interval $I=(a,b)\in\XX$, any $x\in I$ and any $y\in X\setminus I$, it holds that $\E_x\subseteq \{x\}$ and $\E_y\cap I=\emptyset$. So, by Lemma~\ref{lemma_explaining_neighbourhoods}, $I$ is $\tau$-open and, if $y\in \partial_\tau I$, then it must be that either $a\in R_y$ or $b\in L_y$. In particular, by (\ref{itm2: metrizability}) and Hausdorffness (Proposition~\ref{prop_basic_facts_P_x_2}(\ref{itm3:basic_facts_P_x_2})), if $I$ is discrete then $|\partial_\tau I|\leq 1$. 

Now, let $\YY \subseteq \XX$ be the family of all discrete interval subspaces in $\XX$. Let $|\YY|=n$. We prove the theorem by induction on $n$. 

If $n=0$, then $X$ is cell-wise euclidean. In particular, by Remark~\ref{remark_no_homeomorphism}, it contains no definable copy of an interval with either the discrete or the right half-open interval topology and so, applying Theorem~\ref{them_main_2}, $(X,\tau)$ is affine, and in particular it is definably metrizable.

Suppose that $n>0$ and let $\YY=\{I_1,\ldots, I_n\}$. Let $X'=X\setminus I_n$. By induction hypothesis, the space $(X',\tau)$ is definably metrizable with some definable metric $d'$. We extend $d'$ to a definable metric $\hat{d}$ on $X$ such that $\tau_{\hat{d}}=\tau$. Let $I_n=I=(a,b)$. By the argument above, $\partial_{\tau}I = \emptyset$ or $|\partial_{\tau}I| = 1$. We consider each of these two cases in turn. 

\textbf{Case 0:} $\partial_\tau I =\emptyset$. 
In this case, $I$ is a $\tau$-clopen subset of $X$. 
Note that the metric $\min\{1,d'\}$ induces the same topology as $d'$, hence, by passing to the former if necessary, we may assume that $d'\leq 1$. We define the metric $\hat{d}$ on $X$ as follows. 
\begin{itemize}
\item For all $x,y \in X'$, $\hat{d}(x,y)=d'(x,y)$. 
\item For all $x\in I$, $y \in X$, $\hat{d}(x,y)=\hat{d}(y,x)=1$, if $x\neq y$, and $\hat{d}(x,y)=\hat{d}(y,x)=0$ otherwise. 
\end{itemize}
Since the $\tau$-topology on $I$ is discrete, it is easy to check that $\hat{d}$ is a metric that induces the topology $\tau$ on $X$. 
 
\textbf{Case 1:} $|\partial_\tau I|=1$, i.e. $\partial_\tau I=\{x_0\}$ for some $x_0\in X\setminus I$. 
Recall that, by~(\ref{itm1: metrizability}), $\E_{x_0}\cap (a,b)=\emptyset$, and so (see Lemma~\ref{lemma_explaining_neighbourhoods}) it must be that either  $a \in R_{x_0}$ or $b \in L_{x_0}$, but not both (by~(\ref{itm2: metrizability})). We prove the case where $a\in R_{x_0}$. The remaining case, where $b\in L_{x_0}$, is analogous. 

Recall that $X$ is bounded, and so $I$ is a bounded interval. Consider the following definable metric $\hat{d}$ in $X$.
\begin{itemize}
\item For all $x,y \in X'$, $\hat{d}(x,y)=d'(x,y)$. 
\item For all $x, y\in I$, $\hat{d}(x,y)=|x-a|+|y-a|$, if $x\neq y$, and $\hat{d}(x,y)=0$ otherwise. 
\item For all $x\in I$, $y\in X'$, $\hat{d}(x,y)=\hat{d}(y,x)=|x-a|+d'(y,x_0)$. 
\end{itemize}  
It is routine to check that $\hat{d}$ is a metric. We show that $\tau_{\hat{d}}=\tau$ by proving, using Proposition~\ref{prop_cont_lim}, that the identity map $(X,\tau)\rightarrow(X,\tau_{\hat{d}})$ is a homeomorphism. Note that, by definition, $\hat{d}$ induces the corresponding subspace topologies of $\tau$ on $X'$ and $I$, and moreover $I$ is $\hat{d}$-open. In particular, since $I$ is $\tau$-open too, we have that $X'$ is both $\tau$-closed and $\hat{d}$-closed in $X$. Since $\tau|_{X'}=\tau_{\hat{d}}|_{X'}$, we derive that any definable curve in $X'$ $\tau$-converges to a point in $X$ (necessarily inside $X'$) if and only if it $\hat{d}$ converges to that same point. 

Furthermore, note that, by definition of $\hat{d}$, an injective definable curve $\gamma$ in $I$ $\hat{d}$-converges if and only if it $e$-converges to $a$ from the right, $\hat{d}$-converging thus to $x_0$. By the fact that $a\in R_{x_0}$ and Remark~\ref{remark_side_convergence}, $\gamma$ must then also $\tau$-converge to $x_0$. 
Conversely, if an injective definable curve $\gamma$ in $I$ $\tau$-converges then, by the facts that the $\tau$-topology on $I$ is discrete and $\partial_\tau I=\{x_0\}$, $\gamma$ must $\tau$-converge to $x_0$. Recall that, by the assumptions on $I$, we have that $\E_{x_0}\cap (a,b)=\emptyset$, $a \in R_{x_0}$ and $b \notin L_{x_0}$. Hence, by
Remark~\ref{remark_side_convergence}, it must be that $\gamma$ $e$-converges to $a$ from the right, and so, from the definition of $\hat{d}$, it follows that $\gamma$ $\hat{d}$-converges to $x_0$ too. This completes the proof of the theorem.
\end{proof}

It remains open whether or not Theorem~\ref{thm:metrizability} can be generalized to spaces of dimension greater than one.

\begin{question}
Let $\RR$ be an o-minimal expansion of an ordered field satisfying that \mbox{$den_e(R)<|R|$}. Is any $\RR$-metrizable topological space definable in $\RR$ definably metrizable?
\end{question}

\section{A note on an affiness result by Peterzil and Rosel}\label{sec:Peterzil_Rosel}

The majority of the work in this paper already appeared in the first author's doctoral dissertation~\cite{andujar_thesis}. After that work had been completed, the authors learned that Peterzil and Rosel were working on similar questions. Their work resulted in~\cite{pet_rosel_18}. The main theorem (page $1$) in said paper is the following affiness result, which we present in the terminology of the present paper.

\begin{theorem}[\cite{pet_rosel_18}] \label{them:PR}
Suppose that $\RR$ expands an ordered group. Let $(X,\tau)$ be a Hausdorff definable topological space, where $\dim X=1$ and $X$ is a bounded set. The following are equivalent. 
\begin{enumerate}[(1)]
\item  \label{itm:Pet_Rosel_1} $(X,\tau)$ is affine. 

\item  \label{itm:Pet_Rosel_2} There is a finite set $G\subseteq X$ such that the subspace topology $\tau|_{X\setminus G}$ is coarser than the euclidean topology on $X\setminus G$. 

\item \label{itm:Pet_Rosel_3} Every definable subspace of $(X,\tau)$ has finitely many definably connected components. 

\item \label{itm:Pet_Rosel_4} $(X,\tau)$ is regular and has finitely many definably connected components. 
\end{enumerate}

Furthermore, if $\RR$ expands an ordered field then the above is true even without the assumption that the set $X$ is bounded.
\end{theorem}

Their work is in some ways parallel to ours. For example, their notion of set of shadows $S(x)$ of a point $x$ is effectively the $e$-accumulation set $\E_x$ of $x$ (both notions are equivalent for Hausdorff topologies, which are the only topologies studied in~\cite{pet_rosel_18}, while in general it holds that $\E_x \subseteq S(x)$). 
Similarly, \emph{$x$ inhabits the left (respectively right) side of $y$} means $y\in L_x$ (respectively $y\in R_x$). 

Moreover, for a definable topological space $(X,\tau)$, they refer to the property \emph{almost $\tau \subseteq \tau^{af}|_X$} (where $\tau^{af}$ is their notation for the euclidean topology $\tau_e$) to mean the condition that $(X,\tau)$ has a cofinite subspace on which the $\tau$-topology is coarser than the euclidean topology. This property clearly implies being cell-wise euclidean. The converse implication also holds in the case where $\dim X=1$ as follows. Suppose that $\mathcal{X}$ is a finite definable partition of $X$ into subsets where the subspace $\tau$-topology is euclidean. Consider the set $A=\bigcup\{ C \setminus cl_e(X\setminus C) : C\in \XX,\, \dim C = \dim X\}$. 
Since $\tau|_C = \tau_e|_C$ for every $C\in \XX$, observe that it must be that that $\tau|_A = \tau_e|_A$. Furthermore, since the euclidean topology satisfies the \textbf{fdi}, it holds that $\dim X\setminus A < \dim X$. In particular if $\dim X=1$ then $X\setminus A$ is finite, and so $(X,\tau)$ satisfies almost $\tau \subseteq \tau^{af}|_X$. Statement (\ref{itm:Pet_Rosel_2}) in Theorem~\ref{them:PR} can thus be reformulated as stating either that $(X,\tau)$ is almost $\tau \subseteq \tau^{af}|_X$ or that $(X,\tau)$ is cell-wise euclidean. 

We remark here that there appears to be an infelicity in the proof of the first assertion in~\cite[Lemma 3.16]{pet_rosel_18}, of which the authors of~\cite{pet_rosel_18} are aware (per private correspondence). Specifically, the proof relies on the assertion (in our terminology) that, for any Hausdorff one-dimensional definable topological space $(X,\tau)$ and $x\in X$, it holds that $\E_x \setminus \{x\} \subseteq \bigcap_{U\in \BB(x)} \partial_e U $, where $\BB(x)$ denotes a definable basis of $\tau$-neighbourhoods of $x$. However, to see that this statement is false in general suppose that $\E_x\setminus \{x\} \neq \emptyset$ and $X\in \BB(x)$. 
The authors of~\cite{pet_rosel_18} have communicated to us (by private correspondence) a correction to their proof, which we understand may be forthcoming. 
Moreover, observe that the first statement of~\cite[Lemma 3.16]{pet_rosel_18} in fact corresponds to our Lemma~\ref{lemma_2} localised to Hausdorff one-dimensional spaces. Therefore, this issue does not present any reason to suppose that the statement of Theorem~\ref{them:PR} does not hold.

Furthermore, it is straightforward to derive the statement of Theorem~\ref{them:PR}, in the case that $\RR$ expands an ordered field, from the results of this paper.
In this setting, in light of Theorem~\ref{them 1.5}, the implications $(\ref{itm:Pet_Rosel_1})\Leftrightarrow (\ref{itm:Pet_Rosel_2})\Leftrightarrow (\ref{itm:Pet_Rosel_3})$ in Theorem~\ref{them:PR} are equivalent to our Theorem~\ref{them_main_2} and Remark~\ref{rem:(3)implies(1)}. 
The implication $(\ref{itm:Pet_Rosel_1}) \Rightarrow (\ref{itm:Pet_Rosel_4})$ is a classical corollary of o-minimal cell decomposition (see~\cite[Chapter 3, Proposition 2.18]{dries98}) and the easy fact that the euclidean topology is regular. 
Finally, the implication $(\ref{itm:Pet_Rosel_4}) \Rightarrow (\ref{itm:Pet_Rosel_1})$ can be derived using the framework that we introduced to prove  Theorem~\ref{them_ADC_or_DOTS} (in particular Lemmas~\ref{remark_for_two_theorems} and~\ref{lemma_proof_two_thems}). We outline this later in the section (Proposition~\ref{prop:connectedness-cellwise-euclidean}).

Although our approach in this paper leads to an affineness characterization only in the case where $\RR$ expands an ordered field, whereas the equivalence in Theorem~\ref{them:PR} can be shown to hold when $\RR$ expands an ordered group as long as the underlying set $X$ is bounded, our alternative approach does allow us to give answers to a number of questions left open in~\cite{pet_rosel_18}. 

Firstly, Peterzil and Rosel ask (remark at the end of Section 2 in~\cite{pet_rosel_18}) if, given a definable topological space $(X,\tau)$ and $x\in X$, the union of all (definable) definably connected sets containing $x$ is itself definable, i.e. if there exists a (definable) definably connected component containing $x$. We answer this question in the positive in the case that $(X,\tau)$ is $T_3$ with $\dim(X) = 1$ in Proposition~\ref{prop:connected-components} below. In order to prove it we require two lemmas.

\begin{lemma}\label{lem:clopen-tau}
Let $(X,\tau)$ be a Hausdorff regular definable topological space, with $X\subseteq R$. Let $\Xinf$ be as defined by Lemma~\ref{remark_for_two_theorems}, and consider a set $A=\bigcup_{0\leq i < n} f_i(I)$ in $\Xinf$, with $I=(a,b)$. Let $(A^*,\tau_A)$ and $h_A:(A,\tau)\hookrightarrow (A^*,\tau_A)$ be as given by Lemma~\ref{lemma_proof_two_thems}, with $A^*=I\times\{0,\ldots,m\}$. Suppose that $m>0$ and $\tau_A=\taulex$. For any $a<c<d<b$, if we set
\[
A(c,d, m)=((c,d]\times \{0\}) \cup ([c,d)\times \{m\}) \cup \bigcup_{0 <i <m} (c,d)\times \{i\},
\]
then $h_A^{-1}(A(c,d,m))$ is clopen in $(X,\tau)$. 
\end{lemma}
\begin{proof}
Fix $a<c<d<b$. By definition of the $\taulex$ topology, the set $A(c,d,m)$ is clopen in $(A^*,\tau_A)=(A^*,\taulex)$ and so, since $h_A$ is continuous, the set $B=h_A^{-1}(A(c,d,m))$ is clopen in $(A,\tau)$. Since $A$ is $\tau$-open in $X$ it immediately follows that $B$ is $\tau$-open in $X$ too. It remains to show that it is $\tau$-closed in $X$.

Let $x\in cl_\tau B$. We show that $x\in B$. 
By Lemma~\ref{lemma_proof_two_thems}(\ref{itm1.5:lemma_proof_two_theorems}) note that 
\[
B\subseteq\bigcup_{0\leq i<n} f_i([c,d]). 
\]
In particular $cl_e B \subseteq A$. Since $x\in cl_\tau B$ then, by Proposition~\ref{basic_facts_P_x_1}(\ref{itm: basic_facts_4}), we have that $\emptyset\neq \E_x \cap cl_e B \subseteq A$. However, by Lemma~\ref{properties_cofnite_subset_1.5}(\ref{itm:properties_cofnite_subset_1.5_b}) and the fact that $A$ is closed under the equivalence relation addressed in said lemma (i.e. $A=\bigcup_{z\in I} [z]$ by Lemma~\ref{remark_for_two_theorems}), this implies that $x\in A$. Since $B$ is closed in $(A,\tau)$, we conclude that $x\in B$. 
\end{proof}

\begin{lemma}\label{lem:connected-T3}
Let $(X,\tau)$ be a Hausdorff regular definable topological space, with $X\subseteq R$. Let $\Xinf$ be as defined by Lemma~\ref{remark_for_two_theorems}. Suppose that there exists $x\in \bigcup \Xinf$ with $x\notin R_x \cap L_x$. Then, for every $y\in X\setminus \{x\}$, there exists a definable $\tau$-clopen set $B$ with $x\in B$ and $y \notin B$. In particular $\{x\}$ is a maximal definably connected subspace of $(X,\tau)$.
\end{lemma}
\begin{proof}
Let us fix $x\in \bigcup \Xinf$ with $x\notin R_x \cap L_x$ and $y\in X\setminus \{x\}$. Since otherwise the result is obvious we may assume that neither $x$ nor $y$ is $\tau$-isolated in $X$.  

Following the terminology in Lemma~\ref{remark_for_two_theorems}, let $A=\bigcup_{0\leq i < n} f_i(I)$ be the set in $\Xinf$ satisfying that $x\in A$, with $I=(a,b)$. Let $(A^*,\tau_A)$ and $h_A:(A,\tau)\hookrightarrow (A^*,\tau_A)$ be as given by Lemma~\ref{lemma_proof_two_thems}, with $A^*=I\times\{0,\ldots,m\}$. Now since $x\notin R_x \cap L_x$ and moreover $x$ is not $\tau$-isolated, note that this avoids Case 5 in the proof of Lemma~\ref{lemma_proof_two_thems}, and in all the other cases it holds that $m>0$ and $\tau_A=\taulex$. Consequently, by Lemma~\ref{lem:clopen-tau}, in order to prove the lemma it suffices to show that there exists some $a<c<d<b$ such that the set $B=h_A^{-1}(A(c,d,m))$ contains exactly one point among $\{x,y\}$. If $y\notin A$ then this is obvious, and so we assume that $y\in A$. 

Let $h_A(x)=\al x_0, i \ar$ and $h_A(y)=\al x_1, j\ar$, for points $x_0, x_1 \in I$ and $0\leq i, j \leq m$. 
If $x_0\neq x_1$, then suppose without loss of generality that $x_0 < x_1$. Then it suffices to choose any $a< c < x_0$ and any $x_0 < d < x_1$, and the result clearly follows. Now suppose that $x_0=x_1$. In particular, by Lemma~\ref{lemma_proof_two_thems}(\ref{itm1.5:lemma_proof_two_theorems}), note that $\{x,y\} \subseteq \{ f_i(x_0) : 0\leq i <n\}$. 

Since by assumption both $x$ and $y$ are not $\tau$-isolated then, by Lemma~\ref{lemma_explaining_neighbourhoods}, we have that $\E_x\neq \emptyset$ and $\E_y\neq \emptyset$. In particular, by Lemma~\ref{remark_for_two_theorems}, it holds that $[x_0]^E=\{ x_0, f_{n-1}(x_0)\}=\{x, y\}$. Observe that this corresponds in the proof of Lemma~\ref{lemma_proof_two_thems} to Cases 3 and 4. In both these cases it holds that $m=n-1$, and moreover $\{h_A(x), h_A(y)\}=\{h_A(x_0), h_A(f_{n-1}(x_0))\} = \{\al x_0, 0\ar, \al x_0, m\ar\}$. Finally, note that the set $A(x_0, b, m)$ contains the point $\al x_0, m \ar$ but does not contain the point $\al x_0, 0 \ar$, and so we conclude that the set $B=h_A^{-1}(A(x_0,b,m))$ contains exactly one point among $\{x,y\}$, as desired. 
\end{proof}

\begin{proposition}\label{prop:connected-components}
Let $(X,\tau)$ be a Hausdorff regular definable topological space, with $X\subseteq R$. For each $x\in X$ there exists a maximal definably connected definable set $C\subseteq X$ containing $x$. Furthermore, either $C=\{x\}$ or $(C,\tau)$ is an infinite cell-wise euclidean subspace. 
\end{proposition}
\begin{proof}
Let $\Xinf$ be as defined by Lemma~\ref{remark_for_two_theorems}.
Let $Z=\bigcup\Xsing \cup \{ x\in \bigcup \Xinf : x\in R_x \cap L_x\}$. Note that, by Lemma~\ref{lem:connected-T3}, any point in $X\setminus Z$ is not contained in any definably connected set in $(X,\tau)$ besides $\{x\}$. We show that $(Z,\tau)$ is cell-wise euclidean. The proposition follows. 
 
Since $\Xsing$ is finite, to prove that $(Z,\tau)$ is cell-wise euclidean it suffices to show that the set $\{ x\in \bigcup \Xinf : x\in R_x \cap L_x\}$ with the subspace $\tau$-topology is cell-wise euclidean. Let us fix $x$ in this set. Following the terminology in Lemma~\ref{remark_for_two_theorems}, let $A=\bigcup_{0\leq i <n} f_i(I) \in \Xinf$ be such that $x\in A$. We complete the proof by showing that $x \in I\subseteq \{ x\in \bigcup \Xinf : x\in R_x \cap L_x\}$ and $(I,\tau)=(I, \tau_e)$. Let $(A^*, \tau_A)$ and $h_A$ be as given by Lemma~\ref{lemma_proof_two_thems}. Since $x\in R_x \cap L_x$, note that, by Lemma~\ref{remark_for_two_theorems}, $x\in I\cup f_{n-1}(I)$, and moreover this corresponds to Case 5 in the proof of Lemma~\ref{lemma_proof_two_thems}. Observe that, in Case 5, every point in $f_i(I)$ for $0<i<n$ is $\tau$-isolated, and so we derive that $x\in I$.  Furthermore, every point $x'\in I$ satisfies that $x'\in R_{x'}\cap L_{x'}$, and so $I\subseteq  \{ x\in \bigcup \Xinf : x\in R_x \cap L_x\}$. Additionally, again by the fact that we are in Case 5, it holds that $\tau_A=\taualex$ and $h_A(f_i(x')) = \al x', i\ar$ for every $x'\in I$ and $0\leq i <n$. Since $(I\times\{0\}, \taualex)=(I\times \{0\}, \tau_e)$, and $h_A:(A,\tau)\hookrightarrow (A^*, \tau_A)$ is an embedding given by $x'\mapsto \al x', 0\ar$ for every $x'\in I$, we derive that $(I,\tau)=(I,\tau_e)$.
\end{proof}
Using Lemma~\ref{lem:connected-T3}, we may now also explain how implication $(\ref{itm:Pet_Rosel_4})\Rightarrow(\ref{itm:Pet_Rosel_1})$ in Theorem~\ref{them:PR} (in the case that $\RR$ expands an ordered field) can be obtained from our approach to these ideas. Specifically, the implication follows immediately from the proposition below and the proof of Theorem~\ref{them_main_2}, which shows that, if $\RR$ expands an ordered field, then every one-dimensional Hausdorff cell-wise euclidean definable topological space is affine. 

\begin{proposition}\label{prop:connectedness-cellwise-euclidean}
Let $(X,\tau)$ be a Hausdorff regular definable topological space, with $X\subseteq R$. If $(X,\tau)$ has finitely many definably connected components then it is cell-wise euclidean. 
\end{proposition}
\begin{proof}
Let $(X,\tau)$ be a $T_3$ definable topological space with $X\subseteq R$ and suppose that it has finitely many definably connected components. Let $\Xinf$ be as defined by Lemma~\ref{remark_for_two_theorems} and let $A=\bigcup_{0\leq i < n} f_i(I)$ be a set in $\Xinf$. We show that $(A,\tau)=(I,\tau_e)$. Since $\bigcup \Xinf$ is cofinite in $X$ it follows that $(X,\tau)$ is cell-wise euclidean.

Let $(A^*,\tau_A)$ and $h_A$ be as defined by Lemma~\ref{lemma_proof_two_thems}. Since $(X,\tau)$ has finitely many definably connected components, by Lemma~\ref{lem:connected-T3} it follows that there can only be finitely many points in $A$ satisfying that $x\notin R_x \cap L_x$. Observe that this rules out Cases 0, 1, 2, 3 and 4 in the proof of Lemma~\ref{lemma_proof_two_thems}. 
Furthermore observe that, in the only remaining case -- Case 5 -- it holds that every point $x\in f_i(I)$, for $0<i<n-1$, is $\tau$-isolated, and so, since $(X,\tau)$ has finitely many definably connected components (in particular finitely many $\tau$-isolated points), it must be that $n=1$, i.e. $A=f_0(I)=I$. It then follows directly from the fact that we are in Case 5 that $(A^*, \tau_A) =  (I \times \{ 0 \}, \taualex) = (I \times \{0\}, \tau_e)$. Moreover, $h_A$ is a homeomorphism $(I, \tau) \to (I \times \{0 \}, \tau_e)$ which, by Lemma~\ref{lemma_proof_two_thems}(\ref{itm1.5:lemma_proof_two_theorems}), is given by $x \mapsto \al x, 0 \ar$. So $(A, \tau) = (I, \tau) = (I, \tau_e)$, as desired.
\end{proof}

Peterzil and Rosel also raise the question of whether or not an analogue to their result, Theorem~\ref{them:PR} above, could be obtained for Hausdorff definable topological spaces in higher dimensions (\cite{pet_rosel_18}, Section 4.3 (2)).  
They note that being affine (condition (\ref{itm:Pet_Rosel_1})) cannot be equivalent to condition (\ref{itm:Pet_Rosel_2}), namely having a cofinite subspace whose topology is coarser than the euclidean topology, in arbitary dimensions, but they leave as open questions whether or not being affine is equivalent to condition (\ref{itm:Pet_Rosel_3}), namely that every definable subspace has finitely many definably connected components, or to condition (\ref{itm:Pet_Rosel_4}), namely being regular and having finitely many definably connected components. We can answer both of these questions negatively: the non-equivalence of conditions (\ref{itm:Pet_Rosel_1}), (\ref{itm:Pet_Rosel_3}) and (\ref{itm:Pet_Rosel_4}) in dimension greater than one is given by Examples~\ref{example_line-wise_euclidean_not_euclidean} and~\ref{example_space_cell-wise_euclidean_not_metrizable} in the Appendix.

\appendix

\section{Examples}\label{section:examples}

In this appendix we compile examples that witness the heterogeneity of definable topological spaces with reference to their (definable) topological properties, and help frame the results in this paper and their limitations when trying to improve or generalize them. 

The examples are given in the language $(0,1,+,-,\cdot, <)$, where $\RR$ is assumed to expand an ordered group \mbox{$(R,0,+,-,<)$} or an ordered field $(R,0,1,+,-,\cdot, <)$ whenever the corresponding function symbols are involved.   

Since we are working in the generality of an o-minimal structure, it is important to note that we will not address certain classical topological properties of definable topological spaces, because they are dependent on the specifics of the underlying structure $\RR$.  
These include compactness, connectedness, separability, normality or metrizability. We consider however definable versions of these properties (as per the definitions in this paper).

All the examples that are generalizations of classical topological spaces (e.g. definable Split Interval (Example~\ref{example:split_interval}), definable Alexandrov Double Circle (Example~\ref{example:alex_double_circle})) behave, in terms of their definable topological properties, exactly like their classical counterparts. The only exception to this is the Definable Moore Plane (Example~\ref{example:Moore_plane}), which is definably normal. 

We begin by recalling the key examples of definable topological spaces that are used and studied extensively throughout the paper.

\begin{example}[Euclidean topology]\label{example: euclidean topology}
The euclidean topology $\tau_e$ on $R^n$ has definable basis 
\[
\left\{ \prod_{1\leq i \leq n}(x_i,y_i) : x_i<y_i,\, 1\leq i \leq n\right\}.
\]
It is $T_3$, definably separable (Proposition~\ref{prop:def-sep_eucl_disc_llt}(\ref{itm1:sep-spaces-examples})), definably connected and definably metrizable. It is moreover definably compact if and only if it is restricted to a closed and bounded set~\cite{pet_stein_99}.   
\end{example}

\begin{example}[Discrete topology]
The discrete topology $\tauc$ on $R^n$ has definable basis
\[
\{\{x\} : x\in R^n\}.
\]
Note that this topology is definable on any definable set in any model-theoretic structure. It is  $T_3$ and definably metrizable.
\end{example}

\begin{example}[Half-open interval topologies]\label{example:taur_taul}
The right half-open interval topology (or lower limit topology) $\tau_r$ has definable basis
\[
\{ [x,y) : x, y\in R,\, x<y\}.
\]
The space $(\mathbb{R},\tau_r)$ is classically called the Sorgenfrey Line.

The left half-open interval topology (or upper limit topology) $\tau_l$ has definable basis 
\[
\{ (x,y] : x, y\in R,\, x<y\}.
\]
These topologies are $T_3$ and definably separable (Proposition~\ref{prop:def-sep_eucl_disc_llt}(\ref{itm3:sep-spaces-examples})). They are also totally definably disconnected (the only definably connected subspaces are singletons) and not definably metrizable (see Proposition~\ref{prop:tau_r not metrizable}). 
\end{example}

\begin{example}[Definable Alexandrov Double Arrow space or definable Split Interval]\label{example:split_interval}
Let $X=[0,1]\times \{0,1\}$. The definable Alexandrov Double Arrow space (or definable Split Interval) is the space $(X,\taulex)$, where $\taulex$ denotes the topology induced by the lexicographic order on $X$. The space is classically called the Alexandrov Double Arrow space (or Split Interval) when $\RR$ expands $(\R,<)$.  
It is $T_3$, definably compact, definably separable and totally definably disconnected.

It is not definably metrizable since the bottom line $[0,1]\times \{0\}$ is definably homeomorphic to $([0,1], \tau_l)$ and the top line $[0,1]\times \{1\}$ to $([0,1],\tau_r)$. It is also worth noting that $(X,\tau)$ does not satisfy the \textbf{fdi}, since $\partial([0,1]\times\{0\})=[0,1)\times\{1\}$. Moreover, one may show that $[0,1]\times \{0\}$ is not a boolean combination of open definable sets, which was a tameness condition for definable topologies considered by Pillay in~\cite{pillay87}.  
\end{example}

The following two examples, the definable $n$-split interval and definable Alexandrov $n$-line, were already introduced in Definitions~\ref{dfn:lex} and~\ref{example_n_line_0}, and play a crucial role in Theorem~\ref{them_ADC_or_DOTS}. They were motivated by the classical Split Interval (see Example~\ref{example:split_interval} above) and Alexandrov Double Circle (see Example~\ref{example:alex_double_circle} below) respectively.

\begin{example}[Definable $n$-split interval] \label{example:n-split}
Let $n>0$. We call the space $(R\times\{0,\ldots, n-1\}, \taulex)$ the definable $n$-split interval. 

If $n=1$, then note that $(R\times\{0\},\taulex$)$=(R\times\{0\}, \tau_e)$. If $n>1$, then,
by analogy to the definable Split Interval (Example~\ref{example:split_interval}), the definable $n$-split interval is $T_3$, definably separable, totally definably disconnected, not definably metrizable and it does not have the \textbf{fdi}. Furthermore, it is not definably compact but every subspace of the form $I\times \{0,\ldots, n-1\}$, where $I\subseteq R$ is a closed and bounded interval, is definably compact. 

Moreover, for any $0<n<m<\omega$, the spaces $(R\times\{0,\ldots, n-1\}, \taulex)$ and $(R\times\{0,\ldots, m-1\}, \taulex)$ are not definably homeomorphic. In fact, they are not even in definable bijection, since they have different Euler characteristic (see~\cite[Chapter 4]{dries98}). Specifically, observe that every finite cell partition of $R\times\{0,\ldots, n-1\}$ will contain $n$ more cells of dimension one than points. And similarly every cell partition of $R\times\{0,\ldots, m-1\}$ will contain $m$ more cells of dimension one than points. On the other hand, by o-minimal cell decomposition, every definable injection from $R\times\{0,\ldots, n-1\}$ to $R\times\{0,\ldots, m-1\}$  can be decomposed into disjoint definable bijections between singletons and one-dimensional cells, and so any such injection is not surjective. 
\end{example}

\begin{example}[Definable Alexandrov $n$-line] \label{example_n_line}
For any $y<x<z$ in $R$, let 
\[
A(x,y,z)=\{\al x,0 \ar \} \cup (((y,z)\setminus\{x\}) \times R).
\] 
Let $\taualex$ be the topology on $R^2$ with definable basis \[
\{ A(x,y,z) : y<x<z\} \cup \{ \{\al x,y \ar \} : y \neq 0\}.
\]  
Let $n>0$. The definable Alexandrov $n$-line is the definable topological space $(R\times\{0,\ldots, n-1\},\taualex)$. In the case that $n=1$, this is a definable analogue of the classical Alexandrov double of the space $(R,\tau_e)$ (see~\cite{engelking68}).

This space is $T_3$. If $n=1$, then it is simply the euclidean topology on $R\times\{0\}$. Suppose that $n>1$. The subset $R\times \{i\}$ for any $i>0$ contains only isolated points, and so the space is not definably separable. For any closed and bounded interval $I\subseteq R$, the subspace $I\times\{0,\ldots,n-1\}$ is definably compact but not definably separable. It follows that it is not definably metrizable (see~\cite[Lemma 7.4]{walsberg15}, which states that any definably compact definable metric space is definably separable). 
Hence the definable Alexandrov $n$-line, for $n>1$, is not definably metrizable. 

For any $0<n<m<\omega$, the sets $R\times\{0,\ldots, n-1\}$ and $R\times\{0,\ldots, n-1\}$ are not in definable bijection (see Example~\ref{example:n-split}), and so in particular the spaces $(R\times\{0,\ldots, n-1\}, \taualex)$ and $(R\times\{0,\ldots, m-1\}, \taualex)$ are not definably homeomorphic.
\end{example}

The next five examples are provided in order to illustrate the necessity of certain hypotheses in some of the key results of this paper, in particular Proposition~\ref{prop_T2_frontier_ineq_regular}, Proposition~\ref{them_main}, Corollary~\ref{cor:general_them_main}, Theorem~\ref{them 1.5}, Corollary~\ref{cor 1.5} and Theorem~\ref{them_main_2}.

\begin{example}\label{example: t_0 not t_1}
Consider the following definable basis for a topology on $R$.
\[
\{ (-\infty,x] : x\in R \}.
\]
The resulting space is $T_0$ but not $T_1$. Any subspace with more than one element fails to be $T_1$;  
in particular, no interval subspace of this space has the euclidean, discrete or half-open interval topologies. Consequently, the $T_1$ assumption in Proposition~\ref{them_main} and Corollary~\ref{cor:general_them_main} (the definable version of the Gruenhage 3-element-basis Conjecture) cannot be weakened to $T_0$. 
\end{example}

\begin{example}\label{example: t_1 not t_2}
Consider the definable family of sets 
\[
\{ (-\infty,x) \cup (y,z) : x<y<z\}.
\] 
It is a basis for a topology $\tau$ on $R$ that is $T_1$ but not Hausdorff. 

Any finite definable partition of $R$ must include an interval of the form $(-\infty, x)$, whose subspace topology is not Hausdorff. In particular, $(R,\tau)$ cannot be decomposed into finitely many definable subspaces with the euclidean, discrete or half-open interval topologies. It follows that the Hausdorffness assumption in Theorem~\ref{them 1.5} (and \sout{hence} in Corollary~\ref{cor 1.5}) cannot be weakened to $T_1$-ness. 
\end{example}

\begin{example}\label{example: T_1 fdi not regular}
Let $X=[0,1) \cup \{2\}$ and consider a topology $\tau$ on $X$ such that the subspace topology $\tau|_{[0,1)}$ is euclidean and a basis of open neighbourhoods of $\{2\}$ is given by 
\[
\{ (0,x) \cup \{2\} : 0<x<1\}.
\]
This topology is clearly definable and $T_1$ but not Hausdorff, since points $0$ and $2$ fail to have disjoint neighbourhoods. In particular, it is not regular. It is easy to observe that it satisfies the \textbf{fdi}. 
Since it fails to be regular, it illustrates the necessity of the Hausdorffness assumption in Proposition~\ref{prop_T2_frontier_ineq_regular}. Moreover, we note, by considering the partition into subspaces $(0,1)$ and $\{2\}$, that this space is cell-wise euclidean. So it is also not true that every cell-wise euclidean one-dimensional space is Hausdorff, and in particular affine. 
\end{example}

\begin{example}\label{example: dfbly sep not hereditary}
Let $X=R \times\{0,1\}$ and consider a topology $\tau$ on $X$ given by the basis 
\[
\{\{\al x,0\ar\}\cup ((x,y)\times \{1\}) : x<y\} \cup  \{(z,x]\times \{1\} : z<x\}.
\] 
This space is Hausdorff but not regular, since $R\times \{0\}$ is a closed set and, for any $x\in R$ and any neighbourhood $U=(z,x]\times \{1\}$ of $\al x,1\ar$, we have $cl_\tau (U) \cap (R\times\{0\})=[z,x)\times\{0\} \neq \emptyset$. Moreover, note that, since $\partial_\tau (R\times \{1\})=R\times\{0\}$, this space does not satsify the \textbf{fdi}. This example therefore illustrates the necessity of the assumption of satisfying the \textbf{fdi} in Proposition~\ref{prop_T2_frontier_ineq_regular}. 

This example also illustrates a key fact about our notion of definable separability (Definition~\ref{dfn:separable}). Recall that any definable subspace of a definably separable definable metric space is also definably separable (see Lemma~\ref{remark_def_sep}). The space $(X,\tau)$ in the current example is definably separable, but the subspace $R\times \{0\}$ is infinite and discrete, hence $(X,\tau)$ is not definably metrizable, and, moreover, we see that definable separability, much like separability, is not in general a hereditary property. 
\end{example}

\begin{example}\label{example:fdi_Hausdorff_not_regular}
Let $X=\{ \al 0,0 \ar \} \cup [0,1)\times (0,1)$. Consider the topological space $(X,\tau)$, where the subspace $X\setminus \{\al 0,0 \ar \}$ is euclidean, and a basis of open neighbourhoods for $\al 0,0 \ar$ is given by sets 
\[
A(t)=\{\al 0,0\ar \} \cup \left( (0,1)\times (0,t) \right),
\]
for $0<t<1$. 
The topology $\tau$ is clearly definable and Hausdorff. Moreover, for any $0<t<1$, the $\tau$-closure of $A(t)$ is $\{\al 0,0 \ar \} \cup [0,1) \times (0, t]$, and so the space is not regular, since the point $\al 0,0\ar$ and the closed set $\{0\}\times (0,1)$ are not separated by neighbourhoods. 

Since $(X,\tau)$ is $T_1$ and the subspace $X\setminus \{\al 0,0\ar\}$ is euclidean, it easily follows that $(X,\tau)$ satisfies the \textbf{fdi}. Hence $(X,\tau)$ is Hausdorff and satisfies the \textbf{fdi} but fails to be regular, thus is a counterexample to the generalization of Proposition~\ref{prop_T2_frontier_ineq_regular} to spaces of dimension greater than one. 

Moreover, the space can be partitioned into two euclidean subspaces, namely $\{\al 0,0 \ar\}$ and $X\setminus \{\al 0,0\ar\}$; in particular it contains no definable copy of an interval with either the discrete or the right half-open interval topology. However, it is not metrizable, since it is not regular. Hence it is a counterexample to a generalization of Theorem~\ref{them_main_2} to spaces of dimension two. 
\end{example}

For the remaining examples, let $B_2(\al x,y \ar,t)$, for $\al x,y\ar\in R^2$ and $t>0$, denote the ball in the $2$-norm of center $\al x, y\ar$ and radius $t$, namely
\[
B_2(\al x,y\ar,t)=\{\al x',y'\ar \in R^2 : (x-x')^2+(y-y')^2<t^2\}. 
\]
\begin{example}[Definable Moore Plane]\label{example:Moore_plane}
Let $X=\{\al x,y \ar \in R^2: y\geq 0\}$ be the closed upper half-plane.
Let $\BB_e$ be a definable basis for the euclidean topology in $\{\al x,y \ar \in R^2 : y>0\}$ and, for any $x\in R$ and $\epsilon >0$, let 
\[
A(x,\varepsilon)=B_2(\al x,\varepsilon \ar,\varepsilon)\cup\{\al x,0 \ar\}.
\]
The family $\BB=\BB_e \cup \{ A(x,t) : x\in R, t>0\}$ is clearly definable and forms a basis for a topology $\tau$. We call the space $(X,\tau)$ the definable Moore Plane.

This space is $T_3$ and definably separable but not definably metrizable since the subspace $R\times \{0\}$ is infinite and discrete (see Lemma~\ref{remark_def_sep}). When $\RR$ expands the field of reals, the Moore Plane is a classical example of a separable non-normal space (and in particular, it is not metrizable). 
\end{example}

It is worth noting that, even though our definition of definable normality (Definition~\ref{definition:dfbly_normal}) seems the natural adaptation of the classical notion, the classical Moore Plane fails to be normal, but one may show that the definable Moore Plane is definably normal. This suggests that our notion of definable normality might not be adequate. 
Moreover, Fornasiero also considered this same notion of definable normality in unpublished work \cite{fornasiero} (seen in private correspondence) where he showed that a definable topological space that is definably compact (in the sense of condition \eqref{dfn:directed-compact}) and Hausdorff is not necessarily definably normal, in contrast to the classical fact that a compact Hausdorff space is normal. (However, he did also show that if a definably compact, Hausdorff space is given by a definable uniformity (see~\cite{simon_walsberg19} for definitions), then it is indeed definably normal in this sense.)

\begin{example}[Definable Alexandrov Double Circle]\label{example:alex_double_circle}
Let $X=C_1 \cup C_2$, where $C_1$ and $C_2$ denote respectively the unit circle and circle of radius two in $R^2$ centered at the origin. Let $f:C_1\rightarrow C_2$ be the natural $e$-homeomorphism given by $x\mapsto 2x$.  Let 
\[
\BB_1=\{ \left( B_2(x,t)\cap C_1 \right) \cup f(B_2(x,t)\cap C_1 \setminus\{x\}) : x\in C_1, t>0\}
\]
and $\BB_2=\{ \{x\} : x\in C_2\}$. The definable Alexandrov Double Circle is the topology on $X$ generated by the basis $\BB_1\cup \BB_2$.

This space is definably compact and Hausdorff, but not definably separable, since $C_2$ is an infinite definable set of isolated points. It follows (see~\cite[Lemma 7.4]{walsberg15}) that it is not definably metrizable. It also fails to satisfy the \textbf{fdi}, since the outer circle $C_2$ is a dense subset. 

When $\RR$ expands the field of reals this space is simply called the Alexandrov Double Circle and is a classical example of a compact non-separable space (hence one that is not metrizable). 
\end{example}

The following example shows that there exists a Hausdorff two-dimensional definable topological space that does not contain a definable copy of an interval with either the discrete or the lower limit topology but still fails to be cell-wise euclidean. This shows that Remark~\ref{remark_cell-wise_euclidean} cannot be generalized to higher dimensions. In particular, this example is not affine but, by Theorem~\ref{them_main_2}, any one-dimensional subspace is affine (i.e. it is \textquotedblleft line-wise" affine). This is proved below in Proposition~\ref{prop:line-wise_euclidean_not_euclidean}.

\begin{example}[The definable hollow plane]\label{example_line-wise_euclidean_not_euclidean}  

We construct a basis for a topology $\mytau$ on $R^2$ by considering, for each point $x$, a basis of open neighbourhoods given by open euclidean balls without the graph of 
$f(t)=t^2$, for $0<t$, translated to have its origin at $x$. 
That is, for a given $x=\al x_1,x_2 \ar \in R^2$, let $\Gamma_x:=\{\al x_1+t, x_2+t^2 \ar : t>0\}$. Now let $\BB$ be given by sets
\[
A(x,t) =B_2(x,t)\setminus \Gamma_x,
\]
for $x\in R^2$ and $t>0$. We call $A(x,t)$ a $\mytau$-ball of center $x$ and radius $t$. 

We claim that $\mathcal{B}$ is a basis for a topology on $R^2$.
In order to prove this, let $A_0$ and $A_1$ be intersecting sets in $\mathcal{B}$ and let $x\in A_0 \cap A_1$. We show that there exists some $\varepsilon>0$ such that $A=A(x,\varepsilon)$ satisfies that $A(x,\varepsilon)\subseteq A_0 \cap A_1$. 

For any $y\in R^2$ and $t>0$, let $A^*(y,t)=A(y,t)\setminus \{y\}$. Note that, for any $A\in \BB$, the set $A^*$ is $e$-open. 
\begin{description}
\item[Case 1: $x\in  A^*_0 \cap A^*_1$.] Since $A^*_0\cap A^*_1$ is $e$-open, there is some $\varepsilon>0$ such that $B_2(x,\varepsilon)\subseteq A^*_0\cap A^*_1\subseteq A_0\cap A_1$. Hence we may take $A=A(x,\varepsilon)\subseteq B_2(x,\varepsilon)$.
\item[Case 2: $x\notin A^*_0 \cap A^*_1$.] Without loss of generality, suppose that $A_0=A(x,\varepsilon_0)$, for some $\varepsilon_0>0$. If $A_1=A(x,\varepsilon_1)$, for some $\varepsilon_1>0$, let $\varepsilon=\min\{\varepsilon_0, \varepsilon_1\}$ and $A=A(x,\varepsilon)$. Otherwise, by analogy to Case 1, let $\varepsilon_2>0$ be such that $A(x,\varepsilon_2)\subseteq A^*_1$ and let $A=A(x,\varepsilon)$, where $\varepsilon =\min\{\varepsilon_0,\varepsilon_2\}$.
\end{description}
So we may conclude that $\BB$ is a topological basis. Let $\mytau$ be the corresponding topology. We call $(R^2,\mytau)$ the definable hollow plane.

Every $e$-open set in $R^2$ is also $\mytau$-open, i.e. $\tau_e\subsetneq\mytau$. In particular $(R^2,\mytau)$ is Hausdorff. It fails, however, to be regular, since it is easy to check that, for any $x\in R^2$ and $\varepsilon>0$, $cl_{\mytau} A(x,\varepsilon) =cl_e B_2(x,\varepsilon)$, and so, for every $\mytau$-neighbourhood $A$ of $x$, $cl_{\mytau}A \cap \Gamma_x \neq \emptyset$. This space, however, is definably separable, which follows from (\ref{itm:broken_disk_1}) in the following proposition. One may also show that it is definably connected.

\begin{proposition}\label{prop:line-wise_euclidean_not_euclidean}
The following are properties of the definable hollow plane $(R^2,\mytau)$. 
\begin{enumerate}[(1)]
\item \label{itm:broken_disk_1} Any one-dimensional subspace of $(R^2,\mytau)$ is affine. 

\item \label{itm:broken_disk_2} No two-dimensional subspace of $(R^2,\mytau)$ is cell-wise euclidean. In particular, no two-dimensional subspace of $(R^2,\mytau)$ is affine.
\end{enumerate}
\end{proposition}
\begin{proof}
Statement~(\ref{itm:broken_disk_2}) is obvious from the definition. We prove~(\ref{itm:broken_disk_1}). By Theorem~\ref{them_main_2}, it suffices to show that $(R^2,\mytau)$ contains no subspace definably homeomorphic to an interval with either the discrete or the right half-open interval topology. 

Towards a contradiction, let $I\subseteq R$ be an interval and let $f:(I, \mu) \hookrightarrow (R^2, \mytau)$ be a definable embedding, where $\mu\in\{\tau_r,\tauc\}$. By o-minimality, after restricting $f$ if necessary, we may assume that $f$ is an $e$-embedding too. 

Since $\mu=\tau_r$ or $\mu=\tauc$ and $f$ is an embedding, it follows that, for any $t\in I$, $\mytaulim_{s\rightarrow t^-} f(s)\neq f(t)$ (see Proposition~\ref{prop_cont_lim}). So, by o-minimality, for every $t\in I$, there exists some $\varepsilon_t>0$ and $s'<t$ in $I$ such that, for all $s'<s<t$, $f(s)\notin A(f(t),\varepsilon_t)$. However, since $f$ is an $e$-embedding, there is also some $s'<s''<t$ such that, for all $s''<s<t$, $f(s)\in B_2(f(t),\varepsilon_t)$, and so $f[(s'',t)]\subseteq \Gamma_{f(t)}$. For any $t\in I$, let $s_t=\inf \{ s\in I : \,s<t,\, f[(s,t)]\subseteq \Gamma_{f(t)}\}$. This family is definable uniformly in $t\in I$. Since $s_t < t$ for all $t \in I$, by o-minimality, passing to a subinterval on which $t \mapsto s_t$ is continuous, if necessary, there exists an interval $J\subseteq I$ such that, for every $t\in J$, $s_t < J$. In other words, for every $s<t$ in $J$, it holds that $f(s)\in \Gamma_{f(t)}$. 

We now claim that, for any two distinct points $y,z\in R^2$, $|\Gamma_y \cap \Gamma_z|=1$. In that case, we have a contradiction, since we have shown that, for any $s, s', t, t'\in J$, if $s<s'<t<t'$, then $\{f(s), f(s')\} \subseteq \Gamma_{f(t)} \cap \Gamma_{f(t')}$. It therefore remains to prove the claim.  

Let $y= \al y_1,y_2 \ar \in R^2$ and $z=\al z_1,z_2 \ar \in R^2$, with $y \neq z$. Suppose that there exist $t, s >0$ such that 
\[
\al y_1+t, y_2+t^2 \ar= \al z_1+s, z_2+ s^2 \ar.
\]
If $y_1=z_1$, then we would have $t=s$ and hence $y=z$, so in fact we must have $y_1 \neq z_1$. We then substitute $s=t+y_1-z_1$ into $t^2=z_2-y_2+s^2$ in order to get
\[
t^2=z_2-y_2+t^2 +2t(y_1-z_1)+(y_1-z_1)^2,
\]
and hence
\[
t=\frac{y_2-z_2-(y_1-z_1)^2}{2(y_1-z_1)}, \qquad s=\frac{y_2-z_2+(y_1-z_1)^2}{2(y_1-z_1)},
\]
which gives us the unique point in $\Gamma_y \cap \Gamma_z$, which proves the claim. 
\end{proof}
\end{example}

In light of Example~\ref{example_line-wise_euclidean_not_euclidean}, a natural question to ask is whether or not we may instead obtain an analogue to Theorem~\ref{them_main_2} for spaces of all dimensions by substituting the condition of having a definable copy of an interval with either the $\tau_r$ or the $\tauc$ topology (condition~(\ref{itm:them_main_2_1}) in Theorem~\ref{them_main_2}) for simply not being cell-wise euclidean. The answer, even if adding the additional assumption that the space be regular, is no, as witnessed by the following, our final example. 

\begin{example}[Space that is $T_3$ and cell-wise euclidean but not definably metrizable]\label{example_space_cell-wise_euclidean_not_metrizable}

Let $X=\{\langle x,y \rangle \in R^2 : y\geq 0\}$ be the closed upper half-plane. Let $\BB_e$ be a definable basis for the euclidean topology in $\{\langle x, y \rangle\in R^2 : y>0\}$. For any $x\in R$ and $t>0$, let 
\[
A(x,t) = \{\langle x,0 \rangle \}\cup \{ \langle x', y\rangle\in R^2 : |x'-x|<t,\, 0\leq y<t|x'-x| \}. 
\]
Note that, for every $x\in R$ and $t>0$, there is $t'>0$ such that $A(x,t')\subseteq B_2(\al x, 0 \ar,t)$, while the converse is not true (we have that $B_2(\al x, 0 \ar,t') \nsubseteq A(x,t)$ for every $x\in R$ and $t,t'>0$). 
Moreover, for every $x\in R$, the family $\{A(x,t) : t>0\}$ is nested and, for every $t>0$, the set $A(x,t)\setminus\{\langle x,0 \rangle \}$ is $e$-open in $X$. From these three facts it follows, in a manner similar to the case analysis in Example~\ref{example_line-wise_euclidean_not_euclidean} , that the definable family $\BB_{\mytauu}= \BB_e \cup \{A(x,t) : x\in R, t>0\}$ is a basis for a topology $\mytauu$ on $X$.  

Since we have $\tau_e|_X \subsetneq \mytauu$, the topology $\mytauu$ is Hausdorff. Note that, for every $x\in R$ and $t>0$, we have $cl_{\mytauu} A(x,t)=cl_e A(x,t)$, and so $(X,\mytauu)$ is also regular. Moreover, the disjoint subspaces $\{\langle x,y \rangle : x\in R, y>0\}$ and $\{\langle x,y \rangle: x\in R, y=0\}$ are both euclidean, i.e. the space is cell-wise euclidean. In particular, the space is definably separable. Finally, it is also definably connected. 

When $\RR$ expands the field of reals this space is separable but not second countable and thus not metrizable. From the completeness of the theory of real closed fields it follows that there is no metric on $X$ definable in the language of ordered rings that induces $\mytauu$. We show that this holds in greater generality. 

\begin{proposition}
The space $(X,\mytauu)$ is not definably metrizable. 
\end{proposition}
\begin{proof}
Towards a contradiction, suppose that $(X,\mytauu)$ is definably metrizable with definable metric $d$. 
For every $x\in R$, let 
\[
r_x=\sup\{0<t<1 : B_d(\langle x,0 \rangle,t)\cap (\{x\}\times (0,\infty))=\emptyset\}.
\]
Note that, by definition of the neighbourhoods $A(x,t)$, we have necessarily that $r_x>0$, for every $x\in R$. By o-minimality, there exists an interval $I\subseteq R$ and some $r>0$ such that, for every $x\in I$, we have $r\leq r_x$. Now fix $x\in I$ and consider the $d$-ball $B_d(\al x,0\ar,r/2)$. By definition of $\mytauu$, there exists  some $y\in I\setminus \{x\}$ and some $s>0$ such that $\{y\}\times [0,s] \subseteq B_d(\langle x,0 \rangle,r/2)$. But then, by the triangle inequality, $d(\langle y,0 \rangle, \langle y,s \rangle)\leq d(\langle y,0 \rangle, \langle x,0 \rangle)+d(\langle x,0 \rangle, \langle y,s\rangle)<r \leq r_y$. This is a contradiction since, for every $0<t<r_y$, $B_d(\langle y,0\rangle,t)\cap (\{y\}\times (0,\infty))=\emptyset$.
\end{proof}
\end{example}

% \bibliography{mybib_topologies_line}
% \bibliographystyle{acm}

\providecommand{\noopsort}[1]{}\newcommand{\SortNoop}[1]{}

\end{document}